\documentclass{amsart}
\usepackage{graphics}
\usepackage{amssymb}
\usepackage[all]{xy}
\usepackage{color}
\CompileMatrices
\newtheorem*{introthm}{Theorem}
\newtheorem*{introprop}{Proposition}

\newtheorem{theorem}{Theorem}[section]
\newtheorem{lemma}[theorem]{Lemma}
\newtheorem{proposition}[theorem]{Proposition}
\newtheorem{corollary}[theorem]{Corollary}

\theoremstyle{definition}

\newtheorem{definition}[theorem]{Definition}
\newtheorem{example}[theorem]{Example}
\newtheorem{remark}[theorem]{Remark}

\numberwithin{equation}{theorem}
\def\Eff{{\rm Eff}}
\def\Mov{{\rm Mov}}
\def\Mori{{\rm Mori}}
\def\SAmple{{\rm Sample}}
\def\Ample{{\rm Ample}}

\def\div{{\rm div}}

\def\quot{/\!\!/}
\def\mal{\! \cdot \!}

\def\reg{{\rm reg}}

\def\rq#1{\widehat{#1}}
\def\t#1{\widetilde{#1}}
\def\b#1{\overline{#1}}
\def\bangle#1{\langle #1 \rangle}

\def\CC{{\mathbb C}}
\def\KK{{\mathbb K}}
\def\TT{{\mathbb T}}
\def\ZZ{{\mathbb Z}}

\def\NN{{\mathbb N}}
\def\QQ{{\mathbb Q}}
\def\PP{{\mathbb P}}
\def\rel{\operatorname{rlv}}
\def\CDiv{\operatorname{CDiv}}
\def\WDiv{\operatorname{WDiv}}

\def\Cl{\operatorname{Cl}}

\def\Pic{\operatorname{Pic}}

\def\Hom{{\rm Hom}}

\def\Supp{{\rm Supp}}
\def\Spec{{\rm Spec}}

\def\Sing{{\rm Sing}}

\def\cone{{\rm cone}}

\def\lin{{\rm lin}}
\def\cov{{\rm cov}}

\def\rank{\operatorname{rank}}

\def\topto#1{\stackrel{{\scriptscriptstyle #1}}{\longrightarrow}}
\newcounter{itemnumber}

\begin{document}
%

\title[Cox rings and combinatorics]
      {Cox rings and combinatorics}
\author[F. Berchtold]{Florian Berchtold} 
\address{Mathematisches Institut, Universit\"at Heidelberg,
69221 Heidelberg, Germany}
\email{Florian.Berchtold@uni-konstanz.de}
\author[J.~Hausen]{J\"urgen Hausen} 
\address{Mathematisches Forschungsinstitut Oberwolfach, Lorenzenhof, 
77709 Oberwolfach--Walke, Germany}
\email{hausen@mfo.de}
\begin{abstract}
Given a variety $X$ with a finitely generated total 
coordinate ring,
we describe basic geometric properties of $X$ in terms
of certain combinatorial structures living in the 
divisor class group of $X$.
For example, we describe the singularities, 
we calculate the ample cone, and we give simple Fano criteria.
As we show by means of several examples, 
the results allow explicit computations.
As immediate applications we obtain an effective version of
the Kleiman-Chevalley quasiprojectivity criterion, 
and the following observation on surfaces:
a~normal complete surface with finitely generated total 
coordinate ring is projective 
if and only if any two of its non-factorial singularities
admit a common affine neighbourhood.
\end{abstract}

\maketitle

\section*{Introduction}

The remarkable success of Cox's homogeneous coordinate 
ring in toric geometry~\cite{Co}
has stimulated many authors to study related
rings for more general varieties,
see~\cite{BaPo}, \cite{BeHa1}, \cite{ElKuWa}, 
\cite{Tsch} and~\cite{KeHu}.
We consider here the version discussed 
in~\cite{ElKuWa}:
Let $X$ be a normal variety with only constant
invertible global functions
and finitely generated free divisor class group 
$\Cl(X)$. Choose a subgroup $K$ 
of the group of Weil divisors mapping isomorphically
onto $\Cl(X)$. 
Then the modules of global sections fit together
to the {\em total coordinate ring}:
$$
\mathcal{R}(X) := \bigoplus_{D \in K} \Gamma(X,\mathcal{O}(D)).
$$

For a toric variety $X$, the total coordinate 
ring is a polynomial ring~\cite{Co}. 
In general, it turns out to be after all 
a factorial ring, see~\cite{BeHa1} and, 
more generally, \cite{ElKuWa}.  
The present paper is devoted to the study of
varieties $X$ with a finitely generated total 
coordinate ring.
For example, finite generation of 
$\mathcal{R}(X)$ is known to hold for
smooth Fano varieties of dimension at most 
three~\cite{KeHu}, 
and for smooth, unirational varieties with 
a complexity one group action~\cite{Kn}.

In~\cite{KeHu}, Keel and Hu characterized,
for $\QQ$-factorial projective varieties $X$,
finite generation of $\mathcal{R}(X)$ in terms 
of Mori Theory.
Moreover, they observed a 
connection to GIT: if $X$ is $\QQ$-factorial 
projective with $\mathcal{R}(X)$ finitely 
generated, then $X$ is the geometric quotient 
of a set $\rq{X} \subset \b{X}$ of semistable 
points of the affine variety 
$\b{X} = \Spec(\mathcal{R}(X))$
with respect to the natural action of the 
Neron Severi torus.
The small modifications of $X$ are reflected 
in the associated GIT chamber structure, 
and the various cones of divisors arising 
in this context, like the effective cone, 
the closure of the ample cone, and the Mori 
cone are all polyhedral.

In the present paper, we extend the picture 
developed in~\cite{KeHu} to certain normal, 
not necessarily (quasi)projective varieties $X$.
The main aim then is to provide a framework 
that allows explicit descriptions 
of geometric properties of $X$ in terms
of combinatorial data living in the divisor
class group $\Cl(X)$.
The approach generalizes the description 
of toric varieties in terms of
``bunches of cones'' presented in~\cite{BeHa2}.
As we shall see, it allows a systematic 
treatment of many non--toric examples, 
and thus enlarges the ``testing ground of algebraic 
geometry'' provided  by the toric varieties.
For example, many torus quotients of 
Grassmannians fit into our framework.

We now present the results a little more
in detail.
Similarly to~\cite{BeHa2},
we start with a combinatorial input datum, 
and then obtain a variety by means of 
a certain multigraded homogeneous spectrum
construction.  
The input is what we call a {\em bunched ring\/}
$(R,\mathfrak{F},\Phi)$; the precise 
definition is given in Section~\ref{bunchedringsI}.
For the moment, it suffices to know
that $R$ is a finitely generated factorial 
algebra, graded by a lattice $K$, 
that $\mathfrak{F}$ is a system of 
$K$-homogeneous prime generators $f_{1}, \ldots, f_{r}$ of $R$,
and that $\Phi$ is a certain finite collection 
of pairwise overlapping polyhedral cones in 
$K \otimes_{\ZZ} \QQ$,
each of which is generated by the degrees
of some of the $f_{i}$.

The basic construction then
associates to the bunched ring 
$(R,\mathfrak{F},\Phi)$ a 
normal variety $X = X(R,\mathfrak{F},\Phi)$ 
of dimension $\dim(R) - \rank(K)$
having divisor class group $K$ and total 
coordinate ring $R$,
see Section~\ref{bunchedringsI}. 
Geometrically, $X$ is a good quotient 
space of an open subset $\rq{X} \subset \b{X}$
of $\b{X} := \Spec(\mathcal{R}(X))$
by the torus action arising from the 
$K$-grading. The variety $X$ admits 
embeddings into toric varieties, 
and thus inherits the $A_{2}$-property, i.e.,
any pair $x,x' \in X$ has a common affine 
neighbourhood, compare~\cite{Wl}.
In fact $X$ is even {\em $A_2$-maximal\/} 
in the following sense: for any open
embedding $X \subset X'$ into 
an $A_2$-variety $X'$ such that  
$X' \setminus X$ is of codimension 
at least two, we have $X = X'$.
Theorem~\ref{firstresult} shows the
converse:

\begin{introthm} 
Every $A_2$-maximal (e.g., projective), 
normal variety with free finitely 
generated divisor class group and 
finitely generated total coordinate ring 
arises from a bunched ring.
\end{introthm}

As indicated, a central aim is to describe  
geometric properties of the variety $X$ 
arising from a bunched ring
$(R,\mathfrak{F},\Phi)$ in terms of 
the defining data.
A first observation is that the system $\mathfrak{F}$
provides a stratification for the variety $X$.
The generators $f_{i} \in \mathfrak{F}$ induce certain
divisors $D_{X}^{i} \subset X$, 
see Section~\ref{toricambient},
and for any face $\gamma_{0}$ of the positive orthant 
$\gamma := \QQ_{\ge 0}^{r}$,
one may consider the locally closed set
$$
X(\gamma_{0})
\; := \; 
\bigcap_{e_{i} \not\in \gamma_{0}} D_{X}^{i}
\; \setminus \; \bigcup_{e_{i} \in \gamma_{0}}  D_{X}^{i}
\; \subset \; X.
$$

We call a face $\gamma_{0} \preceq \gamma$ 
{\em relevant\/}
if the corresponding $X(\gamma_{0})$ is nonempty,
and we denote the set of relevant faces by 
$\rel(\Phi)$.
Thus, $X$ is stratified by the $X(\gamma_{0})$,
where $\gamma_{0}$ runs through the 
relevant faces.
Moreover, we need the linear map 
$Q \colon \ZZ^{r} \to K$
sending the canonical base vector $e_{i}$ to
the degree of $f_{i}$.
A first result shows that the stratification 
of $X$ is equisingular in a raw sense,
see Theorem~\ref{singularities}:

\begin{introthm}
Consider a stratum $X(\gamma_{0}) \subset X$
and a point $x \in X(\gamma_{0})$.
Then $X$ is $\QQ$-factorial (factorial) at $x$
if and only if $Q(\gamma_{0})$ is of full dimension
($Q$ maps $\lin(\gamma_{0}) \cap E$ onto $K$).
\end{introthm}

In general, there is no way
to describe the
smooth points in combinatorial terms.
Nevertheless, 
if the ring $R$ is ``regular enough'' 
in the sense that $\rq{X}$ 
is smooth, then the factorial 
points of $X$ are even smooth, 
see Proposition~\ref{smoothchar}.
This criterion
applies to all our examples. Furthermore,
in this setting, the singularities
of $X$ are at most rational.

If we pick out the set theoretically minimal 
cones from the relevant faces,
we arrive at the {\em covering collection} 
$\cov(\Phi) \subset \rel(\Phi)$.
In terms of this covering collection, 
we determine in Corollary~\ref{picarddescr}
the Picard group.

\begin{introthm}
In the divisor class group $K = \Cl(X)$, the 
Picard group of $X$ is given by
$$
\Pic(X) = \bigcap_{\gamma_{0} \in \cov(\Phi)} Q(\lin(\gamma_{0}) \cap E).
$$
\end{introthm}

Moreover, in Theorems~\ref{movingcone} 
and~\ref{amplecone}, we describe the 
effective cone, the moving cone,
the semiample cone and, finally, the ample
cone. The statement generalizes corresponding 
results on toric varieties, compare~\cite{BeHa2}
and~\cite{Tsch}.

\begin{introthm}
In the divisor class group $K = \Cl(X)$, we have the
following descriptions:
$$
\begin{array}{ccc}
\displaystyle \Eff(X) \; = \; Q(\gamma),
& \qquad &
\displaystyle \Mov(X) 
\; = 
\bigcap_{\gamma_0 \text{\ facet of \ } \gamma} Q(\gamma_0),
\\
 & & 
\\
\displaystyle \SAmple(X) \; = \; \bigcap_{\tau \in \Phi} \tau,
& \qquad &
\displaystyle \Ample(X) \; = \; \bigcap_{\tau \in \Phi} \tau^{\circ}.
\end{array}
$$
\end{introthm}

An immediate application is the description of the 
Mori cone given in Corollary~\ref{moricone}.
As a further application, we give in Section~\ref{kleichev} 
an effective Kleiman-Chevalley Criterion, 
which in the case of varieties
with finitely generated coordinate ring sharpens 
considerations presented in~\cite{Ko}:
suppose that $X$ is $\QQ$-factorial and complete, 
and let $k$ be the rank of the divisor class group of~$X$.
Then $X$ is projective if and 
only if any $k$ points of $X$ admit a common affine
neighbourhood in~$X$.

If the variety $X$ is sort of an ``intrinsic complete intersection'',
then one can determine its canonical divisor class in $K = \Cl(X)$,
and thus obtains Fano criteria, 
see Propositions~\ref{candivdescr} and ~\ref{fanocrit}:

\begin{introthm}
Let $d := r - \dim(X) - \dim(K)$, and 
suppose that $g_{1}, \ldots, g_{d}$ 
are $K$-homogeneous generators 
for the relations between 
$f_{1}, \ldots, f_{r}$. 
Then 
$$ 
\sum \deg(g_{j}) - \sum \deg(f_{i})
$$
is the canonical divisor class of $X$ in $K = \Cl(X)$.
Moreover, $X$ is a Fano variety if and only if
$$ 
\sum \deg(f_{i}) - \sum \deg(g_{j})
\; \in \;
\bigcap_{\tau \in \Phi} \tau^{\circ}
\cap
\bigcap_{\gamma_{0} \in \cov(\Phi)} Q(\lin(\gamma_{0}) \cap E).
$$
\end{introthm}

The language of bunched rings can as well
be used to produce and study (new) examples.
To indicate how this may run, we consider in 
Section~\ref{intrinsic} the case of coordinate rings 
defined by a single relation.
In that situation, the collection of relevant faces 
can be easily and explicitly determined in terms of 
the coefficients of the defining equation. 
For the case of a divisor class group of rank two
and coordinate rings defined by a quadric, we 
provide a classification extending Kleinschmidt's
classification of toric varieties with few invariant
divisors~\cite{Kl}, see 
Theorem~\ref{quadrklein}:

\begin{introthm}
The smooth intrinsic quadrics of full rank
with divisor class group $\ZZ^{2}$ and only 
constant functions arise from bunched rings 
$(R,\mathfrak{F},\Phi)$
with $\Phi$ given by a figure as follows:
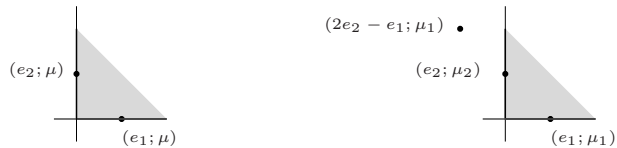
\begin{center}
\begin{picture}(0,0)%
\includegraphics{picdim2.pstex}%
\end{picture}%
\setlength{\unitlength}{1243sp}%
\begingroup\makeatletter\ifx\SetFigFont\undefined%
\gdef\SetFigFont#1#2#3#4#5{%
  \reset@font\fontsize{#1}{#2pt}%
  \fontfamily{#3}\fontseries{#4}\fontshape{#5}%
  \selectfont}%
\fi\endgroup%
\begin{picture}(11733,2793)(901,-2842)
\put(901,-1411){\makebox(0,0)[lb]{\smash{\SetFigFont{6}{7.2}{\familydefault}{\mddefault}{\updefault}$(e_2;\mu)$}}}
\put(3151,-2761){\makebox(0,0)[lb]{\smash{\SetFigFont{6}{7.2}{\familydefault}{\mddefault}{\updefault}$(e_1;\mu)$}}}
\put(7201,-511){\makebox(0,0)[lb]{\smash{\SetFigFont{6}{7.2}{\familydefault}{\mddefault}{\updefault}$(2e_2-e_1;\mu_1)$}}}
\put(11701,-2761){\makebox(0,0)[lb]{\smash{\SetFigFont{6}{7.2}{\familydefault}{\mddefault}{\updefault}$(e_1;\mu_1)$}}}
\put(9001,-1411){\makebox(0,0)[lb]{\smash{\SetFigFont{6}{7.2}{\familydefault}{\mddefault}{\updefault}$(e_2;\mu_2)$}}}
\end{picture}

\end{center}
All these varieties are projective.
They are Fano if and only if $\Phi$ arises
from a figure of the left hand side. 
Moreover, different figures belong to
non-isomorphic varieties.
\end{introthm}

In these pictures, the thick points indicate the 
degrees of the generators $f_i \in \mathfrak{F}$;
note that these degrees may occur with 
multiplicity, which is indicated by the $\mu$'s.
The shadowed cones are the elements of the 
respective ``$\mathfrak{F}$-bunches'' $\Phi$; 
in the above situation, each of these $\Phi$ 
consists of a single cone.

The proof is elementary: it is an interplay
of the combinatorics of bunches with the symmetry
conditions imposed by the defining quadric.
Classification problems for class groups of higher 
rank, and defining equations of higher degree(s), 
seem to pose new, elementary, but nontrivial 
combinatorial challenges.

Finally, we take a closer look at complete surfaces 
$X$ with finitely generated total coordinate ring 
$R := \mathcal{R}(X)$. 
First we discuss general properties. 
It turns out that $X$ is already determined 
by $R$, see Proposition~\ref{isoring2isosurf}.
Thus, in the description $X = X(R,\mathfrak{F},\Phi)$,
there is a certain freedom of choice
for $\mathfrak{F}$ and $\Phi$.
This is the key to the following observation, 
see Corollary~\ref{surfproj}:

\begin{introprop}
Let $X$ be a normal complete surface
with free finitely generated divisor class
group and finitely generated total coordinate ring.
If any two non-factorial singularities of
$X$ admit a common affine neighbourhood
in $X$, then $X$ is projective.
\end{introprop}

Then we consider some explicit examples,
and indicate how the results presented before 
can be used for explicit computations.
We study a certain two-dimensional 
Mumford quotient of the Grassmannian
$G(2,4)$ and determine its singularities,
the Picard group, the (semi)-ample cone,
and the canonical divisor.
Moreover, we consider the minimal resolution 
of the $E_6$ cubic surface, studied recently 
by Hasset and Tschinkel~\cite{Tsch}.

\tableofcontents

\goodbreak

We would like to thank the referee for his careful
reading, and a long list of very helpful comments 
and valuable suggestions for improvement.

\section{Bunches and toric varieties}
\label{bunches}

The combinatorial description of toric varieties
in terms of bunches presented in~\cite{BeHa2}
is an essential ingredient of the present paper.
In this section we recall the necessary concepts.
It should be mentioned that we only need a special case: 
here the word ``bunch'' will refer to what we called
in~\cite{BeHa2} a ``free standard bunch''.

In the whole paper, 
we mean by a variety
an integral
scheme over an algebraically closed field $\KK$
of characteristic zero, and the word ``point''
refers to a closed point. Concerning toric varieties, 
we assume the reader to 
be familiar with the basic theory, see
for example~\cite{Fu}.

To begin, let us fix terminology and notation from 
convex geometry. 
A lattice is a finitely generated free abelian
group. Given a lattice $N$, we denote by 
$N_{\QQ} := N \otimes_{\ZZ} \QQ$ the associated
rational vector space. For a homomorphism
$P \colon F \to N$ of lattices, the induced
map $F_{\QQ} \to N_{\QQ}$ is again denoted by
$P$.

By a cone in a lattice $N$ we mean a 
(not necessarily strictly)
convex polyhedral cone $\sigma \subset N_{\QQ}$.
The relative interior of a cone $\sigma$
is denoted by $\sigma^{\circ}$, its dual cone by 
$\sigma^{\vee}$, and its orthogonal space by 
$\sigma^{\perp}$.
If $\sigma_{0}$ is a face of a cone $\sigma$,
then we write $\sigma_{0} \preceq \sigma$.
A cone is called simplicial, if it is generated
by linearly independent vectors,
and it is called regular, if it is generated
by part of a lattice basis.

A fan in a lattice $N$ is a finite collection
$\Delta$ of strictly convex cones in $N$ such 
that for $\sigma \in \Delta$ also every
$\sigma_{0} \preceq \sigma$ belongs to $\Delta$,
and for any two $\sigma_{1}, \sigma_{2} \in \Delta$
we have $\sigma_{1} \cap \sigma_{2} \preceq \sigma_{i}$.
Given a fan $\Delta$, we denote by 
$\Delta^{\max}$ the set of its maximal cones
and by $\Delta^{(k)}$ the set of its $k$-dimensional
cones.

Let us now recall the basic concepts and constructions
from~\cite{BeHa2}.
By a {\em projected cone\/} we mean here a pair 
$(E \topto{Q} K, \gamma)$,
where $Q \colon E \to K$ is a surjection of lattices, 
and $\gamma \subset E_{\QQ}$ is a regular polyhedral 
convex cone of full dimension. 
A {\em projected face\/} is the image $Q(\gamma_{0})$
of a face $\gamma_{0} \preceq \gamma$.

\begin{definition}
\label{bunchdef}
A {\em bunch\/} in a projected cone 
$(E \topto{Q} K, \gamma)$ is a nonempty
collection $\Theta$ of projected faces 
with the following properties:
\begin{enumerate}
\item A projected face $\sigma$ belongs to $\Theta$ if and only if
$\emptyset 
\ne \sigma^{\circ} \cap \tau^{\circ} 
\ne \tau^{\circ}$
holds for all $\tau \in \Theta$  with 
$\tau \ne \sigma$.
\item For each facet $\gamma_0 \preceq \gamma$, the images of the 
primitive generators of $\gamma_0$ generate the lattice $K$, and 
there is a $\tau \in \Theta$ with $\tau^\circ \subset Q(\gamma_0)^\circ$.
\end{enumerate}
\end{definition}

Note that a collection $\Theta$ of projected 
faces satisfies~\ref{bunchdef}~(i) if and only if 
any two cones of $\Theta$ overlap,
$\Theta$ is irredundant in the sense that no 
two cones of $\Theta$ are strictly contained 
one in the other, 
and $\Theta$ is maximal in the sense that any 
projected face that overlaps with each cone
of $\Theta$ must contain a cone of $\Theta$.

The basic construction of~\cite{BeHa2} 
associates to a given bunch $\Theta$ living in a projected 
cone $(E \topto{Q} K, \gamma)$ a toric variety $Z_{\Theta}$.
This construction is purely combinatorial;
using a certain Gale duality, one transforms 
the bunch $\Theta$ into a fan $\Delta$ and 
thus obtains a toric variety.

Here comes the exact procedure: 
consider the {\em dual projected cone\/} 
$(F \topto{P} N, \delta)$, 
that means that we have $F = \Hom(E,\ZZ)$, 
the map $P \colon F \to N$ 
is dual to the inclusion $M \to E$ of $M := \ker(Q)$,
and $\delta := \gamma^\vee$ is the dual cone of $\gamma$.
Then we have the 
{\em face correspondence}
$$
{\rm faces}(\gamma) \to {\rm faces}(\delta),
\qquad
\gamma_0 \mapsto \gamma_0^* := \gamma_{0}^\perp \cap \delta.
$$

The {\em covering collection\/} $\cov(\Theta)$ of
the bunch $\Theta$ is the set of all faces 
$\gamma_0 \preceq \gamma$ that are minimal 
with respect to the property that $Q(\gamma_0)$ 
contains a cone of $\Theta$.
Via the face correspondence,
we define a fan in $F$:
$$ 
\rq{\Delta} := \{\delta_0; \; \delta_0 \preceq \gamma_0^*
  \text{ for some } \gamma_0 \in \cov(\Theta)\}.
$$

Property~\ref{bunchdef}~(i) implies that $\rq{\Delta}$
is a {\em maximal projectable fan\/} in $(F \topto{P} N, \delta)$
in the sense that it consists of faces of $\delta$ and is 
maximal with respect to the following property: 
any two maximal cones of $\rq{\Delta}$ can be 
separated by linear forms that are invariant with 
respect to $L := \Hom(K,\ZZ)$, 
i.e., that vanish on $L$.

As a consequence of Property~\ref{bunchdef}~(ii),
one obtains that the projected faces $P(\delta_0)$, 
where $\delta_0 \in \rq{\Delta}^{\max}$, are strictly convex.
Thus, the  $P(\delta_0)$ are the maximal
cones of a fan $\Delta$ in the lattice $N$, 
the {\em quotient fan\/} of $\rq{\Delta}$.
In conclusion, transforming $\Theta$ into $\Delta$ 
runs according to the following scheme:
$$
\xymatrix@R-20pt{
0 \ar[r] 
&
L \ar[r] 
&
F \ar[r]^{P}
&
N \ar[r] 
&
0
\\
&
&
{\rq{\Delta}}^{\max} \ar[r]
&
{\Delta^{\max}}
&
\\
&&&&
\\
&
{\Theta}  \ar@{-->}[r]
&
{\cov(\Theta)} \ar[uu]_{\gamma_0 \mapsto \gamma_0^*}
& 
&
\\
0 
& 
K \ar[l]
&
E \ar[l]^{Q}
&
M \ar[l]
&
0 \ar[l]
}
$$
where the two exact sequences are dual to each other,
and the dashed arrow indicates the passage to the 
covering collection. Note that this procedure may 
also be reversed. In fact, in~\cite[Sec.~5]{BeHa2}
it is shown that it provides a one-to-one 
correspondence between bunches and maximal projectable 
fans.

\begin{definition}
The toric variety $Z_{\Theta}$
associated to the bunch $\Theta$
is the toric variety arising from the 
fan $\Delta$ in the lattice $N$.
\end{definition}

In the sequel, we shall
often omit the subscript $\Theta$
and just speak of the toric variety 
$Z$ associated to the bunch $\Theta$.
Recall from~\cite[Theorem~7.4 and Prop.~10.6]{BeHa2}
that $Z_{\Theta}$ is nondegenerate
in the sense that 
$\mathcal{O}^*(Z_{\Theta}) = \KK^*$ 
holds,
and that $Z_{\Theta}$ has free divisor 
class group, canonically isomorphic 
to $K$.

We need a description of the orbit stratification 
of the toric variety $Z$ associated to $\Theta$. 
First recall from~\cite{Fu} that the cones 
of $\Delta$ correspond to the orbits
of the big torus $T_{Z} \subset Z$ via
$\sigma \mapsto T_{Z} \mal z_{\sigma}$.
Here $z_{\sigma} \in Z$ denotes the 
distinguished point of the affine chart 
$Z_{\sigma} = \Spec(\KK[\sigma \cap N])$
of $Z$; it is determined by 
$$
\chi^{u}(z_{\sigma}) = 1 
\quad \text{if } 
u \in \sigma^{\perp},
\qquad
\chi^{u}(z_{\sigma}) = 0 
\quad \text{if } 
u \in \sigma^{\vee} \setminus \sigma^{\perp}.
$$

We want to translate this into the language
of bunches.
Call a face $\gamma_{0} \preceq \gamma$ {\em relevant\/}
if $Q(\gamma_{0})^{\circ} \supset \tau^{\circ}$ 
holds for some $\tau \in \Theta$,
and denote by $\rel(\Theta)$ the collection
of relevant faces.
Then \cite[Prop.~5.11]{BeHa2} tells us
that the cones of the fan $\Delta$ are exactly
the images $P(\gamma_0^*)$, where 
$\gamma_0^* \in \rel(\Theta)$.
This gives the desired description of the orbit 
stratification:

\begin{proposition}
\label{orbitstrat}
There is a canonical bijection
$$
\rel(\Theta) \to \{T_{Z}\text{-orbits of } Z\},
\qquad
\gamma_{0} \mapsto  Z(\gamma_0) := T_{Z} \mal z_{P(\gamma_{0}^{*})}.
$$
\end{proposition}

We shall frequently need the fact 
that the toric morphism
$p_{Z} \colon \rq{Z} \to Z$ arising from the map
$P \colon F \to N$ of the fans 
$\rq{\Delta}$ and $\Delta$ is a 
{\em Cox construction}, 
that means that we have:

\begin{proposition}
\label{coxconstr}
Every ray of $\delta$ occurs in $\rq{\Delta}$,
and the map $P \colon F \to N$ induces bijections
$\rq{\Delta}^{\max} \to \Delta^{\max}$ and
$\rq{\Delta}^{(1)} \to \Delta^{(1)}$.
Moreover, $P$ maps the primitive generator 
of any ray $\rq{\varrho} \in \rq{\Delta}^{(1)}$ 
to a primitive lattice vector in $N$.
\end{proposition}

Let us interprete this also geometrically.
Consider the affine toric variety 
$\b{Z} := \Spec(\KK[\gamma \cap E])$
with its big torus 
$T_{\b{Z}} := \Spec(\KK[E])$.
The $K$-grading of $\KK[\gamma \cap E]$  
defines actions of $T := \Spec(\KK[K])$ on 
$\b{Z}$ and  
on the $T_{\b{Z}}$-invariant open subset 
$\rq{Z} \subset \b{Z}$ corresponding to
the fan $\rq{\Delta}$.
The geometric version of 
Proposition~\ref{coxconstr} is:

\begin{proposition}
\label{coxconstrgeo}
\begin{enumerate}
\item The toric morphism $p_{Z} \colon \rq{Z} \to Z$ is a 
 good quotient for the $T$-action, that means that it is affine, 
 $T$-invariant and defines an isomorphism
 $\mathcal{O}_{Z} \to p_{*}(\mathcal{O}_{\rq{Z}})^{T}$.
\item Let $W_{Z} \subset Z$ and $W_{\b{Z}} \subset \b{Z}$ 
be the unions of all at most one codimensional orbits
of the respective big tori. 
Then $W_{\b{Z}} = p_{Z}^{-1}(W_{Z})$ 
and $W_{\b{Z}} \subset \rq{Z}$ hold, 
and $T$ acts freely on $W_{\b{Z}}$.
\end{enumerate}
\end{proposition}

For later purposes, we note here an elementary observation
on existence of polytopal fans with prescribed set of rays;
compare also~\cite[Prop.~3.4]{Tsch}.

\begin{lemma}
\label{einsgeruest}
Let $N$ be a lattice, and let $v_1, \ldots, v_r \in N$ be 
vectors generating $N$ as a lattice and $N_{\QQ}$ as a cone.
Then there is a polytopal fan $\Delta$ in $N$ having
precisely $\QQ_{\ge 0}v_1, \ldots, \QQ_{\ge 0}v_r$ as its
rays. 
\end{lemma}

\proof
First note that there is a polytopal fan $\Delta_0$ in $N$
the rays of which form a subset of the prescribed choice,
say $\QQ_{\ge 0}v_1, \ldots, \QQ_{\ge 0}v_s$ for some $s \le r$,
see e.g.~\cite{OdPa}.
Then we may construct a new polytopal fan $\Delta_1$ 
from $\Delta_0$ by inserting the ray $\QQ_{\ge 0}v_{s+1}$ 
via stellar subdivision.
Iterating this process, we finally arrive at the desired 
fan.
\endproof

We shall also need the bunch theorical analogue of the 
above observation. It is the following:

\begin{lemma}
\label{1Gbunch}
Let $(E \topto{Q} K, \gamma)$ be a projected cone such
that $Q(\gamma)$ is strictly convex, 
and let $\gamma_1, \ldots, \gamma_r$ be the facets of 
$\gamma$. 
Suppose that each $Q(\lin(\gamma_i) \cap E)$ 
generates $K$ as a lattice,
and that for any two $\gamma_i$, $\gamma_j$ we have 
$Q(\gamma_i)^{\circ} \cap Q(\gamma_j)^{\circ} \ne \emptyset$.
Then 
$Q(\gamma_1)^{\circ} \cap \ldots \cap Q(\gamma_r)^{\circ}
\ne \emptyset$
holds.
\end{lemma}

\proof
First consider the dual projected cone $(F \topto{P} N, \delta)$.
Then the images $v_1, \ldots, v_r \in N$ of the primitive
generators of $\delta$ are as in Lemma~\ref{einsgeruest}.
Hence, we find a polytopal fan $\Delta$ in $N$ having
precisely $\QQ_{\ge 0}v_1, \ldots, \QQ_{\ge 0}v_r$ as its
rays. 

Now consider the bunch $\Theta$ in $(E \topto{Q} K, \gamma)$ 
corresponding to $\Delta$.
According to~\cite[Thm.~10.2]{BeHa2}, the intersection over
all relative interiors $\tau^\circ$, where $\tau \in \Theta$,
is nonempty. 
Since every $\QQ(\gamma_i)^{\circ}$ contains some $\tau^{\circ}$
the assertion follows. 
\endproof

We conclude this section with a few words towards
the visualization of bunches.
Look at the following example:
Let $K := \ZZ$, consider
the elements $w_{1} = 2$, $w_{2} = 4$,
and $w_{3} = 5$ in $K$.
Attach to each $w_{i}$ a multiplicity $\mu_{i}$
by setting $\mu_{1} := 1$, $\mu_{2} := 3$, and
$\mu_{3} := 2$.
Finally, set $\Theta := \{\QQ_{\ge 0}\}$.
Then these data fit into the following figure:

\medskip

\begin{center}

\begin{picture}(0,0)%
\includegraphics{bunchexam.pstex}%
\end{picture}%
\setlength{\unitlength}{2072sp}%
\begingroup\makeatletter\ifx\SetFigFont\undefined%
\gdef\SetFigFont#1#2#3#4#5{%
  \reset@font\fontsize{#1}{#2pt}%
  \fontfamily{#3}\fontseries{#4}\fontshape{#5}%
  \selectfont}%
\fi\endgroup%
\begin{picture}(5895,533)(1339,-511)
\put(6391,-511){\makebox(0,0)[lb]{\smash{\SetFigFont{9}{10.8}{\familydefault}{\mddefault}{\updefault}
$(5;2)$
}}}
\put(5041,-511){\makebox(0,0)[lb]{\smash{\SetFigFont{9}{10.8}{\familydefault}{\mddefault}{\updefault}
$(4;3)$
}}}
\put(4141,-511){\makebox(0,0)[lb]{\smash{\SetFigFont{9}{10.8}{\familydefault}{\mddefault}{\updefault}
$(2;1)$
}}}
\end{picture}

\end{center}

\medskip

We obtain an associated projected cone 
$(E \topto{Q} K, \gamma)$ by taking
the lattice $E := \ZZ^{6}$, 
the map $Q \colon E \to K$, 
defined by 
$e_{1} \mapsto 2$, $e_{2}, e_{3}, e_{4} \mapsto 4$, 
and $e_{5},e_{6} \mapsto 5$,
and finally, the cone $\gamma := \QQ_{\ge 0}^{6}$.
The associated toric variety turns out to be
the weighted projective space $\PP(2,4,4,4,5,5)$.

This principle allows to draw pictures also in general:
Any bunch can represented by 
a set $\{(w_{i};\mu_{i}); \; i \in I\}$ of lattice
vectors $w_{i}$ carrying multiplicities $\mu_{i}$, 
and a collection $\Theta$ of cones, each of which
generated by some of the $w_{i}$, 
such that the whole setup satisfies 
conditions analogue to those
of Definition~\ref{bunchdef}.

\section{Bunched rings and varieties}
\label{bunchedringsI}

In this section, we first put the concept of a bunch into a 
more general framework, and we introduce the notion of a 
bunched ring.
Then the central construction of the paper is presented:
a generalized multigraded proj-construction
that associates to any bunched ring a variety coming along 
with a certain stratification.

Let us fix the setup of the section.
$R$ is a factorial, finitely generated $\KK$-algebra, 
faithfully graded by some lattice $K \cong \ZZ^{k}$.
Here, faithfully graded means that $K$ is 
generated as a lattice by the degrees $w \in K$ admitting 
nontrivial homogeneous elements $f \in R_w$.

Moreover, $\mathfrak{F} =\{f_{1}, \dots, f_{r}\} \subset R$ is a 
system of homogeneous pairwise non associated 
nonzero prime elements generating $R$ as an algebra.
Since we assume the grading to be faithful,
the degrees $w_{i} := \deg(f_{i})$ generate the
lattice $K$.

The projected cone $(E \topto{Q} K, \gamma)$
{\em associated\/} to the system of generators 
$\mathfrak{F} \subset R$ consists of the surjection 
$Q$ of the lattices $E := \ZZ^{r}$ and $K$ sending 
the $i$-th canonical base vector $e_{i}$ to $w_{i}$, 
and the cone $\gamma \subset  E_{\QQ}$ generated by 
$e_{1}, \dots, e_{r}$.

\begin{definition}\label{fdefs}
Let $(E \topto{Q} K, \gamma)$ denote the projected cone associated to
$\mathfrak{F} \subset R$.
\begin{enumerate}    
\item \label{fface}    
We call $\gamma_{0} \preceq \gamma$ an
{\em $\mathfrak{F}$-face\/} if the product
over all $f_{i}$ with $e_{i} \in \gamma_{0}$
does not belong to the ideal
$\sqrt{\bangle{f_{j}; \; e_{j} \not\in \gamma_{0}}} \subset R$.
\item \label{fbunch}
By an {\em $\mathfrak{F}$-bunch\/} we mean a nonempty 
collection $\Phi$ of projected $\mathfrak{F}$-faces such that
\begin{itemize}
\item a projected $\mathfrak{F}$-face $\tau$ belongs to $\Phi$ if 
and only if for each $\tau \neq \sigma \in \Phi$ we have
$\emptyset \neq \tau^{\circ} \cap \sigma^{\circ} \neq 
\sigma^{\circ}$,
\item for each facet $\gamma_{0} \preceq \gamma$,
the image $Q(\gamma_{0} \cap E)$ generates the lattice $K$,
and there is a $\tau \in \Phi$ such that 
$Q(\gamma_{0})^{\circ} \supset \tau^{\circ}$ holds.
\end{itemize}
\end{enumerate}
\end{definition}

If we want to specify these data, then we also speak of the 
$\mathfrak{F}$-bunch in the projected cone 
$(E \topto{Q} K,\gamma)$.
We shall frequently use the following 
more geometric characterization of the $\mathfrak{F}$-faces:

\begin{remark}
\label{geomfface}
Let $\b{X} := \Spec(R)$. 
Then a face $\gamma_{0} \preceq \gamma$
is an $\mathfrak{F}$-face if and only if
there is a point $x \in \b{X}$ satisfying
$e_{i} \in \gamma_{0}
\; \Leftrightarrow \; 
f_{i}(x) \neq 0$.
\end{remark}

We come to the central definition of the paper:
the notion of a bunched ring.
As we shall see, this generalizes 
the usual notion of a bunch in a projected cone.

\begin{definition}
A {\em bunched ring\/} is a triple $(R,\mathfrak{F},\Phi)$, 
where
\begin{enumerate}
\item $R$ is a finitely generated factorial $\KK$-algebra
   with $R^{*} = \KK^{*}$, faithfully 
   graded by a finite dimensional lattice~$K$,
\item $\mathfrak{F} = \{f_{1}, \ldots, f_{r}\}$
   is a system of homogeneous pairwise nonassociated 
   non\-zero prime elements generating $R$ as an
   algebra,
\item $\Phi$ is an $\mathfrak{F}$-bunch in the projected
cone $(E \topto{Q} K, \gamma)$ associated to $\mathfrak{F}$. 
\end{enumerate}
\end{definition}

Let us illustrate this definition by means of 
some examples. 
The second one shows how the usual notion 
of a bunch discussed in the previous section fits 
into the framework of $\mathfrak{F}$-bunches
and bunched rings.

\begin{example}[Trivially bunched rings]
\label{triviallybunchedrings}
Any finitely generated factorial $\KK$-algebra 
$R$ with $R^* = \KK^*$ can be made trivially
into a bunched ring: take the trivial 
grading by $K=\{0\}$, any system 
$\mathfrak{F} = \{f_{1}, \ldots, f_{r}\}$ 
of pairwise nonassociated prime generators,
and the trivial $\mathfrak{F}$-bunch 
$\Phi = \{0\}$ in the projected cone   
$(E \topto{Q} K, \gamma)$, 
where $E := \ZZ^r$, and $\gamma \subset E_{\QQ}$ 
is the positive orthant.
\end{example}

\begin{example}[Bunched polynomial rings]
\label{bunchesasFbunches}
Let $R = \KK[T_{1}, \ldots, T_{r}]$.
Suppose that each variable $T_{i}$ is 
homogeneous with respect to the
$K$-grading of $R$, and let 
$\mathfrak{F} := \{T_{1}, \ldots, T_{r}\}$.
Then any bunch $\Theta$ in the projected cone
associated to~$\mathfrak{F}$ is an 
$\mathfrak{F}$-bunch, and thus defines
a bunched ring $(R,\mathfrak{F},\Theta)$.
\end{example}

\begin{example}%
[Homogeneous coordinates of the Grassmannian
$G(2,4)$]
\label{gzwovierI}
Take the lattice
$K := \ZZ$, and consider the $K$-graded 
factorial ring defined by
$$
R 
\; := \: 
\KK[T_{1}, \ldots, T_{6}]
/
\bangle{T_{1} T_{6} - T_{2}T_{5} + T_{3}T_{4}},
\qquad 
\deg(T_{i}) := 1, \ 1 \le i \le 6.
$$
The projected cone $(E \topto{Q} K,\gamma)$
associated to the system
$\mathfrak{F} := \{T_{1}, \ldots, T_{6}\}$
of generators is given by $E := \ZZ^{6}$, 
the map $Q \colon E \to K$ sending $e_{i}$ to $1$,
and the cone $\gamma := \QQ^{6}_{\ge 0}$.
Setting $\Phi := \{\QQ_{\ge 0}\}$, we
obtain a bunched ring
$(R,\mathfrak{F},\Phi)$.
\end{example}

Let us now present the basic construction of the paper.
As mentioned, it associates to a given bunched ring 
$(R,\mathfrak{F},\Phi)$ a variety $X(R,\mathfrak{F},\Phi)$.
The first step is to generalize the concepts of 
relevant faces and the covering collection defined
in Section~\ref{bunches}:

\begin{definition}\label{fcov}
Let $\Phi$ be an $\mathfrak{F}$-bunch in the
projected cone $(E \topto{Q} K,\gamma)$.
\begin{enumerate}
\item The collection $\rel(\Phi)$ of {\em relevant faces\/} 
    consists of those $\mathfrak{F}$-faces $\gamma_{0} \preceq \gamma$ 
    such that $Q(\gamma_{0})^{\circ} \supset \tau^{\circ}$ 
    holds for some $\tau \in \Phi$.
\item The {\em covering collection\/} $\cov(\Phi)$ of $\Phi$ 
is the set of all minimal cones of $\rel(\Phi)$.
\end{enumerate}
\end{definition}

Note that this gives indeed nothing new for a usual
bunch when we regard it as an $\mathfrak{F}$-bunch in 
the sense of Example~\ref{bunchesasFbunches}.

\begin{remark}\label{fcovoverlap}
Let $\Phi$ be an $\mathfrak{F}$-bunch with
associated projected cone $(E \topto{Q} K,\gamma)$, 
and let $\gamma_{0}, \gamma_{1} \in \cov(\Phi)$. 
Then $Q(\gamma_{0})^{\circ} \cap Q(\gamma_{1})^{\circ} \ne \emptyset$. 
\end{remark}

Let $(R,\mathfrak{F},\Phi)$ be a bunched ring,
and let $(E \topto{Q} K, \gamma)$ be the projected
cone associated to $\mathfrak{F}$.
To construct the variety $X(R,\mathfrak{F},\Phi)$,
we consider the affine variety $\b{X} := \Spec(R)$, 
the action of the torus $T := \Spec(\KK[K])$ on
$\b{X}$ induced by the $K$-grading of $R$,
and the invariant open subset:
\begin{equation}
\label{Xhutdef}
\rq{X}
\; = \;
\bigcup_{\gamma_{0} \in \cov(\Phi)} \b{X}_{{\gamma_{0}}},
\qquad \text{where } 
\b{X}_{\gamma_{0}} = \b{X}_{f^{u}} 
\text{ for some }
u \in \gamma_{0}^{\circ}.
\end{equation}

Here we use the notation 
$f^{u} := f_{1}^{u_{1}} \dots f_{r}^{u_{r}}$ for 
$u \in E = \ZZ^{r}$. 
Note that the open affine subsets $\b{X}_{\gamma_{0}} \subset \b{X}$ 
do not depend on the particular choices of
the lattice vectors $u \in \gamma_{0}^{\circ}$.
The basic observation is the following:

\begin{proposition}\label{xasgoodquot}
There is a good quotient 
$\rq{X} \to \rq{X} \quot T$
for the $T$-action on $\rq{X}$.
\end{proposition}

\proof 
First consider any two cones 
$\gamma_{i}, \gamma_{j} \in \cov(\Phi)$.
Since $Q(\gamma_{i})^{\circ} \cap Q(\gamma_{j})^{\circ}$
is nonempty, we find $u^{i} \in \gamma_{i}^{\circ}$ 
and $u^{j} \in \gamma_{j}^{\circ}$ 
with $Q(u^{i}) = Q(u^{j})$.
Let $f_{i}, f_{j} \in R$ denote the functions corresponding to 
$u^{i}, u^{j}$.
Then we obtain a commutative diagram
$$
\xymatrix{
{\b{X}_{f_{i}}} \ar[d]^{\quot T} 
& {\b{X}_{f_{i}f_{j}}} \ar[l]
                       \ar[r] 
                       \ar[d]^{\quot T} 
& {\b{X}_{f_{j}}} \ar[d]^{\quot T} \\
{{X}_{f_{i}}} 
& {{X}_{f_{i}f_{j}}} \ar[l]
\ar[r]
& {{X}_{f_{j}}}
}
$$ 
where the upper horizontal maps are open embeddings,
the downwards maps are good quotients for the respective
$T$-actions,
and the lower horizontal arrows indicate the induced
morphisms of the affine quotient spaces.

By the choice of $f_{i}$ and $ f_{j}$, the quotient $f_{j}/f_{i}$
is an invariant function on $\b{X}_{f_{i}}$, and
the inclusion $\b{X}_{f_{i}f_{j}} \subset \b{X}_{f_{i}}$
is just the localization by $f_{j}/f_{i}$.
Since $f_{j}/f_{i}$ is invariant, the latter holds
as well for the quotient spaces; that means that
the map $X_{f_{i}f_{j}} \to X_{f_{i}}$ is localization
by $f_{j}/f_{i}$.

In particular, $X_{f_{i}f_{j}} \to X_{f_{i}}$ are open embeddings,
and we may glue the maps $\b{X}_{f_{i}} \to X_{f_{i}}$ along
$\b{X}_{f_{i}f_{j}} \to X_{f_{i}f_{j}}$.
This gives a good quotient $\rq{X} \to \rq{X} \quot T$
for the $T$-action.
Note that the quotient space is separated, because
the multiplication maps 
$\mathcal{O}(X_{f_{i}}) \otimes \mathcal{O}(X_{f_{j}}) 
\to \mathcal{O}(X_{f_{i}f_{j}})$ are surjective.
\endproof

\begin{definition}\label{br2vardef}
In the notation of Proposition~\ref{xasgoodquot}, 
the {\em variety associated to the bunched ring\/} 
$(R,\mathfrak{F},\Phi)$ is the quotient space 
$X(R,\mathfrak{F},\Phi) := \rq{X} \quot T$.
\end{definition}
 
For the subsequent studies and results, 
it is important to observe 
that the variety $X := X(R,\mathfrak{F},\Phi)$ 
associated to a bunched ring $(R,\mathfrak{F},\Phi)$
inherits  a certain stratification
from its construction.

Consider the quotient map $p_X \colon \rq{X} \to X$.
Then, for any $\gamma_0 \in \rel(\Phi)$, 
we have the open affine subset 
$X_{\gamma_0} := \b{X}_{\gamma_0} \quot T$
of $X = \rq{X} \quot T$.
Note that
$p_X^{-1}(X_{\gamma_0}) = \b{X}_{\gamma_0}$
holds, 
and for a further $\gamma_1 \in \rel(\Phi)$ with
$\gamma_0 \preceq \gamma_1$, we have
$X_{\gamma_1} \subset X_{\gamma_0}$ and
$\b{X}_{\gamma_1} \subset \b{X}_{\gamma_0}$.   
We define locally closed subsets
in $X$ and $\rq{X}$:
$$ 
X(\gamma_{0}) := X_{\gamma_{0}} \setminus \bigcup_{\gamma_{0}
\prec \gamma_{1} \in \rel(\Phi)} X_{\gamma_{1}},
\qquad
\rq{X}(\gamma_0) := p_X^{-1}(X(\gamma_0)). 
$$

\begin{proposition}\label{strati}
The variety $X := X(R,\mathfrak{F},\Phi)$ associated 
to the bunched ring $(R,\mathfrak{F},\Phi)$ is a disjoint union of 
nonempty locally closed subsets:
$$ X = \bigcup_{\gamma_{0} \in \rel(\Phi)} X(\gamma_{0}).$$
Moreover, the strata $\rq{X}(\gamma_0)$, $\gamma_0 \in \rel(\Phi)$, 
of the lifted 
stratification on $\rq{X}$ are explicitly given by
$$
\rq{X} (\gamma_0) = \bigl\{x \in \b{X}; \; \prod_{e_{i} \in \gamma_{0}} f_{i}
(x) \ne 0, \; \forall \gamma_0 \prec \gamma_{1} \in \rel(\Phi) \colon  
\prod_{e_{j} \in \gamma_{1}} f_{j} (x) = 0 \bigr\}.
\qed
$$
\end{proposition}

In order to illustrate the construction of the variety
associated to a bunched ring, we look again at 
Examples~\ref{triviallybunchedrings},
\ref{bunchesasFbunches} and~\ref{gzwovierI}.

\begin{example}
Let $(R,\mathfrak{F},\Theta)$ be a trivially bunched
ring as in Example~\ref{triviallybunchedrings}.
Then the associated variety is 
$X(R,\mathfrak{F},\Theta) = \Spec(R)$.
\end{example}

\begin{example}
Let $(R,\mathfrak{F},\Theta)$ be a bunched
polynomial ring as in Example~\ref{bunchesasFbunches}.
Then $X(R,\mathfrak{F},\Theta)$ is the toric variety 
corresponding to $\Theta$, and the stratification 
of $X(R,\mathfrak{F},\Theta)$ is the orbit 
stratification.
\end{example}

\begin{example}\label{gzwovierIII}
Consider the bunched ring 
$(R,\mathfrak{F}, \Phi)$ of Example \ref{gzwovierI}.
Then the associated variety $X(R,\mathfrak{F}, \Phi)$ 
equals the Grassmannian $G(2,4)$.

Indeed, $\b{X} := \Spec(R)$ is the affine cone over $G(2,4)$, 
see~e.g.~\cite[Sec.~5.13]{Wh},
and the action of $T := \Spec(\KK[K])$
on $\b{X}$ is the usual $\CC^{*}$-action.
Moreover, $\rq{X} = \b{X} \setminus \{0\}$,
because $\cov(\Phi)$ consists of 
all rays of $\gamma$, and thus
$\rq{X}$ consists of all $x \in \b{X}$
admitting an $f_{j}$ with $f_{j}(x) \neq 0$. 
Thus, $X(R,\mathfrak{F}, \Phi) = \rq{X} \quot T$ 
equals $G(2,4)$.
\end{example}

\section{Toric ambient varieties}
\label{toricambient}

In the previous section, we presented the basic
construction of the paper, 
which associates to any bunched ring 
a stratified variety.
In this section we show that these varieties admit
closed embeddings into toric varieties such that
the stratification is induced from the orbit 
stratification.

Throughout the whole section, we fix a bunched ring
$(R,\mathfrak{F},\Phi)$, and we denote the projected cone
associated to $\mathfrak{F}$ by $(E \topto{Q} K, \gamma)$. 
The construction of toric ambient varieties is based on 
the following observation:

\begin{lemma}
\label{Phi2Theta}
There exists a bunch $\Theta$ in 
$(E \topto{Q} K,\gamma)$ such that
$\sigma^{\circ} \cap \tau^{\circ} \ne \emptyset$
holds for any $\sigma \in \Theta$ and 
any $\tau \in \Phi$.
Moreover, for every such bunch $\Theta$ one 
has 
$$
\rel(\Phi) 
\; = \; 
\{
\gamma_{0} \in \rel(\Theta); \; 
\gamma_{0} \text{ is an } \mathfrak{F} \text{-face}
\}.
$$
\end{lemma}

Note that for $\Theta$ and $\Phi$ as in this Lemma, 
every cone of $\Phi$ contains in 
its relative interior the relative interior
of some cone of $\Theta$.

\proof[Proof of Lemma~\ref{Phi2Theta}]
To verify the first statement, we construct $\Theta$ as follows:
Let $\Theta_{0} := \Phi$.
If there is a projected face $\tau \notin \Theta_{0}$ such that 
$\tau^{\circ} \cap \sigma^{\circ} \neq \emptyset$
holds for all $\sigma \in \Theta_{0}$, 
then we set $\Theta_{1} := \Theta_{0} \cup \{\tau\}$.
After finitely many such steps, we arrive
at a collection $\Theta_{m}$, which can no longer be enlarged
in the above sense.
Then the set $\Theta$ of minimal cones of $\Theta_{m}$
is a bunch as desired. 

We turn to the second statement.
Let $\gamma_0 \in \rel(\Phi)$. 
Then $\gamma_0$ is an $\mathfrak{F}$-face. 
Moreover,
$Q(\gamma_0)^{\circ} \supset \tau^{\circ}$ 
holds for some $\tau \in \Phi$. 
Hence $Q(\gamma_0)^{\circ}$ intersects
all $\sigma^{\circ}$, where $\sigma \in \Theta$.
By the defining property~\ref{bunchdef}~(i) 
of a bunch, $Q(\gamma_0)^{\circ} \supset \sigma^{\circ}$ 
holds for some $\sigma \in \Theta$.
This in turn implies $\gamma_0 \in \rel(\Theta)$. 
 
Conversely,
let $\gamma_0 \in \rel(\Theta)$ be an 
$\mathfrak{F}$-face.
Then there is a $\sigma \in \Theta$ such that 
$Q(\gamma_0)^{\circ}$ contains $\sigma^{\circ}$.
Thus the open projected 
$\mathfrak{F}$-face $Q(\gamma_0)^{\circ}$ 
meets $\tau^{\circ}$ for every $\tau \in \Phi$. 
By the properties~\ref{fdefs}~(ii) of an 
$\mathfrak{F}$-bunch, we have 
$Q(\gamma_0)^{\circ} \supset \tau^{\circ}$ 
for some $\tau \in \Phi$.
But this means $\gamma_0 \in \rel(\Phi)$.
\endproof

Now, fix any bunch $\Theta$ in 
$(E \topto{Q} K,\gamma)$ as 
provided by Lemma~\ref{Phi2Theta}.
Let $e_{1}, \ldots, e_{r}$ denote the
canonical base vectors of $E = \ZZ^{r}$.
The semigroup algebra $\KK[\gamma \cap E]$ 
becomes $K$-graded by defining the degree of 
$\chi^{e_{i}}$ to be $Q(e_{i})$.
Moreover, we have a $K$-graded surjection
$$ 
\KK[\gamma \cap E] \to R,
\qquad 
\chi^{e_{i}} \mapsto f_{i}.
$$

Geometrically, this corresponds to a closed 
embedding of $\b{X} := \Spec(R)$ into the
toric variety 
$\b{Z} := \Spec(\KK[\gamma \cap E]) = \KK^{r}$,
and this embedding is equivariant with respect 
to the actions of the torus 
$T := \Spec(\KK[K])$ on $\b{X}$ and $\b{Z}$
defined by the $K$-gradings of $R$ and
$\KK[\gamma \cap E]$.

Let $\rq{\Delta}$ be the maximal projectable fan
associated to the bunch $\Theta$, and
let $\rq{Z} \subset \b{Z}$ be the corresponding
open toric subvariety.
Let $Z$ be the toric variety arising from 
the quotient fan $\Delta$ of $\rq{\Delta}$,
i.e., there is a good quotient 
$p_Z \colon \rq{Z} \to Z$. 
For $\gamma_0 \in \rel(\Theta)$, let 
$z(\gamma_0) \in Z$ denote the corresponding 
distinguished point,
that means that the orbit stratification 
reads as
$$
Z 
\; = \; 
\bigcup_{\gamma_0 \in \rel(\Theta)} T_Z \mal z(\gamma_0).
$$

Now, recall from the preceding section that 
$X := X(R,\mathfrak{F},\Theta)$ was defined 
as the good quotient of an invariant open subset 
$\rq{X} \subset \b{X}$ by the $T$-action.
The meaning of the toric varieties $\rq{Z}$
and $Z$ is the following:

\begin{proposition}\label{xintv}
We have $\rq{X} = \b{X} \cap \rq{Z}$,
and there is a commutative diagram,
where the horizontal maps are closed embeddings:
\begin{equation}
\label{bunchedring2variety}
\vcenter{
\xymatrix{
{\rq{X}} \ar[d]^{p_{X}}_{\quot T} \ar[r]^{\rq{\imath}} 
& {\rq{Z}} \ar[d]_{p_{Z}}^{\quot T} \\
X \ar[r]_{\imath} & Z
}}
\end{equation}
Moreover, the stratifications of $X$ and $Z$ are compatible: 
a cone $\gamma_0 \in \rel(\Theta)$ satisfies 
$X \cap T_Z \mal z(\gamma_0) \ne \emptyset$ if and only if
$\gamma_0$ is an $\mathfrak{F}$-face; and in that case one 
has 
$$
X(\gamma_0) = X \cap T_Z \mal z(\gamma_0).
$$
\end{proposition}

Note that a toric variety as in this proposition
is not determined by its properties. 
Nevertheless, the embedding construction will
play an important role, and hence we define: 

\begin{definition}
Let $(R,\mathfrak{F},\Phi)$ be a bunched ring,
and let $\Theta$ be a bunch as in Lemma~\ref{Phi2Theta}.
Then we say that $\Theta$ {\em extends\/} $\Phi$,
and that the toric variety $Z$ defined by $\Theta$
is an {\em ambient\/} toric variety of 
$X(R,\mathfrak{F},\Phi)$.
\end{definition}

\proof[Proof of Proposition~\ref{xintv}] 
Consider the stratification of the toric variety $\rq{Z}$ 
obtained by lifting the orbit stratification of $Z$ with respect
to the quotient map $p_Z \colon \rq{Z} \to Z$:
$$
\rq{Z} 
\; = \; 
\bigcup_{\gamma_0 \in \rel(\Theta)} \rq{Z}(\gamma_0),
\quad
\text{where }
\rq{Z}(\gamma_0) := p_Z^{-1}(T_Z \mal z(\gamma_0)).
$$
One can express this lifted stratification as well explicitly,
in terms of the coordinates 
$\chi^{e_{i}} = z_{i}$ of $\b{Z} = \KK^{r}$.
Namely, one directly verifies
$$
\rq{Z} (\gamma_0) =
\bigl\{ z \in \b{Z}; \; \prod_{e_{i} \in \gamma_{0}} z_{i}
 \ne 0, \; \forall \gamma_{0} \prec \gamma_{1} \in \rel(\Theta): 
\; \prod_{e_{j} \in \gamma_{1}} z_{j} = 0 \bigr\}. 
$$

In view of the description of the stratification of $X$
given in Proposition~\ref{strati}, 
we only have to verify the following claim: 
for $\gamma_0 \in \rel(\Theta)$, one has 
$\b{X} \cap \rq{Z}(\gamma_0) \ne \emptyset$ if and only 
if $\gamma_0$ is an $\mathfrak{F}$-face.
One implication is easy: 
if $\gamma_0$ is an $\mathfrak{F}$-face, 
then the $T_{\rq{Z}}$-orbit of 
the distinguished point 
$\rq{z}(\gamma_0) \in \rq{Z}$ 
corresponding to $\gamma_{0}^{*}$
lies in $\rq{Z}(\gamma_{0})$,
and it contains a point of
$\b{X} \cap \rq{Z}(\gamma_0)$.

For the converse, 
note first that 
$T_{\rq{Z}} \mal \rq{z}(\gamma_0)$
is closed in $\rq{Z}(\gamma_{0})$ and 
that it maps onto 
$T_Z \mal z(\gamma_0)$.
Now, consider a point 
$x \in \b{X} \cap \rq{Z}(\gamma_0)$. 
Then $p_{Z}(x)$ belongs to 
$T_Z \mal z(\gamma_0)$. 
Hence, by basic properties of good quotients,
the closures of $T \mal x$ and
$T_{\rq{Z}} \mal \rq{z}(\gamma_0)$ 
have a point in common, say $x'$.
By $T$-invariance, $x' \in \b{X}$ holds,
and by construction, we have 
$f_i(x') \ne 0$ if and only if
$e_i \in \gamma_0$.
\endproof

\begin{example}\label{gzwovierIV}
Consider the bunched ring
$(R,\mathfrak{F},\Phi)$ of~Example \ref{gzwovierI}. 
Then $\Theta := \{\QQ_{\ge 0}\}$ extends $\Phi$, 
the associated toric variety is $\PP_{5}$, 
and there is a diagram 
$$
\xymatrix{
{\rq{X} = \Spec(R) \setminus \{0\}} \ar[r] \ar[d] & {\KK^{6} \setminus \{0\}}
\ar[d] \\
{X} \ar[r] & {\PP_{5}}
}
$$
Note that the resulting embedding of $X$ is just
the Pl\"ucker Embedding of the Grassmannian $G(2,4) = X$
into the projective space $\PP_{5}$.
\end{example}

We present now an intrinsic 
characterization of the $\mathfrak{F}$-bunch
$\Phi$ and the stratification of $X$.
For this, we first observe that ``most'' of the variety 
$X$ is contained in the ``big'' orbits of the ambient 
toric variety $Z$.
Denote by $W_{Z} \subset Z$ and $W_{\b{Z}} \subset \b{Z}$ 
the unions of all at most one-codimensional orbits of 
the respective big tori.

\begin{lemma}
\label{subsetWX}
Set $W_{X} := X \cap W_{Z}$
and $W_{\b{X}} := \b{X} \cap W_{\b{Z}}$.
Then the respective complements of 
these sets in $X$ and $\b{X}$
are of codimension at least two.
Moreover, we have 
$W_{\b{X}} = p_{X}^{-1}(W_{X})$,
and $T$ acts freely on this set.
\end{lemma}

\proof
First recall from Proposition~\ref{coxconstrgeo}, 
that the corresponding statements hold for the
sets $W_{Z}$ and $W_{\b{Z}}$.
Thus, $T$ acts freely on $W_{\b{X}}$,
and the commutative diagram of
Proposition~\ref{xintv}
shows $W_{\b{X}} = p_{X}^{-1}(W_{X})$.

It remains to verify the statements on
the codimensions.
For this, we describe $W_{\b{Z}}$
in terms of the coordinates 
$\chi^{e_{i}} = z_{i}$ of $\b{Z} = \KK^{r}$.
It is given by
$$ 
W_{\b{Z}}
\; = \; 
\b{Z} \setminus 
\bigcup_{i < j} V(\b{Z}; z_{i}, z_{j}).
$$
According to the fact that $f_{i} = z_{i} \vert_{\b{X}}$ holds, 
we obtain 
$$
W_{\b{X}}
\; = \; 
\b{X} \setminus \bigcup_{i < j} V(\b{X}; f_{i}, f_{j}).
$$

Since the functions $f_{i}$ are pairwise coprime,
we can conclude that the complement 
$\b{X} \setminus p_{X}^{-1}(W_{X})$ 
is of codimension at least two in $\b{X}$.
Using semicontinuity of the fibre dimension,
one obtains the analogous statement for 
$X \setminus W_{X}$ in $X$.
\endproof 

This lemma allows to define a pull back 
$\imath^{*} \colon \WDiv^{T_{Z}}(Z) \to \WDiv(X)$ 
for invariant Weil divisors
(note that by Proposition~\ref{xintv}, $X$ 
meets the big torus of $Z$): 
first restrict to the smooth set 
$W_{Z} \subset Z$, 
then pullback Cartier divisors to $W_{X}$, 
and, finally, apply the unique extension
$\CDiv(W_{X}) \to \WDiv(X)$.

Similarly, we can define a pullback of Weil
divisors
$p_{X}^{*} \colon \WDiv(X) \to \WDiv(\rq{X})$.
As a quotient space of the locally factorial variety
$W_{\b{X}}$ by a free torus action, 
$W_{X}$ is again locally factorial.
Hence we have a well defined pullback
$\WDiv(W_{X}) \to \WDiv(W_{\b{X}})$.
Again, the observation on the 
codimensions of the complements 
gives the desired extension.

\begin{proposition}
\label{intrinsicPhi}
Let $D^{i}_{Z}$ be the invariant prime divisors 
of the ambient toric variety $Z$,
numbered in such a way that 
$p_{Z}^{*}(D^{i}_{Z}) = \div(\chi^{e_{i}})$, 
and set $D^{i}_{X} := \imath^{*}(D^{i}_{Z})$. 
\begin{enumerate}
\item The pullbacks with respect to 
$p_{X} \colon \rq{X} \to X$ satisfy 
$p_{X}^{*}(D_{X}^{i}) = \div(f_{i})$.
\item A face $\gamma_{0} \preceq \gamma$ belongs to
$\rel(\Phi)$ if and only if there is a point $x \in X$
with
$$
\label{intrischar}
x \; \in \; \bigcap_{e_{i} \not\in \gamma_{0}} D_{X}^{i}
\; \setminus \; \bigcup_{e_{i} \in \gamma_{0}}  D_{X}^{i}.
$$
\item For each $\gamma \in \rel(\Phi)$, the stratum 
$X(\gamma_{0}) \subset X$ consists precisely of the points 
$x \in X$ satisfying the condition of~(ii).
\end{enumerate} 
\end{proposition}

\proof The first statement follows immediately 
from functoriality of the pullback and from
$$
\div(f_{i}) = \rq{\imath}^{*}(\div(\chi^{e_{i}})),
\qquad
\div(\chi^{e_{i}}) = p_{Z}^{*}(D_{Z}^{i}).
$$
The second and the third statement follow from
Proposition~\ref{xintv} and the corresponding
description of the toric orbit stratification
in terms of the $D_{Z}^{i}$, see e.g.~\cite{Co}.
\endproof

As mentioned, the ambient toric varieties need 
not be unique. 
One may overcome this problem by admitting
ambient toric varieties that do not come 
from maximal projectable fans.
Consider the following collection of cones
in the lattice $N := \ker(Q)^{\vee}$:
$$ 
\Delta(R,\mathfrak{F},\Phi) 
\; := \; 
\{\sigma \subset N_{\QQ}; \; \sigma \preceq P(\gamma_{0}^{*})
\text{ for some }
\gamma_{0} \in \rel(\Phi)\}.
$$
This is a fan.
Let $Z(R,\mathfrak{F},\Phi)$ denote the corresponding
open toric variety.
Then, by construction, we have:

\begin{proposition}\label{minimalambient}
$Z(R,\mathfrak{F},\Phi)$ is an open toric subvariety
of any ambient toric variety of $X = X(R,\mathfrak{F},\Phi)$, 
it contains $X$ as a closed subvariety,
and every closed orbit of the big torus of
$Z(R,\mathfrak{F},\Phi)$ meets $X$.
\end{proposition}

\section{Geometric characterization}

In this section we characterize the class
of varieties arising from bunched rings
in terms of their geometric properties.
As we will show in Theorem~\ref{firstresult},
among the varieties with finitely generated 
total coordinate ring, one obtains precisely 
those that are in a certain sense maximal with 
the property of being embeddable into a toric 
variety.

We first recall the precise definitions.
Following W\l odarczyk~\cite{Wl}, 
we say that $X$ is an {\em $A_{2}$-variety}, 
if any two points $x,x' \in X$ admit a common 
affine neighbourhood in $X$.
According to~\cite[Theorem~A]{Wl}, a normal variety 
is $A_{2}$ if and only if it admits a closed embedding 
into a toric variety.

\begin{definition}
We say that an $A_{2}$-variety $X$ is {\em $A_{2}$-maximal},
if it does not admit an open embedding $X \subsetneq X'$ into 
an $A_{2}$-variety $X'$ such that $X' \setminus X$ is of 
codimension at least two in $X'$.
\end{definition}

Next, we recall
the {\em total coordinate ring\/} 
studied in~\cite{ElKuWa}: 
let $X$ be a normal variety
with $\mathcal{O}^{*}(X) = \KK^{*}$ and finitely generated 
free divisor class group $\Cl(X)$.
Choose a subgroup $K \subset \WDiv(X)$
of the group of Weil divisors such that
the canonical map $K \to \Cl(X)$ is an
isomorphism.
Then this gives a total coordinate ring
(depending only up to isomorphism on the 
choice made):
$$ 
\mathcal{R}(X)
\; := \;
\bigoplus_{D \in K} 
  \Gamma(X,\mathcal{O}(D)). 
$$

\begin{theorem}
\label{firstresult}
Let $(R, \mathfrak{F}, \Phi)$ be a bunched ring, 
graded by the lattice $K$. 
Then the associated $X := X(R, \mathfrak{F}, \Phi)$
is a normal $A_{2}$-maximal variety, and it satisfies
$$
\mathcal{O}^{*}(X) = \KK^{*},
\qquad
\Cl(X) \cong K, 
\qquad
\mathcal{R}(X) \cong R,
\qquad
\dim(X) = \dim(R) - \rank(K).
$$
Conversely, every normal $A_{2}$-maximal variety
$X$ with $\mathcal{O}^{*}(X)=\KK^{*}$, 
finitely generated free divisor class group $\Cl(X)$,
and finitely generated total coordinate ring $\mathcal{R}(X)$
arises from a bunched ring.
\end{theorem}

An important ingredient of the proof
is factoriality of the total coordinate
ring. 
This holds in fact for any noetherian 
scheme with finitely generated divisor 
class group, see~\cite[Corollary~1.2]{ElKuWa}.
In our situation, it can also can be deduced 
from results of~\cite{BeHa1}:

\begin{proposition}
\label{factorial}
Let $X$ be a normal variety with
$\mathcal{O}^{*}(X) = \KK^{*}$
and finitely generated free divisor 
class group. Then the following holds.
\begin{enumerate}
\item The total coordinate ring $\mathcal{R}(X)$ 
of $X$ is a factorial ring satisfying 
$\mathcal{R}(X)^* = \KK^*$, 
and its grading is faithful.
\item If $\mathcal{R}(X)$ is finitely generated, 
then it admits a system of homogeneous 
pairwise nonassociated nonzero prime generators.
\end{enumerate}
\end{proposition}

\proof
The open set $X' \subset X$ of smooth 
points of $X$ satisfies $\Pic(X') = \Cl(X)$,
and has homogeneous coordinate ring 
$\mathcal{A}(X') = \mathcal{R}(X)$
in the sense of~\cite{BeHa1}.
Thus,~\cite[Remark~8.8 and Prop.~8.4]{BeHa1} 
yield that $\mathcal{R}(X)$ is factorial.
Moreover, $\mathcal{R}(X)$ is faithfully 
graded, because any divisor class is 
represented by a difference of two
effective divisors.

In order to see $\mathcal{R}(X)^* = \KK^*$,
consider an element $f \in \mathcal{R}(X)^*$
and its inverse $g \in \mathcal{R}(X)^*$. 
Since $fg = 1$ is homogeneous, also $f$ 
and $g$ are so, having Weil divisors
$D$ and $E$ as their respective degrees. 
We obtain $D+E=0$, and considering the 
zero set of $fg$ gives
$$
\emptyset
 = 
Z(fg)
 = 
\Supp(\div(fg)+D+E)
 = 
\Supp(\div(f)+D) \cup \Supp(\div(g)+E).
$$

This implies $D = -\div(f)$, and $E = -\div(g)$.
In other words, $D$ and $E$ are principal.
Since the grading group of $\mathcal{R}(X)$
maps isomorphically onto $\Cl(X)$, we can 
conclude $D = E = 0$. Thus, we obtain
$f,g \in \mathcal{O}^*(X) = \KK^*$,
and the first assertion is proved.

As to the second assertion, 
let $\{g_1, \ldots, g_s\}$ be any system 
of nonzero homogeneous generators of 
$\mathcal{R}(X)$.
Since $\mathcal{R}^*(X) = \KK^*$ holds, 
we may assume that no $g_{i}$ is a 
unit in $\mathcal{R} (X)$. 
Decompose each $g_{i}$ 
into prime factors, and fix a 
system $\{f_1, \ldots, f_r\}$ 
of pairwise non associated 
representatives for the prime factors 
occuring in $g_1,\ldots, g_s$. 
Then the $f_j$ are homogeneous, and
they generate $\mathcal{R}(X)$ as well.
\endproof

Now we turn to the proof of the theorem. 
It is split into two parts.
In the first part, we prove in particular
the second assertion of the theorem, 
and then, in the second part, the first 
assertion is verified.

\proof[Proof of Theorem~\ref{firstresult}, Part~1]
Let $X$ be a normal $A_{2}$-variety
$X$ with $\mathcal{O}^{*}(X)=\KK^{*}$, 
finitely generated free class divisor group,
and finitely generated total coordinate 
ring $R := \mathcal{R}(X)$.

\medskip
\noindent
{\em Assertion~\ref{firstresult}~a)}
\enspace
For any system $\mathfrak{F} \subset R$ 
of nonzero homogeneous
pairwise nonassociated prime generators,
there is an $\mathfrak{F}$-bunch
$\Phi$ and an open embedding
$X \subset X(R,\mathfrak{F},\Phi)$
such that the complement is of codimension
at least two.

\medskip

We prove this assertion.
Choose a sublattice $K \subset \WDiv(X)$ 
of the group of Weil divisors 
such that the canonical map $K \to \Cl(X)$ 
is an isomorphism.
Consider the associated sheaf of $K$-graded 
$\mathcal{O}_{X}$-algebras
$$ 
\mathcal{R} := \bigoplus_{D \in K}
\mathcal{O}_{X}(D).
$$

According to Proposition~\ref{factorial}~(i),
the ring $R = \mathcal{R}(X)$ of global 
sections is a unique factorization domain,
and it satisfies $R^* = \KK^*$.
Moreover, by assumption, $R$ is finitely 
generated.
Let $\b{X} := \Spec(R)$,
and consider the relative spectrum 
$\rq{X} := \Spec_X(\mathcal{R})$
with the canonical affine map 
$p_X \colon \rq{X} \to X$.

We show that the relative spectrum $\rq{X}$ 
is an open subvariety of $\b{X}$.
For this, note first that $\rq{X}$ is 
covered by open affine varieties 
of the form $p_X^{-1}(X_f)$, 
where $f \in R$ is homogeneous and 
$X_f$ is defined as 
$$
X_f 
\; := \;
X \setminus \Supp(\div(f)+D).
$$ 

This is due to the fact that every affine 
open $U \subset X$ is of the form 
$U = X_f$, because $X \setminus U$
is the support of some effective Weil 
divisor,
which in turn occurs as the zero set of 
a section 
$f \in \Gamma(X,\mathcal{O}_D) \subset R$ 
of some 
$D \in K \cong \Cl(X)$.

Next note that the canonical morphism 
$\rq{X} \to \b{X}$
maps $p_X^{-1}(X_f)$ isomorphically
onto $\b{X}_f$. 
This follows from the fact that, if
$f \in R$ is homogeneous, say of degree
$D \in K$, then we have 
\begin{equation*}
\mathcal{R}(X_f) 
\; = \; 
R_f.
\end{equation*}

In conclusion, we obtain that 
$\rq{X} \to \b{X}$ is a local 
isomorphism and hence an open embedding
(use e.g. Zariski's Main Theorem).
In particular, we see that
$\rq{X}$ is quasiaffine and that 
$p_X^{-1}(X_f)$ equals
$\b{X}_f$.

We shall make use of the actions of the torus
$T := \Spec(\KK[K])$ on $\b{X}$ and $\rq{X}$
defined by the $K$-grading of the 
$\mathcal{O}_{X}$-algebra $\mathcal{R}$.
First note that the open embedding 
$\rq{X} \to \b{X}$ is $T$-equivariant, and that 
the canonical morphism
$p_{X} \colon \rq{X} \to X$ 
is a good quotient for the $T$-action.

Proposition~\ref{factorial}~(ii)
enables us to fix a system 
$\mathfrak{F} = (f_{1}, \ldots, f_{r})$ 
of pairwise nonassociated homogeneous prime elements that 
generate $R = \mathcal{R}(X)$.
Moreover, let $T$ act diagonally on $\KK^{r}$ such that
the coordinate $z_{i}$ has the same degree as 
the function $f_{i}$.
This gives rise to a $T$-equivariant 
closed embedding
$$
\b{X} \to \KK^{r},
\qquad 
x \mapsto (f_{1}(x), \ldots, f_{r}(x)).
$$

The first major step is to show that the 
good quotient $p_X \colon \rq{X} \to X$ by $T$ 
admits a toric extension,
that means that there is an open toric subvariety 
$\rq{Z} \subset \KK^{r}$ admitting a good quotient 
$p_{Z} \colon \rq{Z} \to Z$
for the $T$-action such that 
$\rq{X} = \rq{Z} \cap \b{X}$ holds.
Note that then restricting $p_{Z}$ to $\rq{X}$
gives back the initial map 
$p_{X} \colon \rq{X} \to X$.

Let $\TT^r \subset \KK^r$ denote the standard $r$-torus.
The candidate for our $\rq{Z}$ is the minimal 
$\TT^{r}$-invariant open subset
of $\KK^r$, which contains $\rq{X}$:
$$ 
\rq{Z} 
:= 
\{z \in \KK^r; \; \b{\TT^r \mal z} \cap \rq{X} \ne \emptyset\}.
$$

Then $\rq{Z}$ satisfies $\rq{X} = \rq{Z} \cap \b{X}$.
Indeed, the inclusion 
``$\subset$'' is clear by definition.
For the reverse inclusion, recall that 
$\rq{X}$ is the union of the subsets
$\b{X}_{f_i}$, and that 
$f_i = z_i \vert_{\b{X}}$ holds.
Thus, we have
$$
\rq{Z} \cap \b{X}
\; \subset \;
\bigcup_{i=1}^r \KK^r_{z_i}
\cap \b{X} 
\; = \; 
\rq{X}.
$$

Our task is to show that $\rq{Z}$ admits a good 
quotient.
To verify this, it suffices to show that 
any two points $z,z' \in \rq{Z}$ 
admit a common $T$-invariant 
affine neighbourhood in $\rq{Z}$, 
see~\cite[Theorem~C]{BBSw}.
The latter property will be established
by refining an argument used in~\cite{Ha1}.
Choose
$$
x \in \b{\TT^{r} \mal z} \cap \b{X},
\qquad 
x' \in \b{\TT^{r} \mal z'} \cap \b{X}.
$$ 

Then it is enough to show that there is  
a common $T$-invariant affine open 
neighbourhood
$U \subset \rq{Z}$ of $x$ and $x'$. 
Namely, any such $U$ intersects the $\TT^r$-orbits
of $z$ and $z'$. 
Thus, as an intersection of nonempty open subsets
of $\TT^r$, the set
$$ 
V 
\; := \;  
\{t \in \TT^r; \; t^{-1} \mal z \in U\}
\cap
\{t \in \TT^r; \; t^{-1} \mal z' \in U\}
$$
is nonempty. Take any element $t_0 \in V$.
Then the translate 
$t_0 \mal U$ contains both points, 
$z$ and $z'$,
and hence it is the desired common $T$-invariant 
affine open neighbourhood 
of $z$ and $z'$ in $\rq{Z}$.

In order to show the existence of a 
common $T$-invariant neighbourhood 
$U$ for the points $x$ and $x'$, we make use of
the $A_2$-property of $X$: we choose a common affine
open neighbourhood $X_{0} \subset X$ of 
the images $p_X(x)$ and $p_X(x')$.

Then the inverse image $p_{X}^{-1}(X_{0})$ 
is affine, and hence its complement in $\b{X}$ 
is the support of an effective Weil divisor.
Thus, by factoriality of $R$, there is a 
$T$-homogeneous $h \in R$ with 
$p_{X}^{-1}(X_{0}) = \b{X}_{h}$.
As the $f_{i}$ generate $R$, we can write 
$h = H(f_{1}, \dots, f_{r})$ 
for some $T$-homogeneous polynomial
$H \in \KK[z_{1}, \ldots, z_{r}]$. 

Consider the good quotient 
$q \colon \KK^{r}_{H} \to \KK^{r}_{H} \quot T$, 
which exists, because $\KK^{r}_{H}$ is an affine
variety. Then we have 
$$ 
q (\b{T \mal x} \cup \b{T \mal x'}) 
\cap 
q(\KK^{r}_{H} \setminus \rq{Z}) 
\; = \;
\emptyset.
$$ 
Take any affine neighbourhood 
$V \subset \KK^{r}_{H} \quot T$ of $q(x)$ and $q(x')$ 
that does not intersect 
$q(\KK^{r}_{H} \setminus \rq{Z})$. 
Then the inverse image $U := q^{-1} (V)$ is 
a $T$-invariant affine neighbourhood
of $x$ and $x'$ in $\rq{Z}$.

So, we obtained the desired toric extension
$p_{Z} \colon \rq{Z} \to Z$ 
of the good quotient 
$p_X \colon \rq{X} \to X$.
Now, consider any 
enlargement
of our toric extension,
that means a toric open
subset $\rq{Z}' \subset \KK^r$ 
admitting a good quotient
$p_{Z'} \colon \rq{Z}' \to Z'$
for the $T$-action such that $\rq{Z}$
is  a $p_{Z'}$-saturated subset 
of $\rq{Z}'$.

Cutting down to $\b{X}$,
we obtain an open inclusion
$\rq{X} \subset \rq{X}'$,
and hence
$X \subset X'$ on the level of 
quotients.
Since $\rq{X}$ and $\rq{X}'$ have 
the same global functions,
the complement 
$\rq{X}' \setminus \rq{X}$ 
is of codimension at least two.
By semicontinuity of the fibre dimension,
the same holds for $X' \setminus X$.

In the above procedure, we may assume 
that $\rq{Z}'$ is $T$-maximal, i.e., cannot
be enlarged itself.
In this situation, we shall construct
an $\mathfrak{F}$-bunch $\Phi$ such
that $X' = X(R,\mathfrak{F},\Phi)$ holds.
For the sake of readable notation, 
we perform a reset:
we replace $X'$ with $X$ and $Z'$ with $Z$ etc..

The crucial point is to show that 
the toric morphism $p_{Z} \colon \rq{Z} \to Z$
is a Cox construction for $Z$.
For this, it is most convenient 
to verify the conditions of 
geometric characterization 
Proposition~\ref{coxconstrgeo}.
The first property is clear by 
construction.
Thus, we are left with checking
Condition~\ref{coxconstrgeo}~(ii).

Let  $V_{X} \subset X$ denote the 
set of smooth points,
and consider the inverse image
$V_{\rq{X}} = p_{X}^{-1}(V_{X})$.
This set has a 
complement of codimension at least two
in $\b{X}$, because, by normality,
the same holds for the subset $V_{X} \subset X$,
and hence we obtain for the global functions 
$$
\mathcal{O}(V_{\rq{X}})
=
\mathcal{R}(V_{X}) 
= 
R
=\mathcal{O}(\b{X}).
$$

Moreover, $V_{\rq{X}}$ meets the big torus
$\TT^r \subset \KK^r$.
It meets also every $\TT^r$-orbit
of codimension one:
These orbits are the sets $B_i$
given by $z_i = 0$
for some $i$ and $z_j \ne 0$
for all $j \ne i$.
Thus, the intersection of 
such an orbit with $\b{X}$
is given by
$$
B_i \cap \b{X} 
\; = \;
\{x \in \b{X}; \; 
f_i(x) = 0, \; 
f_j(x) \ne 0 \text{ for all } j \ne i\}. 
$$
Since the $f_k$ are relatively 
coprime, these intersections are 
nonempty, and hence of codimension one 
in $\b{X}$.
Since $V_{\rq{X}} \subset \b{X}$ 
has a complement of codimension at least two,
every $B_i \cap V_{\rq{X}}$ is nonempty.

Eventually, we claim that $T$ acts freely on the 
set $V_{\rq{X}}$.
This is most easily seen in terms of homogeneous 
functions: 
on $V_{X}$ every divisor $D \in K$ 
is Cartier.
Hence, locally on $V_{\rq{X}}$, there 
exists an invertible
homogeneous function in every degree 
$D \in K$.
This implies freeness of the $T$-action
on $V_{\rq{X}}$.

In conclusion, for the set 
$W_{\b{Z}} \subset \b{Z} := \KK^r$
consisting of the at most one codimensional 
$\TT^r$-orbits,
the preceding observations imply that 
$T$ acts freely on $W_{\b{Z}}$, 
that $W_{\b{Z}} \subset \rq{Z}$ holds, 
and that every fibre
of $p_Z^{-1}(p_Z(z))$, where $z \in W_{\b{Z}}$,
consists of a single (free) $T$-orbit.
Thus we verified Condition~\ref{coxconstrgeo}~(ii).

Knowing that $\rq{Z}$ is $T$-maximal
and that $p_Z \colon \rq{Z} \to Z$ is
a Cox construction,
we may apply~\cite[Lemma~7.8]{BeHa2}.
It tells us that $p_Z \colon \rq{Z} \to Z$
arises from a bunch
$\Theta$ in the projected cone
$(\ZZ^{r} \topto{Q} K, \QQ^{r}_{\ge 0})$, 
where $Q$ maps each $e_{i}$ to the degree of the 
coordinate $z_{i}$.

Now we are ready to define the $\mathfrak{F}$-bunch
$\Phi$.
Consider the collection $\Omega \subset \rel(\Theta)$
of all $\mathfrak{F}$-faces.
Let $\Phi$ be the collection of all minimal
cones among the images $Q(\gamma_0)$, 
where $\gamma_0 \in \Omega$.
Then $\Phi$ is an $\mathfrak{F}$-bunch,
extended by $\Theta$.
Using Proposition~\ref{xintv},
it is straightforward to check that 
$X = X(R,\mathfrak{F},\Phi)$ holds.
\endproof

\proof[Proof of Theorem~\ref{firstresult}, Part~2]
We verify the claimed properties for 
the variety $X$ arising from a bunched ring 
$(R,\mathfrak{F},\Phi)$.
We work in terms of $\b{X} := \Spec(R)$,
the torus $T := \Spec(\KK[K])$, the subset
$\rq{X} \subset \b{X}$ and the good quotient
$p_{X} \colon \rq{X} \to X$ as introduced in
the construction of $X$.

First of all note that, as a good quotient space 
of the normal variety $\rq{X}$, the variety $X$ is 
normal. 
Moreover, choosing any ambient toric variety 
$X \subset Z$, we obtain that $X$ inherits the 
$A_{2}$-property.

To determine the coordinate ring $\mathcal{R}(X)$,
consider the subset $W_{X} = X \cap W_{Z}$ as 
discussed in Lemma~\ref{subsetWX}.
Then $p_{X}^{-1}(W_{X})$ is locally factorial, and
$T$ acts freely on this set.
Thus, the results of~\cite{BeHa1} give us the
Picard group $\Pic(W_{X})$ and the homogeneous coordinate ring 
$\mathcal{A}(W_{X})$ as defined in~\cite{BeHa1}: 
$$
\Pic(W_{X}) \cong K,
\qquad
\mathcal{A}(W_{X}) \cong \mathcal{O}(p_{X}^{-1}(W_{X})) \cong 
\mathcal{O}(\b{X}) = R. 
$$

Since $X \setminus W_{X}$ is of codimension at least two,
we obtain that $\Pic(W_{X})$ is isomorphic to $\Cl(X)$,
and that $\mathcal{A}(W_{X})$ is isomorphic to 
$\mathcal{R}(X)$.
In fact, we shall later provide an explicit
isomorphism $K \cong \Cl(X)$, making the present
proof independent from the reference~\cite{BeHa1}. 

It remains to show that $X$ is $A_{2}$-maximal.
So, assume that there is an open embedding 
$X \subset X'$
with an $A_2$-variety $X'$ 
such that the complement $X' \setminus X$ is of 
codimension at least two in $X'$.
We may assume that $X'$ is normal; otherwise,  
take its normalization instead.  
Then the inclusion $X \subset X'$ gives rise to commutative
diagrams of isomorphisms
$$
\xymatrix{
& K \ar[dl] \ar[dr] & \\
{\Cl(X)} & & {\Cl(X')} \ar[ll]
}
\qquad
\qquad
\xymatrix{
& R \ar[dl] \ar[dr] & \\
{\mathcal{R}(X)} & & {\mathcal{R}(X') \ar[ll]}
}
$$

According to Proposition~\ref{factorial}~(ii)
and Assertion~\ref{firstresult}~a), we 
find an $\mathfrak{F}$-bunch $\Phi'$
defining a commutative diagram 
where the lower row is an open embedding
having a complement of codimension at
least two:
$$
\xymatrix{
& X \ar[dl] \ar[dr] & \\
X' \ar[rr] & & X(R,\mathfrak{F},\Phi')
}
$$

Hence, we may even assume 
that $X'$ equals $X(R,\mathfrak{F},\Phi')$.
Now choose extending bunches
$\Theta$ and $\Theta'$ for $\Phi$ 
and $\Phi'$.
Consider the corresponding 
ambient toric varieties
$Z$ and $Z'$.
Using Lemma~\ref{subsetWX}
(and its notation), 
we obtain the following 
commutative diagram
$$
\xymatrix@!0{
& 
{W_{\b{Z}}} \ar@{=}[rrr]  \ar'[d][dd]
& & &
W_{\b{Z}'} \ar[dd] 
\\
p_X^{-1}(W_X) \ar@{=}[rrr]  \ar[dd] \ar[ur]
& & &
p_{X'}^{-1}(W_{X'})  \ar[dd] \ar[ur]
\\
&
W_Z \ar@{=}'[rr][rrr] 
& & &
W_{Z'}
\\
W_X \ar[rrr] \ar[ur] 
& & &
W_{X'} \ar[ur]
}
$$

In particular, we can conclude that
the induced morphism $W_X \to W_{X'}$
is the identity map.
Consequently, the divisors 
$D^{i}_{X}$ and $D^{i}_{X'}$
considered in Proposition~\ref{intrinsicPhi} 
satisfy
$$D^{i}_{X'} \vert_{X} = D^{i}_{X}.$$
Now we may use the intrinsic characterization 
of relevant faces provided in 
Proposition~\ref{intrinsicPhi}~(ii).
It yields
$\rel(\Phi) \subset \rel(\Phi')$.
By the maximality properties of $\mathfrak{F}$-bunches, we obtain
$\Phi = \Phi'$.
This eventually gives $X = X'$.
\endproof

We conclude this section with two immediate 
applications concerning toric varieties.
For $\QQ$-factorial projective varieties 
$X$, the first statement was also observed 
in~\cite[Cor.~2.10]{KeHu}.

\begin{corollary}
Let $X$ be a normal, $A_2$-maximal variety 
with $\mathcal{O}(X) = \KK$, 
finitely generated free divisor class group
and finitely generated total coordinate ring
$\mathcal{R}(X)$. Then the following
statements are equivalent:
\begin{enumerate}
\item $X$ is a toric variety.
\item $\mathcal{R}(X)$ is a polynomial ring.
\end{enumerate} 
\end{corollary}

\proof 
By~\cite{Co}, the total coordinate ring of 
a toric variety is a polynomial ring.
Conversely, let us assume that the total 
coordinate ring $R := \mathcal{R}(X)$ 
of $X$ is a polynomial ring. Then the set 
$\mathfrak{F} = \{T_1,\ldots, T_r\}$
consisting of the indeterminates
is a system of generators for $R$.

Since we assumed $\mathcal{O}(X)=\KK$,
we have $R_0 = \KK$, and 
thus the torus action on $\KK^r = \Spec(R)$ 
induced by the $\Cl(X)$-grading is 
linearizable, see~\cite{KaRu}.
In other words, we may assume that 
the $T_i$ are homogeneous.
Assertion~\ref{firstresult}~a) 
then says that $X$ arises from a bunched
polynomial ring $(R,\mathfrak{F},\Phi)$,
and hence is a toric variety.
\endproof

Note that $A_2$-maximality is essential 
in this statement: the projective plane 
with four different points removed has
$\KK[T_1,T_2,T_3]$ as its total coordinate 
ring, but it can never be made into a 
toric variety.

Here comes a characterization of
$A_2$-maximality for a toric variety
in terms of its fan 
(inside the category of toric varieties, 
the statement is obvious).

\begin{corollary}
\label{toric2a2max}
Let $\Delta$ be a fan in a lattice $N$,
such that $N$ is generated as a lattice
by the primitive vectors of the rays of
$\Delta$, and $\Delta$ cannot
be enlarged without adding new rays.
Then the corresponding toric variety $X$
is $A_2$-maximal.
\end{corollary}

\proof
The assumptions on the fan $\Delta$ guarantee
that $X$ arises from a bunched (polynomial)
ring, compare~\cite[Theorem~7.3~(ii)]{BeHa2}.
Hence, Theorem~\ref{firstresult} applies.
\endproof

\section{Divisor classes and singularities}

In this section, we take a closer look at
the divisors of the variety $X$ associated 
to a given bunched ring $(R,\mathfrak{F},\Phi)$.
We provide an explicit isomorphism 
$K \to \Cl(X)$ 
from the grading lattice onto the divisor 
class group,
and we investigate the behaviour of divisors
in terms of this isomorphism and the stratification 
of $X$.
As an application, we study the singularities 
of $X$.

For the definition of the isomorphism
$K \to \Cl(X)$, recall that $X$ is the 
good quotient of an open subset 
$\rq{X}$ of $\b{X} := \Spec(R)$ by the action
of $T := \Spec(\KK[K])$.
Moreover, we showed in Lemma~\ref{subsetWX}
that there is an open set $W_{X} \subset X$ 
such that $T$ acts freely on 
the inverse image 
$W_{\b{X}} := p_{X}^{-1}(W_{X})$
under the quotient map 
$p_{X} \colon \rq{X} \to X$,
and both sets have small complement in $X$
and $\b{X}$ respectively.

\begin{lemma}\label{isodescr}
For every $w \in K$, choose a $w$-homogeneous 
function $h_{w} \in \KK(\rq{X})^{*}$
and cover the set $W_{\b{X}}$
by $T$-invariant open sets $W_{j}$ admitting
$w$-homogeneous functions
$h_{j} \in \mathcal{O}^{*}(W_{j})$. 
Then there is a (welldefined) isomorphism
$$
\b{\mathfrak{D}}_{X} \colon K \to \Cl(X),
\qquad 
w \mapsto \b{\mathfrak{D}}_{X}(w),
\text{ where }
\b{\mathfrak{D}}_{X}(w)\vert_{p_{X}(W_{j})}
=
[\div(h_{w}/h_{j})].
$$
Moreover, in the notation of Proposition \ref{intrinsicPhi}, we have
$\b{\mathfrak{D}}_{X} (w_{i}) = [D^{i}_{X}]$, where $w_{i} = \deg(f_{i})$
for $f_i \in \mathfrak{F}$.
\end{lemma}

\proof
First note that the open subset 
$W_{\b{X}} \subset \Spec(R)$ is locally factorial,
and that $\Pic(W_{\b{X}})$ is trivial.
Both statements hold, because $R$ is a factorial ring.
{}From the first one we infer that, as a free 
$T$-quotient of the locally factorial $W_{\b{X}}$, 
the variety $W_{X}$ is again locally factorial. 

To proceed, we work in terms of $T$-linearized 
line bundles.
Any $w \in K$ corresponds to a character 
$\chi^{w} \colon T \to \KK^{*}$, and 
this character defines a linearization of the 
trivial bundle, namely
$$ 
t \mal (x,z) 
\; := \; 
(t \mal x, \chi^{w}(t)^{-1} z).
$$

Thus, every $w \in K$ defines an element
in the group  $\Pic_{T}(W_{\b{X}})$
of equivariant isomorphy classes of $T$-linearized
line bundles on $W_{\b{X}}$. 
Since $\Pic(W_{\b{X}}) = 0$ and 
$\mathcal{O}^*(W) = \KK^*$ hold,
this assignment even defines an isomorphism, 
use~\cite[Lemma~2.2]{KKV}:
$$\alpha \colon K \to \Pic_{T}(W_{\b{X}}).$$

Moreover, according to~\cite[Proposition~4.2]{KKV}
one obtains a pullback homomorphism
$\Pic(W_{X}) \to \Pic_{T}(W_{\b{X}})$
by endowing each pullback bundle with
the trivial linearization.
Denoting by $\beta$ the inverse of this pullback
isomorphism, we obtain an isomorphism
$$
\beta \circ \alpha \colon K \to  \Pic(W_{X}). 
$$

Since the complement $X \setminus W_{X}$ is of codimension 
at least two in $X$, the isomorphism $\beta \circ \alpha$
induces an isomorphism $K \to \Cl(X)$.
One directly verifies that this is the map $\b{\mathfrak{D}}_{X}$
defined in the assertion.

In order to verify the last statement, recall
from Proposition~\ref{intrinsicPhi}(i)
that $p_{X}^{*}(D_{X}^{i})$ equals $\div(f_{i})$.
Thus, locally, we have
$$ 
p_X^*(D_{X}^{i})
\; = \; 
\div(f_i)
\; = \; 
\div(f_i/h_j)
\; = \; 
p_X^*(\div(f_i/h_j)),
$$
where the $h_j$ are functions of degree $w_i$
as in the definition of $\b{\mathfrak{D}}_X$.
Since pullback is injective on the level of  
divisors, the assertion follows. 
\endproof

Note that the isomorphism 
$\b{\mathfrak{D}}_{X} \colon K \to \Cl(X)$
does not depend on the choices of $h_{w}$ and 
the $h_{j}$ made in its definition.

Moreover, $\b{\mathfrak{D}}_X$ 
is compatible with any ambient toric 
variety $Z$ of $X$:
In terms of the projected cone
$(E \topto{Q} K, \gamma)$
associated to $\mathfrak{F}$,
we have a (well defined) isomorphism
compare~\cite[Section~10]{BeHa2}:
$$
\b{\mathfrak{D}}_{Z} \colon K \to \Cl(Z),
\qquad
w \mapsto \sum_{i=1}^r 
\bangle{\rq{w}, e_{i}^*} 
[ D_i^Z ],
\qquad 
\text{where } Q(\rq{w}) = w.
$$

\begin{proposition}
\label{toricCl2XCl}
The closed embedding $\imath \colon X \to Z$
induces a commutative diagram of isomorphisms:
$$
\xymatrix{
K \ar@{=}[r] \ar[d]_{\b{\mathfrak{D}}_{Z}} 
&
K \ar[d]^{\b{\mathfrak{D}}_{X}}
\\
{\Cl(Z)} \ar[r]_{\imath^{*}} 
& 
\Cl(X)
}
$$
\end{proposition}

\proof 
As noted earlier, Lemma~\ref{subsetWX} ensures that 
the pullback map 
$\imath^{*} \colon \Cl(Z) \to \Cl(X)$ 
is in fact well defined.
The rest is an immediate consequence of the last statement 
of Proposition~\ref{isodescr}. 
\endproof

In terms of the map $\b{\mathfrak{D}}_{X}$,
we can say when a given Weil divisor 
$D \in \WDiv(X)$ is $\QQ$-Cartier or even Cartier 
at a point $x \in X$:

\begin{proposition}
\label{divisors}
Let $\gamma_{0} \in \rel(\Phi)$, and let $x \in X(\gamma_{0})$
be a point in the corresponding stratum.
\begin{enumerate}
\item A class $\b{\mathfrak{D}}_{X}(w)$ is $\QQ$-Cartier at $x$ if and only if 
      $w \in \lin(Q(\gamma_{0}))$ holds.
\item A class $\b{\mathfrak{D}}_{X}(w)$ is Cartier at $x$ if and only if 
      $w \in Q(\lin(\gamma_{0}) \cap E)$ holds.
\end{enumerate}
\end{proposition}

\proof It suffices to verify the second assertion.
For this let $D$ represent the class 
$\b{\mathfrak{D}}_{X}(w)$, and note that $D$ is Cartier at $x$ 
if and only if the stalk $\mathcal{O}(D)_{x}$ is generated 
by a single element $g_{x} \in \mathcal{O}(D)_{x}$.
This is equivalent to saying that
the fibre $Y := p_{X}^{-1}(x)$ admits
an invertible homogeneous element
$g \in \mathcal{O}(Y)_{w}$. 

Let $Z$ be any ambient toric variety of $X$,
and let $\rq{Z}$ be as in Proposition~\ref{xintv}.
Then the fibre $Y$ is closed in the toric 
variety $\rq{Z}$, and the affine toric chart
$\rq{Z}_{\gamma_{0}^{*}} \subset \rq{Z}$
corresponding to the cone $\gamma_0^*$ 
is the smallest $T_{\rq{Z}}$-invariant
affine open subset containing $Y$.

Consequently, there exists an invertible
element $g \in \mathcal{O}(Y)_{w}$
if and only if the toric neighbourhood 
$\rq{Z}_{\gamma_{0}^{*}} \subset \rq{Z}$ 
admits an invertible character function
$\chi^{u} \in \mathcal{O}(\rq{Z}_{\gamma_{0}^{*}})_{w}$.
But the latter is equivalent to the existence of 
a linear form
$u \in \lin(\gamma_{0}) \cap E$ with $Q(u) = w$.
\endproof

As a first application, we give basic statements
on the singularities of 
$X = X(R,\mathfrak{F},\Phi)$.
In particular, we will see that the stratification 
by relevant faces is ``equisingular'' in a certain
raw sense.
We say that a point $y$ of a variety $Y$ 
is {\em $\QQ$-factorial\/}
if it is normal and, near $y$, every 
closed subset of codimension one
is the set of zeroes of 
a regular function.

\begin{corollary}
\label{singularities}
Let $\gamma_{0} \in \rel(\Phi)$, and consider a
point $x \in X(\gamma_{0})$.
\begin{enumerate}
\item The point $x$ is factorial if and only if 
$Q$ maps $\lin(\gamma_{0}) \cap E$ onto $K$.
\item The point $x$ is $\QQ$-factorial if and only if 
 $Q(\gamma_{0})$ is of full dimension.
\end{enumerate}
\end{corollary}

\proof
This is an immediate consequence of two things: firstly
a point $x \in X$ is factorial ($\QQ$-factorial) if and 
only if every Weil divisor is Cartier ($\QQ$-Cartier)
at $x$; secondly, the characterization of the latter 
properties in terms of relevant faces performed in 
Proposition~\ref{divisors}. \endproof

Let us have a look at the case
that the open set $\rq{X} \subset \b{X}$ occuring
in the construction of $X$ is smooth.
This situation will occur in all the examples 
studied later.
First we note a general statement.

\begin{proposition}
Suppose that $\rq{X}$ is smooth. Then $X$ has at most
rational singularities. In particular, $X$ is 
Cohen-Macaulay.
\end{proposition}

\proof 
This follows from general theorems saying 
that quotients of smooth varieties by 
reductive groups have at most rational
singularities, see~\cite{Bo} and~\cite{HR}.
\endproof

In general, there is no hope to decide in terms of 
combinatorial data whether or not a given $x \in X$ is 
a smooth point. 
But under the additional assumption that $\rq{X}$ is smooth,
we have the following.

\begin{proposition}\label{smoothchar}
Suppose that $\rq{X}$ is smooth,
let $\gamma_{0} \in \rel(\Phi)$,
and $x \in X(\gamma_{0})$.
Then $x$ is a smooth point if and only if
$Q$ maps $\lin(\gamma_{0}) \cap E$ onto $K$.
\end{proposition}

\proof
The ``only if'' part is clear by
Corollary~\ref{singularities}.
So, let $Q$ map $\lin(\gamma_{0}) \cap E$ onto $K$.
Then the quotient presentation 
$p_Z \colon \rq{Z} \to Z$ 
of an ambient toric variety,
is a principal 
bundle over the affine chart 
$Z_{\gamma_{0}} := p_Z(\rq{Z}_{\gamma_0^*})$ 
of $Z$ corresponding to $\gamma_{0}$. 
The latter holds as well for $\rq{X} \to X$ over 
$X_{\gamma_{0}} = X \cap Z_{\gamma_{0}}$.
Since we have $X(\gamma_{0}) \subset X_{\gamma_{0}}$,
we see that $x$ is smooth. 
\endproof

According to Proposition~\ref{divisors}
the behaviour of singularities is reflected 
in the ambient toric varieties
studied in Section~\ref{toricambient}.

\begin{proposition}
\label{singulemb}
The embedding $X \subset Z$ into an ambient 
toric variety $Z$ has the following
property: 
a point $x \in X$ is a factorial ($\QQ$-factorial)
point of $X$ if and only if it is a smooth 
($\QQ$-factorial) point of $Z$.
\endproof
\end{proposition}

Recall from~\cite{Co}, that in Cox's quotient 
presentation $p_Z \colon \rq{Z} \to Z$ of a toric 
variety $Z$, a fibre $p_Z^{-1}(z)$ over a point 
$z \in Z$ consists of a single closed orbit if and 
only if $z$ is a $\QQ$-factorial point. 
Thus, Corollary~\ref{singularities},
and Propositions~\ref{singulemb} and~\ref{xintv}
fit together to the following statement.

\begin{corollary}
For $X = X(R,\mathfrak{F},\Phi)$,
the following statements are equivalent:
\begin{enumerate}
\item  $p_X \colon \rq{X} \to X$ is a geometric quotient;
\item $X$ is $\QQ$-factorial;
\item $\Phi$ consists of cones of full dimension.
\end{enumerate}
\end{corollary}

As before, if $\rq{X}$ is smooth, then we may
replace ``factorial'' with ``smooth'' in 
Proposition~\ref{singulemb}.
Moreover, if, in this setup, we take the embedding 
$X \subset Z(\Phi)$ discussed in 
Proposition~\ref{minimalambient},
then we arrive at the statement that
$X$ is smooth ($\QQ$-factorial) if and only if $Z(\Phi)$
is so.

\section{Functoriality properties}

In this section, we make the assignment 
from bunched rings to varieties into a
functor, and study basic properties of
this functor. For example, we figure 
out when two bunched rings define 
isomorphic varieties, 
see Corollary~\ref{iso42iso}.

Firstly, we have to introduce the concept
of a morphism of bunched rings.
In doing this, we follow closely the ideas 
of~\cite{BeHa1}.
To any bunched ring $(R,\mathfrak{F},\Phi)$ 
we associate an {\em irrelevant ideal}:
$$
I(R,\mathfrak{F},\Phi)
\; := \; 
\sqrt{
\bangle{f^{u}; \; u \in \gamma_{0}^{\circ}, %
                     \; \gamma_{0} \in \rel(\Phi)}%
},
$$ 
where $f^{u} := f_{1}^{u_{1}} \ldots f_{r}^{u_{r}}$ as
in~\ref{Xhutdef}. Geometrically, the ideal $I(R,\mathfrak{F},\Phi)$ 
is the vanishing ideal of the complement $\b{X} \setminus \rq{X}$, where
$\b{X} := \Spec(R)$ and $\rq{X} \subset \b{X}$ is the open set introduced
in~\ref{Xhutdef}.

Now, let $(R_{1},\mathfrak{F}_{1}, \Phi_{1})$, 
and $(R_{2},\mathfrak{F}_{2}, \Phi_{2})$ be bunched 
rings. 
We set $\b{X}_i := \Spec(R_i)$ etc..
The definition of a morphism is built up in 
two steps. The first one is

\begin{definition}
\label{mograburi}
By a {\em graded homomorphism\/} of the bunched rings
$(R_{2},\mathfrak{F}_{2},\Phi_{2})$ and 
$(R_{1},\mathfrak{F}_{1},\Phi_{1})$ 
we mean a $\KK$-algebra homomorphism 
$\mu \colon R_{2} \to R_{1}$ such that
\begin{enumerate}
\item there is a homomorphism  $\t{\mu} \colon
K_{2} \to K_{1}$ of the respective grading groups with
$$ \mu(R_{2,w}) 
\; \subset \; 
R_{1,\t{\mu}(w)}.$$
\item for the respective irrelevant ideals we have 
$$
I(R_{1},\mathfrak{F}_{1},\Phi_{1})
\; \subset \;
\sqrt{ \bangle{\mu(I(R_{2},\mathfrak{F}_{2},\Phi_{2}))}}.$$
\end{enumerate}
\end{definition}

Geometrically such a graded homomorphism
is a morphism $\b{\varphi} \colon \b{X}_{1} \to \b{X}_{2}$ 
that maps $\rq{X}_{1}$ into $\rq{X}_{2}$ and is equivariant
in the sense that $\b{\varphi}(t \mal x) = m(t) \mal \b{\varphi}(x)$
holds with $m \colon T_{1} \to T_{2}$ denoting the homomorphism
of tori arising from $\mu \colon K_{2} \to K_{1}$.

\begin{definition}
We say that two graded homomorphisms $\mu$, $\nu$
of the bunched rings 
$(R_{2},\mathfrak{F}_{2},\Phi_{2})$ and 
$(R_{1},\mathfrak{F}_{1},\Phi_{1})$ 
are {\em equivalent\/} if there is a 
homomorphism $\kappa \colon K_{2} \to \KK^{*}$ such that
always $\mu|_{R_{2,w}} = \kappa(w) \nu|_{R_{2,w}}$ 
holds.
\end{definition}

From the geometric point of view, two graded
homomorphisms are equivalent if and only if
the associated equivariant morphisms differ
by multiplication with a torus element.

Note that equivalence of graded homomorphisms 
of bunched rings is compatible with composition, 
that is $\mu \circ \nu$ and $\mu' \circ \nu'$ are 
equivalent provided both $\mu$, $\mu'$ and 
$\nu$, $\nu'$ are. 
Thus, the following definition of morphism
makes the class of bunched rings into a category:

\begin{definition}
\label{moburi}
A {\em morphism of bunched rings\/}
is an equivalence class of graded homomorphisms
of bunched rings.
\end{definition}

Consider a morphism of the bunched rings
$(R_{2},\mathfrak{F}_{2},\Phi_{2})$ and 
$(R_{1},\mathfrak{F}_{1},\Phi_{1})$,
and let it be represented by 
$\mu \colon R_{2} \to R_{1}$.
Then the associated morphism
$\b{\varphi} \colon \b{X}_{1} \to \b{X}_{2}$ 
defines a commutative diagram
\begin{equation}
\label{liftingdiag}
\vcenter{
\xymatrix{
{\rq{X}_{1}} \ar[r]^{\rq{\varphi}} \ar[d]_{{p_{1}}} 
& 
{\rq{X}_{2}} \ar[d]^{{p_{2}}} 
\\
{X_{1}} \ar[r]_{\varphi} 
& {X_{2}}
}
}
\end{equation}
where $X_{i} := X(R_{i},\mathfrak{F}_{i},\Phi_{i})$.
Note that the morphism $\varphi \colon X_{1} \to X_{2}$ does
not depend on the choice of the representative $\mu$.
Moreover, this construction is compatible with composition,
and hence we obtain:

\begin{proposition}
The assignments 
$(R,\mathfrak{F},\Phi) \rightsquigarrow X(R,\mathfrak{F},\Phi)$
and $[\mu] \mapsto \varphi$ define a covariant functor 
from the category of bunched rings to the category of
algebraic varieties.
\end{proposition}

In general this functor is neither injective nor 
surjective on the level of morphisms; counterexamples
occur already in the case of toric varieties,
see \cite[Ex.~9.7 and~9.8]{BeHa2}.
For a more detailed study, the following concept
is crucial, compare~\cite{Be} for the toric case:

\begin{definition}
\label{geopull}
Let $\varphi \colon X_{1} \to X_{2}$ be a dominant morphism. 
We say that a homomorphism 
$\eta \colon \WDiv(X_{2}) \to \WDiv(X_{1})$ 
of the respective groups of Weil divisors 
is a {\em $\varphi$-pullback\/} if
\begin{enumerate}
\item for every principal divisor $D$ on $X_{2}$ we have
  $\eta(D) = \varphi^{*}(D)$;
\item for every effective divisor $D$ on $X_{2}$, 
   the image $\eta(D)$ is again effective;
\item for the supports, we have 
  $\Supp(\eta(D)) \subset \varphi^{-1}(\Supp(D))$.
\end{enumerate}
\end{definition}

In terms of this notion, we can describe when a given
dominant morphism $X_{1} \to X_{2}$ arises from a 
morphism of bunched rings. 
The result generalizes a corresponding statement
on toric varieties, see \cite[Theorem 2.3]{Be}. 

\begin{theorem}
\label{liftchar}
Let $\varphi \colon X_{1} \to X_{2}$ be a dominant morphism. 
Then the morphisms of the bunched rings 
$(R_{2},\mathfrak{F}_{2},X_{2})$ and 
$(R_{1},\mathfrak{F}_{1},X_{1})$ 
that map to $\varphi$ are in  one-to-one
correspondence with the $\varphi$-pullbacks 
$\WDiv(X_{2}) \to \WDiv(X_{1})$.
\end{theorem}

Before giving the proof of this statement, we list a couple
of immediate consequences.

\begin{corollary}
\label{isolift}
Let $\varphi \colon X_{1} \to X_{2}$ be a dominant morphism
such that $\varphi^{-1}(\Sing(X_{2}))$ is of 
codimension at least two in $X_{1}$.
Then $\varphi$ arises from a unique morphism
of the bunched rings $(R_{2},\mathfrak{F}_{2},X_{2})$ 
and $(R_{1},\mathfrak{F}_{1},X_{1})$. 
\endproof
\end{corollary}

This applies in particular to isomorphisms. Thus, we can
answer the question when two bunched rings define isomorphic
varieties:

\begin{corollary}
\label{iso42iso}
Every isomorphism $X_{1} \to X_{2}$ arises from a unique
isomorphism of the bunched 
rings $(R_{2},\mathfrak{F}_{2},X_{2})$ and 
$(R_{1},\mathfrak{F}_{1},X_{1})$.
In particular, $X_{1}$ and $X_{2}$ are isomorphic
if and only if their defining bunched rings are so.
\endproof
\end{corollary}

For a smooth variety $X_{2}$, the following result 
generalizes~\cite[Theorem 5.7]{BeHa1}:

\begin{corollary}
If $X_{2}$ is smooth, then every dominant morphism 
$X_{1} \to X_{2}$ arises from a unique morphism of 
the defining bunched rings.
\endproof
\end{corollary}

\proof[Proof of Theorem~\ref{liftchar}]
We assign to each $\varphi$-pullback
$\eta \colon \WDiv(X_2) \to \WDiv(X_1)$ 
a morphism of the bunched rings
$(R_2,\mathfrak{F}_2,\Phi_2)$ and  
$(R_1,\mathfrak{F}_1,\Phi_1)$ that 
induces $\varphi$. This morphism
will satisfy the property
\begin{equation}
\label{pullbed}
p_1^* \circ \eta 
\; = \; 
\rq{\varphi}^* \circ p_2^*, 
\end{equation}
where 
$p_i^* \colon \WDiv(X_i) \to \WDiv(\rq{X}_i)$
denotes the pullback discussed in 
Section~\ref{toricambient}, and 
$\rq{\varphi} \colon \rq{X}_1 \to \rq{X}_2$
denotes the map arising from any representative
$\mu \colon R_2 \to R_1$ of our morphism
of bunched rings.

First, suppose that $X_2$ is a toric 
variety.
In this step, $\varphi \colon X_{1} \to X_{2}$
need not be dominant;
we only require that $\varphi(X_{1})$
intersects the big torus 
$T_{X_2} \subset X_{2}$.
Moreover, we replace $\WDiv(X_{2})$
with the group of invariant Weil divisors
of $X_{2}$.
Then it makes again sense to speak of 
a $\varphi$-pullback.

Recall that $p_2 \colon \rq{X}_2 \to X_2$ 
is a Cox construction. 
Consider the invariant prime divisors
$D_{X_2}^i$ and the corresponding coordinate
functions $f_i^2 \in \mathfrak{F}_2$ 
on $\b{X}_2$.
Since $R_1$ is factorial, we may choose
for each $i$ an element 
$h_i \in R_1$ with 
$$
\div(h_i) 
\; = \; 
\b{p_1^*(\eta(D_{X_2}^i))}.
$$ 
As its divisor is invariant, this function
is homogeneous, and, by Lemma~\ref{isodescr},
the degree of $h_i$ corresponds to the class
of $\eta(D_{X_2}^i)$ under the identification
$K_1 \cong \Cl(X_1)$.
Using condition~\ref{geopull}~(i),
we obtain a family of graded 
homomorphisms: for  $\alpha \in (\KK^*)^{r_2}$,
we set 
$$
\mu_{\alpha} \colon R_2 \to R_1,
\qquad
f_i^2 \mapsto \alpha_i h_i.
$$

We now adjust $\alpha$ in such a way that
the map 
$\mathcal{O}(T_{\rq{X}_2}) \to \KK(\rq{X}_1)$ 
induced by
$\mu_{\alpha}$ coincides in 
degree zero with the pull back 
$\varphi^* \colon \mathcal{O}(T_{X_2}) \to \KK(X_1)$.
This amounts to solving the equations
$$
\prod \alpha_i^{u_i} 
= 
p_1^* \circ \varphi^* \Bigl( \prod (f_i^2)^{u_i} \Bigr)  \prod h_i^{-u_i},
$$
where the vectors 
$u =(u_{1}, \ldots, u_{r_{2}}) \in \ZZ^{r_2}$ run
through a basis of the space of exponents 
of the monomials of degree $0 \in K_2$.
Note that the expressions on the right hand side 
are in fact nonvanishing constants, 
because we have $R_{1}^* = \KK^*$.

Since the vectors $u$ form a matrix of 
maximal rank, 
the above system of equations is solvable.
Thus, we may fix a solution $\alpha_0$.
We denote the resulting graded homomorphism 
by $\mu \colon R_2 \to R_1$, and the 
corresponding map of the spectra by 
$\b{\varphi} \colon \b{X}_1 \to \b{X}_2$.
Note that, by construction, any invariant
divisor satisfies
\begin{equation}
\label{ratpullbed}
\b{p_1^* \circ \eta(D)} 
\; = \; 
\b{\varphi}^* \bigl(\b{p_2^*(D)}\bigr).
\end{equation}

If $X_2$ is general, then we first construct
$\b{\varphi} \colon \b{X}_1 \to \b{Z}_2$ for
the map $X_1 \to Z_2$ with an ambient toric
variety $Z_2$ of $X_2$. 
By our adjustment of $\alpha$, we then see
that $\b{\varphi}$ gives an equivariant 
map $\b{X}_1 \to \b{X}_2$. 
Note that this map satisfies 
condition~\ref{ratpullbed}, now for any divisor.
Let $\mu \colon R_2 \to R_1$ denote the 
corresponding homomorphism.

We show that $\mu \colon R_{2} \to R_{1}$ 
is a graded homomorphism of bunched rings.
For this, we only have to show that
the corresponding morphism 
$\b{\varphi} \colon \b{X}_{1} \to \b{X}_{2}$
maps $\rq{X}_{1}$ into $\rq{X}_{2}$.

Assume the converse. Then there is an $x \in \rq{X}_{1}$ 
with $\b{\varphi}(x) \in \b{X}_{2} \setminus \rq{X}_{2}$. 
Translated into terms of
$\mathfrak{F}_{2} = \{f_{1}^{2}, \dots, f_{s}^{2}\}$ and the 
corresponding divisors $D^{i}_{X_{2}}$ 
of Proposition \ref{intrinsicPhi}, 
this means that there
is a nonempty 
set of indices $I \subset \{1, \dots, s\}$ 
fulfilling
$$
f^{2}_{i} (\b{\varphi}(x)) = 0
\; \Leftrightarrow \; 
i \in I,
\qquad 
\bigcap_{i \in I} D^{i}_{X_{2}} = \emptyset.$$
The first equation shows 
$x \in \b{\varphi}^{*} (\b{p_{2}^{*} (D^{i}_{X_{2}})})$
for every $i \in I$. 
Condition~\ref{ratpullbed} and $x \in \rq{X}_{1}$
give 
$x \in p_{1}^{*} (\eta(D^{i}_{X_{2}}))$. 
But then Property~\ref{geopull}~(iii)
leads to a contradiction:
$$ 
\varphi(p_{1}(x)) 
\; \in \; 
\bigcap_{i \in I} D_{X_{2}}^{i} 
\; = \; 
\emptyset.
$$ 

The remaining task is to construct the inverse 
mapping from the set of morphisms of bunched rings 
inducing $\varphi$ to the set of 
$\varphi$-pullbacks $\WDiv(X_{2}) \to \WDiv(X_{1})$.
So, let $\mu$ be a graded ring homomorphism 
inducing $\varphi$, and let 
$\rq{\varphi} \colon \rq{X}_{1} \to \rq{X}_{2}$ 
be as in~\ref{liftingdiag}.

Consider $D \in \WDiv(X_{2})$.
Then $E := \rq{\varphi}^{*} (p_{2}^{*}(D))$ 
is well defined, because $\varphi$
is dominant.
Moreover, this divisor is principal,
say $E = \div(h)$ with 
$h \in \KK(\rq{X}_{1})$.
Since $E$ is $T_{1}$-invariant,
the function $h$ is $T_{1}$-homogeneous,
say of degree $w \in K_{1}$.

Since the degrees of the
$f_{i}^{1} \in \mathfrak{F}_{1}$ 
generate $K_{1}$,
the function $ \rq{g} := h/\prod (f_{i}^{1})^{u_{i}}$ 
is of degree zero for suitable 
$u_{i} \in \ZZ$. 
But this implies 
$\rq{g} = p_{1}^{*} (g)$ 
for some $g \in \KK(X_{1})$. Set
$$
\eta(D) 
\; := \; 
\div(g) + \sum_{i=1}^{r} u_{i} D_{X_{1}}^{i}.
$$ 
This gives a well defined map 
$\eta \colon \WDiv(X_{2}) \to \WDiv(X_{1})$,
for which the properties of a $\varphi$-pullback 
are directly verified.
By construction, condition~\ref{pullbed}
is satisfied.

To conclude the proof, we have to check 
that the two assignments are inverse to 
each other. 
Starting with a $\varphi$-pullback,
constructing the associated morphism of
bunched rings, and then returning to
the $\varphi$-pullbacks gives 
nothing new, because condition~\ref{pullbed}
determines the $\varphi$-pullbacks.
For the other way round, we may argue
similarly: the class of a graded
homomorphism is fixed by 
condition~\ref{pullbed}.
\endproof

\section{Picard group, semiample, and ample cone}
\label{ample}

We present explicit descriptions of the Picard group,
the semiample and the ample cone of the variety $X$
arising from a bunched ring $(R,\mathfrak{F},\Phi)$.
The descriptions are given in terms of the isomorphism
$\b{\mathfrak{D}}_{X} \colon K \to \Cl(X)$
provided in Lemma~\ref{isodescr};
they generalize the corresponding results on toric 
varieties of~\cite[Sec.~10]{BeHa2}.

\begin{proposition}
\label{picarddescr}
Let $X$ be the variety arising from a bunched ring
$(R,\mathfrak{F},\Phi)$,
and let $(E \topto{Q} K,\gamma)$ be the projected cone
associated to $\mathfrak{F}$.
Then, in $K \cong \Cl(X)$, the Picard group 
is given by
$$
\Pic(X) 
\; = \; 
\bigcap_{\gamma_{0} \in \cov(\Phi)} 
              Q(\lin(\gamma_{0}) \cap E).
$$
\end{proposition}

\proof 
We have to figure out when a $w \in K$ 
represents a Cartier divisor. 
According Proposition~\ref{divisors}, this 
is the case if and only if
$$
w
\; \in \; 
\bigcap_{\gamma_{0} \in \rel(\Phi)} 
              Q(\lin(\gamma_{0}) \cap E).
$$
Now recall that $\cov(\Phi)$ is a subset of $\rel(\Phi)$,
and that for any element $\gamma_0 \in \rel(\Phi)$,
there is a face $\gamma_1 \preceq \gamma_0$ with
$\gamma_1 \in \cov(\Phi)$.
So, it suffices to take the intersection over
$\cov(\Phi)$.
\endproof

We now describe some cones of 
rational divisor classes on $X$ 
that are of general interest. 
First we consider the cone 
$\Eff(X) \subset \Cl_{\QQ}(X)$ 
generated by the classes of the
effective divisors and the 
moving cone
$\Mov(X) \subset \Cl_{\QQ}(X)$.
For the definition of the latter, 
recall that the base locus
of a divisor $D \in \WDiv(X)$ is
$$
\bigcap_{f \in \Gamma(X,\mathcal{O}(D))}
Z(f),
\qquad
\text{where }
Z(f) := \Supp(\div(f)+D).
$$ 
Then the moving cone $\Mov(X) \subset \Cl_{\QQ}(X)$ 
is the cone generated by the classes of 
those divisors, 
the base loci of which are of codimension at 
least two.
Our result generalizes corresponding statements
in the toric case;
for the description of the toric moving cone 
see e.g.~\cite[Prop.~3.1]{Tsch}.

\begin{proposition}
\label{movingcone}
Let $(R,\mathfrak{F},\Phi)$ be a bunched ring
with associated projected cone 
$(E \topto{Q} K,\gamma)$,
and let $X := X(R,\mathfrak{F},\Phi)$.
\begin{enumerate}
\item Let  $e_1, \ldots, e_r$ be the primitive generators of $\gamma$.
Then, in $K_{\QQ} = \Cl_{\QQ}(X)$, the effective cone of $X$ is given by
$$\Eff(X) = \cone(Q(e_1), \ldots, Q(e_r)).$$
\item Let $\gamma_1, \ldots, \gamma_r$ denote the facets of $\gamma$.
Then the moving cone of $X$ is of full dimension in 
$K_{\QQ} = \Cl_{\QQ}(X)$, and it is given by
$$ \Mov(X) = Q(\gamma_1) \cap \ldots \cap Q(\gamma_r). $$ 
\end{enumerate}
\end{proposition}  

\proof
A divisor $D$ is linearly equivalent to an effective one, 
if and only if $D$ admits a global section. 
Since $R$ is the total coordinate ring of $X$, the latter
holds, if and only if we have $R_w \ne 0$ for the degree 
$w \in K$ corresponding to the class of $D$.
In other words, the effective cone in $K$ is precisely the 
grading cone of $R$.
This proves~(i).

To verify~(ii), let $D$ be any effective divisor on $X$.
We have to show that $D$ has no fixed components if and
only if its corresponding degree $w$ lies in every
$Q(\gamma_i)$.
According to~(i), we may assume that $D$ is a nonnegative
linear combination of the divisors $D_X^i$ introduced in
Proposition~\ref{intrinsicPhi}.  

Now, by definition, $D_X^i$ is a fixed component of $D$ 
if and only if for every section 
$h \in \Gamma(X,\mathcal{O}(D))$, the prime divisor
$D_X^i$ occurs in $\Supp(\div(h) + D))$.
This holds if and only if the pullback divisors with
respect to $p_X \colon \rq{X} \to X$ as defined
in Section~\ref{toricambient} satisfy 
$$
p_X^*(D_X^i) 
\; \le \; 
p_X^*(D) + \div(p_X^*(h)).
$$

Using factoriality of $R$, we can conclude 
that $D_X^i$ is a fixed component of $D$ if and 
only if $f_i$ divides 
every $f \in R_w$, where $w \in K$ denotes the 
class of $D$.
But the latter can be expressed in terms
of degrees: it holds if and only if
the class $w \in K$ can not be written as 
nonnegative linear combination of the 
$w_j = Q(e_j)$ with $j \ne i$.

In other words, $D_X^i$ is not a fixed component 
of $D$ if and only if $D$ corresponds to an element of
$\cone(Q(e_j); \; j \ne i)$. 
This verifies the explicit description of the moving cone.
The fact that $\Mov(X)$ is of full dimension then is
an immediate consequence of Lemma~\ref{1Gbunch}.
\endproof

Next we describe the semiample cone 
$\SAmple(X)$ and the ample 
cone $\Ample(X)$ in $\Cl_{\QQ}(X)$.
Recall that, by definition, 
the semiample cone is 
generated by the classes of the
divisors having an empty base locus.
Moreover, the ample cone is generated
by the classes of divisors 
$D \in \WDiv(X)$ that provide a cover
of $X$ by affine open sets of the form
$X \setminus Z(f)$ with 
$f \in \Gamma(X,\mathcal{O}(D))$.

The descriptions of the semiample and the 
ample cone can be done purely in
terms of the $\mathfrak{F}$-bunch $\Phi$:

\begin{theorem}
\label{amplecone}
Let $X$ be the variety arising from a bunched ring
$(R,\mathfrak{F},\Phi)$ with grading lattice $K$.
Then, in  $K_{\QQ} \cong \Cl_{\QQ}(X)$, the 
semiample and the ample cone are given by
$$
\SAmple(X)
\; = \;  
\bigcap_{\tau \in \Phi} \tau,
\qquad
\Ample(X)
\; = \; 
\bigcap_{\tau \in \Phi} \tau^{\circ}.
$$
\end{theorem}

\proof Consider the embedding 
$X \subset Z := Z(R,\mathfrak{F},\Phi)$ 
into the minimal toric variety presented 
in Proposition~\ref{minimalambient}.
Recall that every closed $T_{Z}$-orbit has 
nontrivial intersection with $X$. 
Moreover, $X \subset Z$ defines a pullback
isomorphism $\Cl(Z) \to \Cl(X)$.
This is due to Proposition~\ref{toricCl2XCl}
and the fact $Z$ is an open toric subvariety
of any ambient toric variety of $X$.

We claim that the semiample cone $\SAmple(Z)$ 
is mapped onto $\SAmple(X)$. Indeed, given a
class $[D] \in \SAmple(X)$, this class extends 
to a class $[D'] \in \Cl(Z)$, represented by some 
$T_{Z}$-invariant divisor $D'$ on $Z$.
Since the base locus of $D'$ is $T_{Z}$-invariant 
and $X$ meets every closed $T_{Z}$-orbit it follows
that $D'$ is semiample.
This proves $\SAmple(Z) = \SAmple(X)$.
Thus, we only have to verify 
$$ 
\SAmple(Z) 
\; = \;  
\bigcap_{\tau \in \Phi} \tau.
$$

Consider the projected cone $(E \topto{Q} K, \gamma)$
associated to $\mathfrak{F} \subset R$, and
its dual projected cone $(F \topto{P} N, \delta)$.
We denote by $\rq{\Delta}$ the fan in $F$ having
$\gamma_0^*$, where $\gamma_0 \in \cov(\Phi)$, as
its maximal cones.
Then, by definition, the fan 
$\Delta := \Delta(R,\mathfrak{F},\Phi)$
of $Z$ has the images $P(\delta_0)$,
where $\delta_0 \in \rq{\Delta}^{\max}$,
as its maximal cones.

In the sequel, we follow the lines of the proof
of~\cite[Theorem~10.2]{BeHa2}.
Let $e_{\varrho}$ denote the primitive
generator of $\delta$, which
is mapped to the primitive
generator of $\varrho \in \Delta^{(1)}$.
Then~\cite[Prop.~10.1 and~10.6]{BeHa2},
provide an isomorphism
$$
\mathfrak{D}_{Z} \colon E \to \WDiv^{T_{Z}}(Z),
\qquad
\rq{w} \mapsto \sum_{\varrho \in \Delta^{(1)}} 
\bangle{\rq{w}, e_{\varrho}} \b{T_{Z} \mal z_{\varrho}}.
$$

Next recall from~\cite{Fu} that an invariant divisor 
$D$ on $Z$ is semiample if and only if 
it is described by a convex conewise linear
function $(u_{\sigma})$
living on the support of the fan $\Delta$.
Here convexity means that
$u_{\sigma} - u_{\sigma'}$ is nonnegative
on $\sigma \setminus \sigma'$ for any two
$\sigma, \sigma' \in \Delta^{\max}$. 
Here we do not restrict ourselves to integral 
linear forms $u_{\sigma}$, but we also admit 
rational ones. 

Now, let $w \in K$ represent a semiample class.
This class comes from a semiample rational
divisor $D = \mathfrak{D}_{Z}(\rq{w})$ with 
$\rq{w} \in E_{\QQ}$ such that $w = Q(\rq{w})$.
Let $(u_{\sigma})$ be a convex function describing~$D$.
We have to interpret $(u_{\sigma})$ in terms 
of the fan $\rq{\Delta}$.
For $\sigma \in \Delta^{\max}$, let 
$\rq{\sigma} \in \rq{\Delta}$ be the 
(unique) cone with $P(\rq{\sigma}) = \sigma$.

Then, for $\ell_{\sigma} := \rq{w} - P^{*}(u_{\sigma})$,
each $\ell_{\sigma'} - \ell_{\sigma}$ 
is nonnegative on $\rq{\sigma} \setminus \rq{\sigma}'$.
Since $\ell_{\sigma} \in \rq{\sigma}^{\perp}$
holds, this is equivalent to the nonnegativity
of $\ell_{\sigma'}$ on every $\rq{\sigma} \setminus \rq{\sigma}'$.
Since all rays of the cone $\delta$ occur in the fan $\rq{\Delta}$,
the latter is valid if and only if
$\ell_{\sigma} \in \rq{\sigma}^{*}$ 
holds for all $\sigma$. This in turn implies
that for every $\sigma \in \Delta^{\max}$ we have
$$
w 
= 
Q(\rq{w}) 
= 
Q(\ell_{\sigma}) 
\in Q(\rq{\sigma}^{*}).
$$

By the definition of the fan $\Delta$, 
the cones of $\Phi$ are the minimal 
ones among the images $Q(\rq{\sigma}^{*})$, where 
$\sigma \in \Delta^{\max}$.
Consequently, $w$ lies in the intersection over the cones
of $\Phi$.

Conversely, if $w$ belongs to the intersection 
of all cones in $\Phi$, then $w$ is as well 
an element of every image $Q(\rq{\sigma}^{*})$
with $\sigma \in \Delta^{\max}$.
In particular, we find for every 
$\sigma \in \Delta^{\max}$ an 
$\ell_\sigma \in \rq{\sigma}^*$ mapping to $w$. 
Reversing the above arguments, 
we see that $u_{\sigma} := \rq{w} - \ell_{\sigma}$
is a convex support function describing 
the representative 
$D := \mathfrak{D}_{Z}(\rq{w})$
of the class
$\b{\mathfrak{D}}_{Z}(w)$.

The description of the ample cone is obtained in an
analogous manner, and therefore is left to the reader. 
For the equality 
$\Ample(Z)= \Ample(X)$ 
we just note that any $T_{Z}$-invariant divisor 
has a $T_{Z}$-invariant ample locus, i.e., 
the set of points admitting an affine neighbourhood
$Z_{f}$ with a global section $f$ of $\mathcal{O}(D)$.
\endproof

For a variety $X$ with finitely generated free 
divisor class group, the Picard group $\Pic(X)$ 
is equal to the N\'{e}ron-Severi group $N^1(X)$.
In the associated vector space 
$N^1_{\QQ}(X) = \Pic_{\QQ}(X)$, one 
studies the numerically effective cone,
i.e., the cone generated by the classes of 
the Cartier divisors having nonnegative
intersection with all (effective) curves.

For projective varieties with free 
finitely generated divisor class group 
and finitely generated total coordinate ring,
we obtain as an immediate consequence
of the preceding results and Kleiman's
ampleness criterion that the numerically 
effective cone and the semiample cone coincide; 
compare~\cite[Prop.~1.11]{KeHu} for the 
$\QQ$-factorial case. 

\begin{corollary}
\label{projamplecone}
Let $X$ be a normal projective variety 
with free finitely generated divisor class group
and finitely generated total coordinate ring.
\begin{enumerate}
\item The ample and the semiample cone of $X$
are of full dimension in $\Pic_{\QQ}(X)$.
\item The ample cone of $X$ is the relative interior
of the semiample cone of $X$.
\item The numerically effective cone 
and the semiample cone of $X$ coincide.
\end{enumerate} 
\end{corollary}

\proof
We may assume that $X$ arises from a bunched ring
$(R,\mathfrak{F},\Phi)$.
Since $X$ is projective, it admits ample divisors.
Thus, Theorem~\ref{amplecone} tells us that the
intersection of all relative interiors 
$\tau^{\circ}$, where $\tau \in \Phi$, is 
nonempty.
Consequently, 
Proposition~\ref{picarddescr} and,
once more, Theorem~\ref{amplecone} 
give the first two assertions;
the first one is also a special case 
of~\cite[Thm.~IV.2.1]{Kleim}.
The third assertion follows from the
the second one and  
fact that for projective varieties,
the numerically effective 
cone is the closure of the ample cone,
see~\cite[Thm.~IV.2.1]{Kleim}. 
\endproof

The Mori cone is the cone 
$\Mori(X)$ in the dual space of
$N_{\QQ}^1(X)$ generated by the 
numerical equivalence classes of 
effective curves.
In our setup, the Mori cone lives
in the dual vector space of 
$\Pic_{\QQ}(X)$,
and we obtain an explicit description 
for it.

\begin{corollary}
\label{moricone}
Let $X$ be a projective variety 
arising from a bunched ring 
$(R,\mathfrak{F},\Phi)$.
Then, in $\Pic_{\QQ}(X)^*$,
the Mori cone of $X$ is the 
strictly convex cone given by
$$
\Mori(X)
\; = \;
\sum_{\tau \in \Phi} (\Pic(X)_\QQ \cap \tau)^\vee.
$$ 
\end{corollary} 

\proof
The Mori cone is dual to
the numerically effective cone
in $\Pic_{\QQ}(X)$, by definition
of the latter. 
Hence Corollary~\ref{projamplecone},
Proposition~\ref{picarddescr},
and Theorem~\ref{amplecone} 
give the assertion.
\endproof

Finally, we note that 
the proofs of the results of this section 
show that the various cones of divisor 
classes are inherited from suitable toric 
ambient varieties.

\begin{remark}
Let $X := X(R,\mathfrak{F},\Phi)$ be the variety arising 
from a bunched ring $(R,\mathfrak{F},\Phi)$.
\begin{enumerate}
\item For any toric ambient variety $Z$ of $X$ 
in the sense of Proposition~\ref{xintv} we have
$$ 
\Eff(X) = \Eff(Z), 
\qquad
\Mov(X) = \Mov(Z). 
$$
\item For the minimal toric ambient variety 
$Z := Z(R,\mathfrak{F},\Phi)$ of $X$ we have
$$
\SAmple(X) = \SAmple(Z),
\quad
\Ample(X) =  \Ample(Z).
$$
\end{enumerate}
\end{remark}

It should be mentioned that 
the second statement may become false if 
one enlarges the minimal toric 
ambient variety, e.g., if one makes it $A_2$-maximal;
see~\cite{Sze} and Example~\ref{hatschi}.  

\section{Projectivity and Fano criteria}
\label{kleichev} 

In this section, we apply the description of the ample
cone given in the previous one to obtain projectivity and 
Fano criteria.
The first result is an effective version of the 
Kleiman-Chevalley Criterion.
In our setup, it sharpens considerations
performed in~\cite[p.~307]{Ko}.

\begin{theorem}
\label{kleimanchevalley}
Let $X$ be a (normal) $A_2$-maximal 
$\QQ$-factorial variety with 
$\mathcal{O}(X) = \KK$,
free divisor class group of rank $k$
and finitely generated total coordinate 
ring.
If any collection $x_1, \ldots, x_k \in X$ 
admits a common affine neighbourhood in $X$,
then $X$ is projective. 
\end{theorem}

For the proof, and also for later purposes,
we need the following characterization of 
existence of nontrivial global functions.

\begin{proposition}
\label{globfunconst}
Let $X$ arise from a bunched ring $(R,\mathfrak{F}, \Phi)$,
and let $(E \topto{Q} K, \gamma)$ be the 
projected cone associated to $\mathfrak{F}$. 
Then the following statements are equivalent:
\begin{enumerate}
\item $\mathcal{O}(X) = \KK$,
\item the cone $Q(\gamma)$ is strictly convex and $Q(e_{i}) \neq 0$ 
holds for each primitive generator $e_{i}$ of $\gamma$.
\end{enumerate}
\end{proposition}

\proof
Choose a bunch $\Theta$ extending $\Phi$, and 
consider the corresponding ambient toric variety 
$Z$ of $X$.
We claim that $\mathcal{O}(X) = \KK$ is equivalent
to $\mathcal{O}(Z) = \KK$.
In terms of $\b{X} = \Spec(R)$ and 
$\b{Z} = \Spec(\KK[E \cap \gamma])$,
and the torus $T = \Spec(\KK[K])$,
this means to show that
$$
\mathcal{O}(\b{X})^{T} = \KK
\; \Leftrightarrow \;
\mathcal{O}(\b{Z})^{T} = \KK.
$$

The implication ``$\Leftarrow$'' is easy.
The reverse implication follows from the fact 
that $\b{X}$ meets the big torus and every 
coordinate hyperplane of $\b{Z} = \KK^{r}$:
if $\b{Z}$ admits a nonconstant
$T$-invariant function,
then it admits also a nonconstant $T$-invariant
monomial.
This in turn restricts to a nonconstant function on $\b{X}$.

Thus we verified that $\mathcal{O}(X) = \KK$ is equivalent
to $\mathcal{O}(Z) = \KK$.
But the latter property was characterized in 
\cite[Proposition~8.4]{BeHa2} in terms of the bunch $\Theta$:
It is valid if and only if $Q$ contracts no ray of $\gamma$,
and $Q(\gamma)$ is strictly convex.
\endproof

Another preparatory result concerns the question,
whether or not an $A_2$-maximal variety $X$ with 
finitely generated free divisor class group and 
finitely generated total coordinate ring is
necessarily complete if it satisfies 
$\mathcal{O}(X) = \KK$.
First note that, in general, the answer to 
this question is no.

In fact, there exist examples of non-complete 
fans as in Corollary~\ref{toric2a2max},
the support of which generates
the whole ambient vector 
space as a cone, see~\cite{Ew}.
The corresponding toric 
varieties are $A_2$-maximal and non-complete.
However, such varieties can never be 
quasiprojective, as we will see now.

\begin{proposition}
\label{quasiproj2proj}
Let $X$ be an 
$A_2$-maximal variety 
with free finitely generated divisor class group 
and finitely generated total coordinate ring.
If $X$ is quasiprojective and 
$\mathcal{O}(X) = \KK$ holds,
then $X$ is even projective. 
\end{proposition}

\proof
We may assume that $X$ arises from a bunched
ring $(R,\mathfrak{F},\Phi)$.
Then $X$ is in particular a good quotient of 
an open subset $\rq{X}$ of $\b{X} := \Spec(R)$ 
by the action of a torus $T$ with $\dim(T)$ 
equal to $\dim(\b{X})-\dim(X)$.

Since $X = \rq{X} \quot T$ is quasiprojective,
it admits an open embedding into the quotient 
$X' = \rq{X}'(w)  \quot T$ of the set 
$\rq{X}' \subset \rq{X}$ of
semistable points of some
$T$-linearization of the trivial bundle over 
$\b{X}$, compare~\cite[Converse~1.13]{Mu}. 
The assumption 
$\mathcal{O}(\b{X})^T = \mathcal{O}(X) = \KK$
implies that the quotient space $X'$ is 
projective. 
 
Since $\rq{X}$ has a complement at least two
in $\b{X}$, 
the same holds for $\rq{X}$ in $\rq{X}'$,
and hence for $X$ in $X'$. 
By $A_2$-maximality, we have $X=X'$,
and hence $X$ is projective.
\endproof

\proof[Proof of Theorem~\ref{kleimanchevalley}]
If $X$ is as in the assertion, then 
$X = X(R,\mathfrak{F},\Phi)$ holds
with a bunched ring $(R,\mathfrak{F},\Phi)$.
According to Proposition~\ref{quasiproj2proj}
and Theorem~\ref{amplecone},
we have to show that the intersection over
all relative interiors $\tau^{\circ}$,
where $\tau \in \Phi$, is nonempty.
As usual, let $(E \topto{Q} K, \gamma)$ denote
the projected cone associated to $\mathfrak{F}$.

We perform a first reduction step.
Consider the convex hull
$\vartheta$ of all $\tau \in \Phi$.
Then Proposition~\ref{globfunconst}
tells us that $\vartheta$ is a strictly convex cone, 
because we assumed $\mathcal{O}(X) = \KK$. 
Thus, choosing a $w \in \vartheta^{\circ}$
and a suitable affine hyperplane $H \subset K_\QQ$ 
orthogonal to $w$, 
we obtain for every $\tau \in \Phi$ 
a nonempty polytope $B(\tau) := H \cap \tau$.
By Corollary~\ref{singularities},
we have $\dim(\tau) = k$,
so $B(\tau)$ is of dimension $k-1$.

Clearly it suffices to show that the intersection
of all $B(\tau)^{\circ}$, 
where $\tau \in \Phi$, is nonempty.
According to Helly's Theorem, see~\cite[Satz~I.7.1]{Le},
the latter holds if any intersection 
$B(\tau_1)^{\circ} \cap \ldots \cap B(\tau_k)^{\circ}$
is nonempty.
Thus is suffices to that for any collection
$\tau_1, \ldots, \tau_k$,
the intersection of their relative interiors is
nonempty.

Choose a bunch $\Theta$ extending $\Phi$,
and consider the corresponding ambient toric 
variety $Z$ of $X$.
Given $\tau_1, \ldots, \tau_k$ as before, 
choose $\gamma_i \in \rel(\Phi)$ with 
$Q(\gamma_i) = \tau_i$, and let $Z' \subset Z$ 
be the associated open toric subvariety, i.e., 
$$ Z' = \bigcup_{i=1}^k \Spec(\KK[\ker(Q) \cap \gamma_i]). $$

We claim that $Z'$ is quasiprojective.
First note that $Z'$ is $\QQ$-factorial.
Moreover, the big torus $T' \subset Z'$
has at most $k$ closed orbits in $Z'$, 
and any such orbit contains a point of $X$.
Thus we find an affine open subset $W \subset Z'$
that intersects every closed $T'$-orbit of $Z'$.

Since $Z'$ is also $\QQ$-factorial,
the complement $Z' \setminus W$ is the support 
of a Cartier divisor $D'$ on $Z'$.
Since the translates $t' \mal W$ cover
$Z'$, and $D'$ is linearly equivalent to
an invariant divisor, we can conclude that 
$D'$ is ample. Thus $Z'$ is quasiprojective. 
As the toric part of the proof of 
Theorem~\ref{amplecone} shows,
this implies that
$\tau_1^{\circ} \cap \ldots \cap \tau_k^{\circ}$
is nonempty.
\endproof

\begin{corollary}
Every (normal) $\QQ$-factorial 
$A_2$-variety 
$X$ with $\mathcal{O}(X) = \KK$,
free divisor class group of rank two
and finitely generated total coordinate 
ring is quasiprojective.
\end{corollary}

\proof
Assertion~\ref{firstresult}~a) and 
Theorem~\ref{firstresult} tell us that 
$X$ is an open subset of some $A_2$-maximal
variety with the same properties.
Thus, Theorem~\ref{kleimanchevalley} gives 
the desired result.
\endproof

In this corollary, the assumption that
$X$ is an $A_2$-variety is essential:
one obtains complete nonprojective threefolds
with divisor class group isomorphic to 
$\ZZ^2$ as geometric quotient spaces of 
a certain $\KK^*$-action on the Grassmannian
$G(2,4)$, see for example~\cite{BBSw1} 
and~\cite{Sw}.

Our next aim is to give some Fano criteria for
the variety $X$ arising from a bunched ring 
$(R,\mathfrak{F},\Phi)$.
For this, we first have to determine the 
canonical divisor class of $X$.
This can be easily done if $X$ is a complete 
intersection in some of its ambient toric 
varieties.

Let $(E \topto{Q} K, \gamma)$ be the projected cone 
associated to $\mathfrak{F}$. Denote the primitive
generators of $\gamma$ by $e_{1}, \dots, e_{r}$, 
and set $w_{i} := Q(e_{i})$.

\begin{proposition}\label{candivdescr}
Suppose that the relations of $\mathfrak{F}$ 
are generated by $K$-homogeneous polynomials 
$g_{1}, \ldots, g_{d} \in \KK[T_{1}, \ldots, T_{r}]$,
where $d := r - \rank(K) - \dim(X)$.
Then, in $K \cong \Cl(X)$,
the canonical divisor class of $X$ is given by
$$
D^{c}_{X}
\; = \;
\sum_{i=1}^{d} \deg(g_{i})
- 
\sum_{i=1}^{r} w_{i}
$$
\end{proposition}

\proof
Consider the embedding $X \to Z$ into any
ambient toric variety $Z$ of $X$.
Then Proposition~\ref{singulemb} tells us
that the respective embeddings of the smooth
loci fit into a commutative diagram
$$
\xymatrix{
{X_{\reg}} \ar[r] \ar[d] & {Z_{\reg}} \ar[d] \\
X \ar[r] & Z
}
$$

Note that all maps induce isomorphisms on 
the respective divisor class groups. 
Moreover, the restrictions of the canonical divisors $K_{X}$ and
$K_{Z}$ of $X$ and $Z$ give the canonical divisors of $X_{\reg}$ and
$Z_{\reg}$, respectively, see for example \cite[p.~164]{Ma}. 
Thus, we may assume that $X$ and its ambient toric 
variety $Z$ are smooth.

Let $\mathcal{I} \subset \mathcal{O}_Z$ be the 
ideal sheaf of $X$.
Then the normal sheaf of $X$ in $Z$ is
the rank $d$ locally free sheaf
$\mathcal{N}_X := (\mathcal{I} / \mathcal{I}^2)^*$,
and a canonical bundle on $X$ can be 
obtained as follows, see for 
example~\cite[Prop.~II.8.20]{Hart}: 
$$
\mathcal{K}_X
= 
\mathcal{K}_Z \vert_X \otimes 
\Bigl({\bigwedge}^{d} \mathcal{N}_X \Bigl).
$$

Choose a cover of $Z$ by open subsets $U_{i}$ such that  
$\mathcal{I} / \mathcal{I}^2$ is free over 
$U_{i}$. 
Then the $g_{l}$ generate the relations 
of the $f_{j}$ over 
$\rq{U}_{i} := p_{Z}^{-1}(U_{i})$.
Thus, after suitably refining the cover,
we find functions 
$h_{il} \in \mathcal{O}^{*}(\rq{U}_{i})$ of degree
$\deg(g_{l})$ 
such that $\mathcal{I} / \mathcal{I}^2 (U_{i})$ 
is generated by $g_{1}/h_{i1}, \ldots, g_{d}/h_{id}$.
Therefore over $U_{i}$, the $d$-th exterior power 
of $\mathcal{I} / \mathcal{I}^2$ is generated by 
the function
$$ \frac{g_{1}}{h_{i1}} \wedge \dots \wedge \frac{g_{d}}{h_{id}}.$$ 

Now, Lemma~\ref{isodescr} tells us that the class of
$\bigwedge^{d} \mathcal{I} / \mathcal{I}^{2}$ in $K$ 
is minus the sum of the degrees of the $g_{j}$. 
As $\bigwedge^{d} \mathcal{N}_{X}$ is the dual sheaf, 
its class is $\deg(g_{1}) + \ldots + \deg(g_{l})$. 
Furthermore, from~\cite[Section~10]{BeHa2} we know that 
the class of the canonical divisor of $Z$ in $K$ is
given by $-(w_{1} + \ldots + w_{r})$. 
Putting all together we arrive at the
assertion.
\endproof

A variety is said to be ($\QQ$-)Gorenstein if (some multiple of) 
its anticanonical divisor is Cartier. In our setting, we obtain

\begin{corollary}
In the setting of Proposition~\ref{candivdescr}, the variety $X$ is
\begin{enumerate}
\item $\QQ$-Gorenstein if and only if
$$ 
\sum_{i=1}^{r }w_{i} - \sum_{j=1}^{d} \deg(g_{j})
\; \in \; 
\bigcap_{\tau \in \Phi} \lin(\tau),
$$
\item Gorenstein if and only if
$$ 
\sum_{i=1}^{r }w_{i} - \sum_{j=1}^{d} \deg(g_{j})
\; \in \; 
\bigcap_{\gamma_{0} \in \rel(\Phi)} 
              Q(\lin(\gamma_{0}) \cap E).
$$
\end{enumerate}
\end{corollary}

Finally, as announced, we determine, 
when $X$ is a ($\QQ$-)Fano variety,
what means that (some multiple of) 
its anticanonical divisor is ample.

\begin{corollary}\label{fanocrit}
In the situation of Proposition~\ref{candivdescr},
the variety $X$ is
\begin{enumerate}
\item  $\QQ$-Fano if and only if we have
$$
\sum_{i=1}^{r }w_{i} - \sum_{j=1}^{d} \deg(g_{j})
\; \in \; 
\bigcap_{\tau \in \Phi} \tau^{\circ},
$$
\item  Fano if and only if we have
$$
\sum_{i=1}^{r }w_{i} - \sum_{j=1}^{d} \deg(g_{j}),
\; \in \; 
\bigcap_{\tau \in \Phi} \tau^{\circ}
\cap
\bigcap_{\gamma_{0} \in \rel(\Phi)} 
              Q(\lin(\gamma_{0}) \cap E).
$$
\end{enumerate}
\end{corollary}

\section{Intrinsic quadrics I}
\label{intrinsic}

In the following two sections, we
consider varieties that have a total 
coordinate ring defined by a single 
nontrivial relation.
We focus on the case of a quadratic 
relation, 
and study these varieties in terms 
of their bunched rings.
For example, we shall classify the
smooth complete varieties of this 
type having a divisor class group 
of rank two.

\begin{definition}
\label{quadricaldef}
By an {\em intrinsic hypersurface},
we mean a normal $A_{2}$-maximal 
variety $X$ with 
$\mathcal{O}^*(X) = \KK^*$ and
free finitely generated divisor 
class group $K := \Cl(X)$
such that the total coordinate ring 
of $X$ admits a representation
$$
\mathcal{R}(X)
\; \cong \;
\KK[T_{1}, \ldots, T_{r}] / \bangle{g},
\qquad
\text{with}
\quad
\KK[T_{1}, \ldots, T_{r}]
\; = \;
\bigoplus_{w \in K}
\KK[T_{1}, \ldots, T_{r}]_w, 
$$
where the relation $g$ is homogeneous and
nontrivial, and the classes of 
$T_{1}, \ldots, T_{r}$ are homogeneous
pairwise nonassociated 
prime elements, forming a system of
generators of minimal length for
$\mathcal{R}(X)$.
\end{definition} 

By Theorem~\ref{firstresult}, 
every intrinsic hypersurface is of the form
$X = X(R,\mathfrak{F},\Phi)$ with a bunched 
ring $(R,\mathfrak{F},\Phi)$.
Moreover, if $\mathfrak{F} \subset R$ is a system of 
generators of minimal length, then 
Proposition~\ref{xintv} says that $X$
is a hypersurface in each of its ambient 
toric varieties; this explains the name.

We begin with some general observations. 
Let $R$ be any factorial ring
satisfying $R^* = \KK^*$, and suppose that there is a 
representation
$$
R
\; = \;
\KK[T_{1}, \ldots, T_{r}] / \bangle{g}
$$
for some nontrivial $g \in \KK[T_{1}, \ldots, T_{r}]$.
The following criteria are 
often helpful for deciding whether or not the 
class of $T_{i}$ is prime in $R$: 

\begin{remark}\label{primetest}
Let $g \in \KK[T_{1}, \ldots, T_{r}]$ be such that 
$R := \KK[T_{1}, \ldots, T_{r}] / \bangle{g}$
is factorial.
\begin{enumerate}
\item The class of $T_{i}$ is prime in $R$,
if and only if  
$g(T_{1}, \dots,T_{i-1},0,T_{i+1}, \dots, T_{r})$ 
is irreducible in 
$\KK[T_{1}, \dots,T_{i-1},T_{i+1}, \dots, T_{r}]$.
\item If $R$ admits an $\NN$-grading such that 
$R_0^* = \KK^*$ and $\deg(T_i) = \pm 1$ hold,
then the class of $T_i$ is prime in $R$.
\end{enumerate}
\end{remark}

In the sequel, we will assume that $R$ is 
as in Definition~\ref{quadricaldef}; 
that means that $\KK[T_{1}, \ldots, T_{r}]$
is graded by some lattice $K$,
the $T_i$ and $g$ are homogeneous, 
and the classes $f_i \in R$ of $T_i$ 
are pairwise nonassociated 
prime elements, forming a 
system of generators of minimal length
for $R$.
Moreover, we require that the degrees 
$$
w_i 
\; := \; 
\deg(f_i) 
\; = \; 
\deg(T_i)
$$
generate $K$ as a lattice. 
Let $(E \topto{Q} K,\gamma)$ denote the
projected cone associated to $\mathfrak{F}$.
Thus, $E = \ZZ^r$ and $Q(e_i) = w_i$
hold, and $\gamma \subset E_{\QQ}$ is the
positive orthant. 
We want to determine the $\mathfrak{F}$-faces for 
$\mathfrak{F} := \{f_{1}, \ldots, f_{r}\}$.
For a given face $\gamma_{0} \preceq \gamma$,
consider the number
$$
\nu(\gamma_{0}) := \# (u \in \gamma_{0} \cap E; \; \alpha_{u} \ne 0),
\qquad 
\text{where }
g = \sum \alpha_{u} T^{u}.
$$

\begin{proposition}\label{hypersfrel}
A face $\gamma_{0} \preceq \gamma$ is an $\mathfrak{F}$-face
if and only if $\nu(\gamma_{0})$ differs from one.
\end{proposition}    

\proof
We treat as an example a face of the form 
$\gamma_{0} := \cone(e_{1}, \ldots, e_{s})$. 
Consider the polynomial
$$
h(T_{1}, \ldots, T_{s})
:=
g(T_{1}, \ldots, T_{s},0,\ldots,0).
$$
Then, by Remark~\ref{geomfface}, 
the cone $\gamma_{0}$ is an $\mathfrak{F}$-face 
if and only if there is a point $z \in \KK^{r}$
with
$$
z_{1}, \ldots, z_{s} \in \KK^{*},
\qquad
z_{s+1}, \ldots, z_{r} = 0,
\qquad
g(z) = h(z_{1}, \ldots, z_{s}) = 0.
$$
Observe that the polynomial $h$ is a sum of exactly 
$\nu(\gamma_{0})$ (nontrivially scaled) monomials. 
Thus, the desired 
$z_{1}, \ldots, z_{s}$ exist if and only if
$\nu(\gamma_{0})$ differs from one.
\endproof

We now specialize to the case that $R$ has a 
defining 
relation $q := g \in \KK[T_1, \ldots, T_r]$ 
of degree two.
Let $\Phi$ be an $\mathfrak{F}$-bunch in 
the projected cone $(E \topto{Q} K,\gamma)$,
and let $X := X(R,\mathfrak{F},\Phi)$ be the
associated variety.
We refer to the variety $X$ as an
{\em intrinsic quadric}.

\begin{proposition}
\label{quadricsetup}
Let $(R,\mathfrak{F},\Phi)$ and 
$X := X(R,\mathfrak{F},\Phi)$ be
as above. 
Then the following statements
hold.
\begin{enumerate}
\item 
If $\mathcal{O}(X)=\KK$ holds, then 
the defining quadric $q$ is 
homogeneous in the usual sense:
\begin{equation*}
q(T_{1}, \ldots, T_{r}) 
\; = \; 
\sum_{i \le j} \alpha_{ij} T_{i}T_{j}.
\end{equation*}
\item 
If $q$ is homogeneous with a representation
as in~(i), then its degree $\deg(q) \in K$
satisfies 
\begin{equation*}
\deg(Q) = w_{i} + w_{j}  \qquad \text{whenever } \alpha_{ij} \ne 0.
\end{equation*}
\end{enumerate}
\end{proposition}

\proof
Only for~(i) there is something to show.
We first show that $q$ contains no linear
term. 
Assume the contrary.
Then, after suitably renumbering the 
indeterminates, we find $\alpha_1 \in \KK^*$ 
and a quadric $q'$ not containing $T_1$
as a monomial such that we have  
$$
q(T_1, \ldots, T_r) 
\; = \;
\alpha_1T_1 + q'(T_1, \ldots, T_r).   
$$

We claim that $q'$ does not even depend 
on $T_1$. 
Otherwise, there were a monomial $P$ in 
$q'(T_1, \ldots, T_r)$ such that $P = T_1P'$ 
holds with a nontrivial monomial $P'$. 
This implies $\deg(P') = 0 \in K$, which, 
in view of Proposition~\ref{globfunconst},
is a contradiction to the assumption 
$\mathcal{O}(X)=\KK$. 

Having convinced ourselves that the quadric $q'$ does 
not depend on $T_1$, we see that there is a 
well defined isomorphism
$$ 
R \to 
\KK[T_2, \ldots, T_r],
\qquad
T_1 \mapsto -q'(T_2, \ldots, T_r),
\
T_2 \mapsto T_2,
\ldots,
T_r \mapsto T_r.
$$
Consequently, $R$ is again a polynomial 
ring. But this contradicts the 
assumption
that the shortest systems of generators 
of $R$ are of length $r = \dim(R) + 1$.

Thus, we saw that $q$ contains no nontrivial
linear terms.
If $q$ would contain a constant term 
$c \ne 0$, then $q$ would be of degree
zero in $K$. This is excluded by
$\mathcal{O}(X) = \KK$ and  
Proposition~\ref{globfunconst}.
\endproof

We now restrict to the case that 
$\mathcal{O}(X) = \KK$ 
holds, and that $q$ is of full rank;
we then refer to $X$ as a 
{\em full intrinsic quadric}. 
Note that factoriality of $R$ implies
$r \ge 5$ for any full intrinsic 
quadric $X$.

Examples are provided as follows: 
take any quadric 
$q \in \KK[T_1, \ldots, T_r]$,
homogeneous in the usual sense and 
of full rank $r \ge 5$,
any $K$-grading making $q$ and 
the $T_i$ homogeneous such that the
degrees of $T_i$ are all nontrivial, 
generate a strictly convex cone, 
and generate $K$ as a lattice,
and, finally, choose an 
$\mathfrak{F}$-bunch for the system
$\mathcal{F} = \{f_1, \ldots, f_r\}$
of the classes $f_i$ of the indeterminates
$T_i$.

\begin{proposition}
\label{quadricprop}
Let $(R,\mathfrak{F},\Phi)$ be as above, 
suppose that the associated 
intrinsic quadric
$X = X(R,\mathfrak{F},\Phi)$ 
is full, and let
$q = \sum \alpha_{ij}T_iT_j$
be the defining equation.
\begin{enumerate}
\item For $\gamma_{0} \in \rel(\Phi)$, 
the stratum $X(\gamma_{0})$ consists of regular points
if and only if $Q(\lin(\gamma_{0}) \cap E) = K$ holds.
\item A face $\gamma_{0} \preceq \gamma$ is an $\mathfrak{F}$-face
if and only if the number of pairs $e_{i},e_{j} \in \gamma_{0}$
with $\alpha_{ij} \ne 0$ differs from one.
\end{enumerate}
\end{proposition}

\proof
Since $q$ is of full rank, the origin
is the only singular point of 
$\b{X} = \Spec(R) \subset \KK^r$.
Since  $\mathcal{O}(X) = \KK$ holds, 
the open set 
$\rq{X} \subset \b{X}$ 
defined in~\ref{Xhutdef} does not 
contain the origin. 
Thus, the first assertion follows from 
Proposition~\ref{smoothchar}.
The second assertion is an immediate consequence 
of Proposition~\ref{hypersfrel}. 
\endproof

\begin{proposition}\label{rankest}
Let $X$ be a full intrinsic quadric with total
coordinate ring $\mathcal{R}(X)$. 
Then we have 
$$ 
\rank (\Cl(X)) \le \dim(X) + 3,
\qquad
\dim (\mathcal{R}(X)) \le 2 \dim(X) + 3.
$$
\end{proposition}

\proof
We may work in the above setup, 
that means that we have
$X = X(R,\mathfrak{F},\Phi)$, where
$R \cong \KK[T_{1}, \dots, T_{r}]/\bangle{q}$ 
with a homogeneous quadric $q$ of full rank, etc..
Recall that, counted with multiplicities, we have
precisely $r$ degrees $w_i = \deg(T_i) \in K$.
 
Using Proposition~\ref{quadricsetup}~(ii) and the fact
that $q$ is of rank $r$, we see that for any  
$w_{i}$ there is a counterpart $w_j$ such
that $w_{i} + w_{j} = \deg(q)$ holds. 
Since the $w_i$ generate $K \cong \Cl(X)$
as a lattice, this gives
$$
\rank(\Cl(X))
\; \le \;
\frac{r}{2} + 1.
$$
Now, by Theorem~\ref{firstresult},
the rank of $\Cl(X)$ equals
$\dim(\mathcal{R}(X)) - \dim(X)$.
Moreover, $R \cong \mathcal{R}(X)$ 
is of dimension $r-1$.
Together, this gives the assertion.
\endproof

We shall now study intrinsic quadrics having
small divisor class group.
In the statements, we visualize the occuring
$\mathfrak{F}$-bunches according to the 
discussion at the end of Section~\ref{bunches}.
The case of a one dimensional divisor class 
group is simple:

\begin{proposition}
\label{quadricalrankone}
Let $X$ be a full intrinsic quadric with 
$\Cl(X) \cong \ZZ$.
Then $X$ arises from a bunched ring 
$(R,\mathfrak{F},\Phi)$ 
with $\Phi$ given by
\begin{center}
\begin{picture}(0,0)%
\includegraphics{picdim1.pstex}%
\end{picture}%
\setlength{\unitlength}{1657sp}%
\begingroup\makeatletter\ifx\SetFigFont\undefined%
\gdef\SetFigFont#1#2#3#4#5{%
  \reset@font\fontsize{#1}{#2pt}%
  \fontfamily{#3}\fontseries{#4}\fontshape{#5}%
  \selectfont}%
\fi\endgroup%
\begin{picture}(6345,527)(889,-511)
\put(3601,-511){\makebox(0,0)[lb]{\smash{\SetFigFont{8}{9.6}{\familydefault}{\mddefault}{\updefault}
$(w_i,\mu_i)$
}}}
\put(5401,-511){\makebox(0,0)[lb]{\smash{\SetFigFont{8}{9.6}{\familydefault}{\mddefault}{\updefault}
$(w_{n-i},\mu_{n-i})$
}}}
\end{picture}

\end{center}
where $w_{1}, \ldots, w_{n}$ is a strictly increasing 
sequence of positive integers, such that $w_{i}+w_{n-i+1} =w$ 
holds, for some fixed $w \in \NN$, and
the weights $\mu_{i}$ satisfy 
$$
\mu_i \ge 1,
\qquad
\mu_{i} = \mu_{n-i},
\qquad 
\mu_1 + \ldots + \mu_n \ge 5.
$$
The variety $X$ is always 
$\QQ$-factorial, projective, 
and $\QQ$-Fano, 
and for its dimension, we have
$$
\dim(X)
\; = \;
\mu_1 + \ldots + \mu_n - 2.
$$  
\end{proposition}

Note that for a full intrinsic quadric 
$X = X(R,\mathfrak{F},\Phi)$ 
with one dimensional divisor class group,
the minimal toric ambient variety
$Z(R,\mathfrak{F},\Phi)$ as defined
in~\ref{minimalambient} is a 
weighted projective space.
Moreover, $X$ is smooth 
if and only if all numbers
$w_{i}$ equal one, and in this 
case it is merely a smooth projective
quadric.

An explicit example of a smooth projective 
full intrinsic quadric with divisor class group one 
is the Grassmannian $G(2,4)$,
studied in Examples~\ref{gzwovierI},~\ref{gzwovierIII}, 
and~\ref{gzwovierIV}. 
One may also consider ``weighted'' 
versions of this Grassmannian;
this provides singular varieties.

\begin{example}
Let $R := \KK[T_1,\ldots,T_6]/\bangle{T_1T_6-T_2T_5+T_3T_4}$,
and let $\mathfrak{F}$ consist of the classes of the $T_i$. 
Consider the $\mathfrak{F}$-bunch $\Phi$ as in 
Proposition~\ref{quadricalrankone}
defined by the weights $w_i := i$ 
and the multiplicities $\mu_i = 1$, 
where $i=1,\ldots,6$.

Then the corresponding variety $X := X(R,\mathfrak{F},\Phi)$
is not smooth, its Picard group is of index~$60$
in its divisor class group $\Cl(X) = \ZZ$, and its
canonical divisor generates a sublattice of index $14$
in $\Cl(X)$. In particular, $X$ is not Fano.
\end{example}

Now we turn to the case of smooth 
full intrinsic quadrics
with divisor class group $\ZZ^{2}$.
We give the ``intrinsic quadrics version'' 
of Kleinschmidt's classification~\cite{Kl}
of smooth complete 
toric varieties with Picard number two,
compare also~\cite[Prop.~11.1]{BeHa2}.

\begin{theorem}
\label{quadrklein}
Let $X$ be a smooth full intrinsic quadric 
with $\Cl(X) \cong \ZZ^{2}$.
Then $X$ arises from a bunched ring 
$(R,\mathfrak{F},\Phi)$ with an 
$\mathfrak{F}$-bunch
$\Phi$ given by one of the following figures:
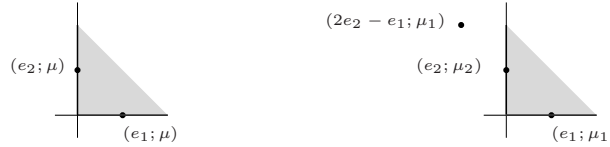
\begin{center}
\begin{picture}(0,0)%
\includegraphics{picdim2.pstex}%
\end{picture}%
\setlength{\unitlength}{1243sp}%
\begingroup\makeatletter\ifx\SetFigFont\undefined%
\gdef\SetFigFont#1#2#3#4#5{%
  \reset@font\fontsize{#1}{#2pt}%
  \fontfamily{#3}\fontseries{#4}\fontshape{#5}%
  \selectfont}%
\fi\endgroup%
\begin{picture}(11733,2793)(901,-2842)
\put(901,-1411){\makebox(0,0)[lb]{\smash{\SetFigFont{6}{7.2}{\familydefault}{\mddefault}{\updefault}$(e_2;\mu)$}}}
\put(3151,-2761){\makebox(0,0)[lb]{\smash{\SetFigFont{6}{7.2}{\familydefault}{\mddefault}{\updefault}$(e_1;\mu)$}}}
\put(7201,-511){\makebox(0,0)[lb]{\smash{\SetFigFont{6}{7.2}{\familydefault}{\mddefault}{\updefault}$(2e_2-e_1;\mu_1)$}}}
\put(11701,-2761){\makebox(0,0)[lb]{\smash{\SetFigFont{6}{7.2}{\familydefault}{\mddefault}{\updefault}$(e_1;\mu_1)$}}}
\put(9001,-1411){\makebox(0,0)[lb]{\smash{\SetFigFont{6}{7.2}{\familydefault}{\mddefault}{\updefault}$(e_2;\mu_2)$}}}
\end{picture}

\end{center}
where in the left hand side case,
$\mu \ge 3$ holds, 
and in the right hand side case, 
one has $\mu_i \ge 1$ and
$2\mu_1+\mu_2 \ge 5$.
Any such $X$ is projective, and its 
dimension is given by
$$
\dim(X)
=
2\mu-3,
\quad
\text{or}
\quad
\dim(X)
= 
2\mu_1+\mu_2-3,
$$
where the first equation corresponds to 
the l.h.s. case, and the second one to
the r.h.s. case.
The variety $X$ is Fano if 
and only if $\Phi$ 
belongs to the l.h.s. case. 
Moreover, different figures define
non-isomorphic varieties.
\end{theorem}

The proof is presented in the next section;
it is an interplay of the combinatorics of
bunches, and the symmetry condition~\ref{quadricsetup}~(ii). 
For divisor class groups of higher rank and/or equations
of higher degree, this type of classification 
problems seems to pose new and elementary 
but nontrivial combinatorial challenges.

\section{Intrinsic quadrics II}

In this section we prove Theorem~\ref{quadrklein}.
The setup is the following:
$X$ is a smooth full intrinsic quadric with
divisor class
group of rank two, arising from
a bunched ring $(R,\mathfrak{F},\Phi)$. 
According to Remark~\ref{quadricsetup},
the defining quadric then 
can be written as 
$$q = \sum_{i \le j} \alpha_{ij} T_{i}T_{j}.$$
As usual, $(E \topto{Q}K, \gamma)$
is the projected cone associated to $\mathfrak{F}$;
so, the lattice $K$ is of rank two.
We set $w := \deg(q)$, and
by $e_{1}, \ldots, e_{r}$ we denote the primitive 
generators of $\gamma$.
We call the images $Q(e_i)$ the {\em weight vectors}. 
Here comes a short list of arguments used frequently 
in the proof.

\begin{lemma}\label{multfface}
In the above situation, we have the following statements.
\begin{enumerate}
\item Let $e_i$ and $e_j$ such that $2Q(e_i)$, $2Q(e_j)$ 
and  $Q(e_i) + Q(e_j)$ 
are all different from $w = \deg(q)$. 
Then $\cone(e_i,e_j)$ is an $\mathfrak{F}$-face. 
\item Let $e_i,e_j$ and $e_k$ be pairwise different such that 
$Q(e_j)$ and $Q(e_i)$ lie on a common ray and $2Q(e_{k}) \ne w$ holds. 
Then $Q(e_i)$ and $Q(e_k)$ generate a projected $\mathfrak{F}$-face.
\item Let $e_{i}$ and $e_{k}$ be such that $2Q(e_{i}) = w$ and
$Q(e_{k}) \notin \QQ_{\ge 0} Q(e_{i})$ hold.
Then, $Q(e_{i})$ and $Q(e_{k})$ generate a projected $\mathfrak{F}$-face
if and 
only if there is a $j \ne i$ with $Q(e_{j}) \in \QQ_{\ge 0} Q(e_{i})$.
\item For every index $i$ there exists an index $j$, not necessarily 
different from $i$, such that $Q(e_i) + Q(e_j) = w$ holds.
\item If $\cone(e_i; \; i \in I)$ is a relevant $\mathfrak{F}$-face,
then the images $Q(e_i)$, $i \in I$, generate the lattice.
\end{enumerate}
\end{lemma}

\proof
For~(i), just note that in this case $\alpha_{ii}$, $\alpha_{ij}$ 
and $\alpha_{jj}$ vanish.
Thus, the assertion follows from Proposition~\ref{quadricprop}~(ii).

For (ii), note that $\alpha_{kk} = 0$ holds. To define
an $\mathfrak{F}$-face $\gamma_0 \preceq \gamma$ projecting
onto $\cone(Q(e_i),Q(e_k))$, we go through the possible cases
and apply  Proposition~\ref{quadricprop}~(ii):
$$
\begin{array}{lllll}
\alpha_{ii} = 0, & 
\alpha_{ik} = 0: & 
& 
& 
\gamma_0 :=  \cone(e_i,e_k), \\
\hfill \text{---} \hfill, & 
\alpha_{ik} \ne 0, & 
\alpha_{ij} \ne 0: & 
&
\gamma_0 :=  \cone(e_i,e_j,e_k), \\
\hfill \text{---} \hfill, & 
\hfill \text{---} \hfill, & 
\alpha_{jk} \ne 0: & 
&
\gamma_0 :=  \cone(e_i,e_j,e_k), \\
\hfill \text{---} \hfill,  & 
\hfill \text{---} \hfill,  &
\alpha_{ij} = \alpha_{jk} = 0, &
\alpha_{jj} = 0: & 
\gamma_0 :=  \cone(e_j,e_k), \\
\hfill \text{---} \hfill,  & 
\hfill \text{---} \hfill,  &
\hfill \text{---} \hfill, &
\alpha_{jj} \ne 0: & 
\gamma_0 :=  \cone(e_i,e_j,e_k). \\
\end{array}
$$
By symmetry, this settles as well the case $\alpha_{jj} = 0$.
If both $\alpha_{ii}$ and $\alpha_{jj}$ differ from zero,
then we can take $\gamma_0 := \cone(e_i,e_j,e_k)$.

We prove (iii).
The ``if'' part is a direct application of~(ii).
To see the ``only if'' part, 
assume that only $e_{i}$ is mapped into
$\QQ_{\ge 0} Q(e_{i})$.
Then, by regularity of $q$, we must have 
$\alpha_{ii} \ne 0$. 

Moreover, for all $j,l$ different from $i$,
with $Q(e_j)$ and $Q(e_l)$ lying in the 
cone generated by $Q(e_i)$ and $Q(e_k)$,
we have $Q(e_j) + Q(e_l) \ne w$, 
and hence $\alpha_{jl} = 0$.
Thus, Proposition~\ref{quadricprop}~(ii)
tells us that there is no
$\mathfrak{F}$-face $\gamma_0 \preceq \gamma$
with $Q(\gamma_0) = \cone(Q(e_i),Q(e_k))$.
Hence the latter cone is not a projected 
$\mathfrak{F}$-face.

We turn to~(iv). Suppose that some index $i$ satisfies 
$Q(e_i) + Q(e_j) \ne w$
for all indices~$j$.
Then this implies $\alpha_{ij}=0$ for all $j$.
But this contradicts the regularity of the defining 
quadric~$q$.

Finally, assertion~(v) is a direct consequence of smoothness of $X$
and Proposition~\ref{quadricprop}~(i).
\endproof

\proof[Proof of Theorem~\ref{quadrklein}]
We have to show that $\Phi$ admits a representation
as in the figures shown in~\ref{quadrklein}.
First of all, as $X$ admits only constant functions, 
Proposition~\ref{globfunconst} shows
that $\vartheta := Q(\gamma)$ is strictly convex,
and that all weight vectors are different from zero.

\medskip\noindent
{\em Step 1}. The $\mathfrak{F}$-bunch $\Phi$ consists 
of a single cone $\tau$.
\medskip

Assume the converse. Then we can choose two
different cones $\tau_1, \tau_2 \in \Phi$. 
Both are images of relevant faces, and thus 
they are two-dimensional.
By the defining properties of an 
$\mathfrak{F}$-bunch, we have 
$\tau_1^{\circ} \cap \tau_2^{\circ} \ne \emptyset$.
Hence, $\tau_{12} := \tau_1 \cap \tau_2$ is
of dimension two.
Moreover, $\tau_{12} \ne \tau_i$ holds.
Thus, there are pairwise different weight vectors 
$x_1$, $x_2$, $y_1$ and $y_2$ such that
$$
\tau_1 = \cone(x_1,y_1),
\qquad
\tau_2 = \cone(x_2,y_2),
\qquad
\tau_{12} = \cone(x_2,y_1).
$$

We claim that $\tau_{12}$ is not a projected $\mathfrak{F}$-face.
Indeed, for any given $\sigma \in \Phi$, we have 
$\sigma^\circ \cap \tau_i^{\circ} \ne \emptyset$.
Using the fact that $\sigma$ and the $\tau_i$ are contained
in the strictly convex cone $\vartheta = Q(\gamma)$,
we obtain $\sigma^{\circ} \cap \tau_{12}^{\circ} \ne \emptyset$.
Thus, if $\tau_{12}$ were a projected $\mathfrak{F}$-face,
it would contain an element of $\Phi$. 
This contradicts $\tau_{12} \subsetneq \tau_i \in \Phi$.

Next we claim that there exist two-dimensional $\mathfrak{F}$-faces
$\gamma_i \preceq \gamma$ with $Q(\gamma_i) = \tau_i$.
For this, we first exclude $2x_2 = w$: otherwise, 
since $\tau_2$ is a projected $\mathfrak{F}$-face,
Lemma~\ref{multfface}~(iii) provides a weight vector
$x_2' \in \QQ_{\ge 0} x_2$.
But then, using again Lemma~\ref{multfface}~(iii), we
see that $\tau_{12}$ is a projected $\mathfrak{F}$-face.
This contradicts the previous claim.
Similarly, one excludes $2y_1 =w$.

Since $\tau_{12}$ is not a projected $\mathfrak{F}$-face, 
Lemma~\ref{multfface}~(i) gives $x_2+y_1 = w$.
As an immediate consequence, we obtain $x_i + y_i \ne w$ 
and $2y_2 \ne w \ne 2x_1$.
Hence, Lemma~\ref{multfface}~(i) tells us that
the $\tau_i$ are images of two-dimensional $\mathfrak{F}$-faces
$\gamma_i$ respectively, and the claim is verified.

As remarked before, the last claim implies 
that the $\tau_i$ are regular, and that we may assume 
$x_{1} = (1,0)$ and $y_{1} = (0,1)$. 
Write $x_{2} = (a,b)$ with $a, b > 0$.
Then, by strict convexity of $\vartheta = Q(\gamma)$,
we have $y_{2} = (-c ,d)$ with $c, d > 0$. 
The regularity of $\tau_{2}$
implies
$$ 1 = ad + bc \ge 2. $$
But this is impossible. 
Thus, the $\mathfrak{F}$-bunch $\Phi$ cannot contain 
more than one cone.
This ends the proof of Step~1.

\medskip

In the sequel, we denote by $\tau$ the unique cone of $\Phi$,
and we denote by $x$ and $y$ the weight vectors of minimal 
length satisfying $\tau = \cone(x,y)$.

\medskip\noindent
{\em Step~2.} If $v \in \tau^{\circ}$ is a weight vector,
then it satisfies $2v = w$.
\medskip

Assume that $2v \neq w$ holds. 
Then we can exclude $2x = w$.
If so, the fact that $\tau$ is a projected $\mathfrak{F}$-face
and Lemma~\ref{multfface}~(iii) provide
$e_i \ne e_j$, both mapping into $\QQ_{\ge 0} x$.
Again Lemma~\ref{multfface}~(iii) shows
that $\cone(x,v)$ is a projected $\mathfrak{F}$-face. 
This contradicts the defining property of the 
$\mathfrak{F}$-bunch $\Phi = \{\tau\}$.
Similarly we exclude $2y = w$. 

Since we cannot have $x + v = w$ and $y + v = w$
simultaneously, 
Lemma~\ref{multfface}~(i) shows that at least one of
$\cone(x,v)$ and $\cone(v,y)$ is a projected 
$\mathfrak{F}$-face. 
As before, this contradicts the defining property of the 
$\mathfrak{F}$-bunch $\Phi = \{\tau\}$.  
Thus, we proved that any weight vector $v$ contained in
$\tau^{\circ}$ must fulfill $2v = w$.

\medskip\noindent
{\em Step~3.} 
The cone $\tau$ is regular.
\medskip

Suppose that $\tau$ is not regular.
Since $\tau$ is a projected $\mathfrak{F}$-face,
the lattice $K$ is generated by certain weight 
vectors lying in $\tau$.
By non-regularity, at least one of these weight 
vectors, say $v$, has to lie in $\tau^{\circ}$.
By the previous step, there is only one
possibility: $v$ is determined by $2v =w$.

Consider any pair $x'$, $y'$ of weight
vectors generating the cone $\tau$.
Then $2x' \ne w$ and $2y' \ne w$ hold. 
Using this, we can conclude $x' + y' = w$,
because otherwise, Lemma~\ref{multfface}~(i) and~(v)
would show that $\tau$ is regular.
This implies that $x$, $y$, and $v$ are the 
only weight vectors in $\tau$.

Consequently, $x,y,v$ form a system of generators 
of $K$ satisfying $x+y=2v$. 
But then $x,v$ as well as $v,y$ are bases.
Thus, we may prescribe $x$, $v$ and then
determine $w$ and $y$ as follows:
$$
x = (1,0),
\qquad
v = (1,1),
\qquad
w = (2,2),
\qquad
y = (1,2).
$$ 

Now, suppose that there exists a weight vector $u$ different from
$x$, $y$ and $v$. 
Then, because of Lemma~\ref{multfface}~(iv) and $w=2v$, 
there must exist a weight vector $u' \ne u$ with $u+u' = w$.
Since $u$ as well as $u'$ do not belong 
to $\tau$, and $u+u' \in \tau^{\circ}$ 
holds, we have $\tau \subset \cone(u,u')$.
By possibly interchanging $u$ and $u'$, we achieve that 
$\cone(x,u)$ contains $\tau$. 

Lemma~\ref{multfface}~(i) and~(v) say that 
$x$ and $u$ generate the lattice.
Consequently, $u = (a,1)$ for some integer $a$.
Note that $a \le 0$ holds, because of $u \not\in \tau$.
But this leads to a contradiction:
$$ 
u' 
\; = \; 
w - u
\; = \; 
(2-a,1)
\; \in \;
\tau.
$$

So, we arrived at the following picture:
the only weight vectors are
$x$, $v$ and $y$, 
and we have $\tau = \cone(x,y)$.
The multiplicities must satisfy 
$\mu(x) =  \mu(y)$ and $\mu(v) = 1$.
The latter is, once more, a consequence
of Lemma~\ref{multfface}~(iii): 
a higher multiplicity $\mu(v)$ would
imply that $\cone(x,v)$ and $\cone(v,y)$
are projected $\mathfrak{F}$-faces, 
which contradicts
$\Phi=\{\tau\}$.

Since $q$ has at least 5 variables,
we see that $\mu(x) = \mu(y) \ge 2$  
holds. 
Having in mind that $q$ is a regular 
quadric, this implies that there are at 
least two different pairs $i,j$ and $i',j'$
with $Q(e_i) = Q(e_{i'}) = x$ and
$Q(e_j) = Q(e_{j'}) = y$
such that $\alpha_{ij} \ne 0 \ne \alpha_{i'j'}$
holds.
Thus, Proposition~\ref{quadricprop} 
and Lemma~\ref{multfface}~(v) tell us that
$x$ and $y$ generate the lattice.
A contradiction.

\medskip

\medskip
\noindent
{\em Step~4.} The weight vectors $x$ and $y$ form a lattice basis. 
\medskip

First consider the case that there is a weight vector 
$v \in \tau^{\circ}$.
Then we know by step~2 that $2v = w$ holds.
Consider any pair $x'$, $y'$ of weight vectors 
generating $\tau$ as a convex cone.
If we have $x'+y' \ne w$, then Lemma~\ref{multfface}~(i) and~(v) 
shows that $x'$ and $y'$ generate the lattice,
hence $x' = x$ and $y' = y$, and we are done.

So, we may assume $x'+ y' = w$. But this
fixes $x'$ and $y'$. Thus, we have again
that $x$ and $y$ are the only weight vectors
lying in the boundary of $\tau$.
Lemma~\ref{multfface}~(v) says that $x$, $y$ and $v$
generate the lattice.

Because of $x+y=2v$, the pair $x$, $v$ is a basis.
Thus, we may assume that $x=(1,0)$ and $v=(1,1)$ holds.
This implies 
$$ y = w-x = 2v-x = (1,2). $$
So, we arrived at a contradiction to regularity of
$\tau$.
That means $x+y=w$ cannot happen, 
and the case that there exists a weight vector
$v \in \tau^{\circ}$ is settled.

Let us consider the case that all weight vectors of 
$\tau$ lie in the boundary.
If $x$ and $y$ are the only ones, then we succeed
with Lemma~\ref{multfface}~(v).
So, suppose that on some ray, there are at least
two different weight vectors, 
say $y', y'' \in \QQ_{\ge 0}y$.

Then, if $2y' = w$, then 
Lemma~\ref{multfface}~(i) and~(v) show that
$x$ and $y''$ generate the lattice.
This implies $y'' = y$, and we are done.
Similarly, we settle the case $2y'' = w$.
Next, if $2x$, $2y'$ and $2y''$ all differ
from $w$, then 
Lemma~\ref{multfface}~(i) and~(v) show
that either the pair $x,y'$ or the pair
$x,y''$ is a lattice basis.

Finally, suppose that $2x = w$ holds.
Then, by a suitable choice of coordinates,
we may achieve that 
$$
x = (a,0),
\qquad
w = (2a,0),
\qquad
y' = (0,b'),
\qquad
y'' = (0,b'')
$$
hold with positive integers $a$, and $b' \ne b''$.
According to Lemma~\ref{multfface}~(iv), 
we have weight vectors
$$
z' = w-y' = (2a,-b'),
\qquad
z'' = w-y'' = (2a,-b'').
$$
By Lemma~\ref{multfface}~(i) and~(v), the vectors
$y'$ and $z''$ generate the lattice.
This implies $2ab' = \pm 1$,
which is impossible.

\medskip

So, having passed step 4, 
we know that $x$ and $y$ form a 
lattice basis.
According to Lemma~\ref{multfface}~(iv),
there exist weight vectors $x'$ and 
$y'$ with
$
x + x'
 = 
y + y' 
= 
w
$.

\medskip\goodbreak
\noindent
{\em Step~5.} The case $x' \neq x$ and $y \neq y'$. 
\medskip

We shall show that the present situation
amounts to the left hand side figure of Theorem~\ref{quadrklein}.
Note that we have
$$
2x \ne w, 
\qquad
2x' \ne w, 
\qquad
2y \ne w, 
\qquad
2y' \ne w, 
\qquad
x + y' \ne w, 
\qquad
x' + y \ne w.
$$
In particular, using step~2, we can conclude from 
this list that neither $x'$ nor $y'$ belong to 
$\tau^{\circ}$.

We claim that $\tau \subset \cone(x',y')$ holds.
For the verification of this claim, we may assume 
that $x=(1,0)$ and $y=(0,1)$ hold. 
Since all weight vectors belong to the 
strictly convex cone $\vartheta = Q(\gamma)$, 
we see that $\tau \not\subset \cone(x',y')$ can 
only happen if $x'$ and $y'$ both lie 
either in the upper or the right half plane.

By symmetry, it suffices to treat the 
case that $x'$ and $y'$ 
both lie in the upper half plane.
There are integers $a$ and $b$ such that
$$
x' = (a,b),
\qquad 
w = x + x' = (a+1,b),
\qquad
y' = w-y  = (a+1,b-1).
$$
Note that $b \ge 1$ holds because $y'$ belongs 
to the upper half plane. 
Moreover, we have $a \le 0$ because $x'$
does not belong to $\tau^{\circ}$.

We settle the case $b=1$. 
Strict convexity of $\vartheta$
and $x,y' \in \vartheta$ imply $a \ge -1$.
If $a=0$, then $x'=y$
and $y'=x$ holds, and we are done.
If $a =-1$ holds, then we obtain $y'=0$,
which contradicts
the characterization~\ref{globfunconst} of 
$\mathcal{O}(X) = \KK$.
So we may assume $b \ge 2$ from now on.

Since $2x$, $2y'$ and $x+y'$ all differ from $w$,
Lemma~\ref{multfface}~(i) yields that $\cone(x,y')$ 
is the image of a twodimensional $\mathfrak{F}$-face.
By the defining properties of $\Phi$,
we can conclude that $\cone(x,y')$ contains $\tau$.
This means $a \le -1$. 
Moreover, by Lemma~\ref{multfface}~(v), the vectors
$x$ and $y'$ generate the lattice. 
In particular, $b=2$ holds.
Summing up, we obtained
$$
x' = (a,2),
\qquad 
w =  (a+1,2),
\qquad
y' = (a+1,1),
\qquad
a \le  -1.
$$

In view of step~2, this excludes in particular,
the existence of weight vectors inside $\tau^{\circ}$.
We can even exclude existence of weight vectors different
from $x$ and $y$ in the whole right half plane:
Such a vector is of the form 
$u = (\zeta,\eta)$ with $\zeta \ge 0$ and 
$\eta \le 0$. 
By Lemma~\ref{multfface}~(iv), we have a further
weight vector
$$
u' = w -u = (a+1 - \zeta, 2 - \eta).
$$

Since $2x$, $2u'$ and $x+u'$ all differ from $w$,
we may apply Lemma~\ref{multfface}~(i) and~(v), 
and see that $x$ and $u'$ generate the lattice $K$.
This contradics that fact that the second coordinate
of $u'$ is at least two.
So, there are no weight vectors different from $x$
in the fourth quadrant.

Using the second property of Definition~\ref{fdefs}~(ii), 
we obtain that there are at least two different 
$e_i, e_j \in \gamma$ mapping into the ray through $x$.
Thus, applying Lemma~\ref{multfface}~(ii) to $e_i$, $e_j$
and any $e_k$ mapping to $x'$ shows that 
$x'$ and the images of $e_i$ and $e_j$ generate
the lattice, which is a contradiction.

This ends the verification of the first claim of step 5, 
and we know now that $\tau \subset \cone(x',y')$ holds.

Next we claim that $x'$ and $y'$ form a lattice basis.
Indeed, if $x' + y' = w$ holds, then we have 
$x' = y$ and $y' = x$, and we are done.
If $x' + y' \ne w$ holds, then Lemma~\ref{multfface}~(i)
and~(v) tell us that $x'$ and $y'$
form a lattice basis, 
and the claim is verified.

By a suitable change of coordinates,
we achieve that $x' = (0,1)$ and $y'=(1,0)$ hold.  
Write $x = (k,l-1)$. Then $y = (k-1,l)$ 
holds with $k,l > 0$.  
The fact that $x$ and $y$ form a lattice
basis translates to
$$ 
\pm 1 
\; = \; 
kl - (k-1)(l-1) 
\; = \; 
k+l-1.
$$ 
Since $k,l$ are positive integers, we
eventually obtain $k = l = 1$ and therefore 
$x = y'$ and $x' = y$. 
Thus, we have $w=(1,1)$.
According to step 2, this shows in particular,
that $\tau^{\circ}$ contains no weight vectors.

In order to arrive at the left hand side figure 
of Theorem~\ref{quadrklein}, we still have to
show that $x$ and $y$ are the only weight
vectors.

Suppose that there is a weight vector $u$,
which differs from $x$ and from $y$.
By Lemma~\ref{multfface}~(iv), we then find 
a weight vector $u'$ with 
$u+u'=w=(1,1)$.
Note that neither $u$ nor $u'$ belong to 
$\tau$. Moreover, they lie in different
quadrants, because of $u+u' = w \in \tau^{\circ}$.

By possibly interchanging
$u$ and $u'$, we achieve that
$u = (a,-b)$ holds 
with positive integers $a$ and $b$. 
Now, Lemma~\ref{multfface}~(i) applies to
$\cone(u,y)$. It shows that $u$ and $y$
form a lattice basis. 
Thus, we obtain $a=1$.
This leads to a contradiction:
$$
u'
\; = \;
w-u
\; = \; 
(0,1+b)
\; \in \; 
\tau.
$$

\medskip
\noindent
{\em Step~6.} The case $x=x'$ or $y=y'$.
\medskip

We treat as an example the case $y = y'$.
Note that then $2y=w$ holds.
Thus, according to step~2, there are no weight vectors
in $\tau^{\circ}$.
By a suitable choice of coordinates
we achieve
$$
x = (1,0),
\qquad
y = (0,1),
\qquad
w = (0,2),
\qquad
x' = (-1,2).
$$  

We show that we are in the situation of the
right hand side figure.
That means that we have to exclude existence of
further weight vectors.
Assume to the contrary that there is a weight 
vector $u$, which is different from $x$, $x'$ 
and $y$. 

Lemma~\ref{multfface}~(iv) provides a 
counterpart $u'$ satisfying $u+u' = w$.
Since all weight vectors lie in the 
strictly convex cone $\vartheta=Q(\gamma)$,
neither $u$ nor $u'$ can belong to the negative 
quadrant.
Consequently, after possibly interchanging 
$u$ and $u'$, there are nonnegative integers 
$a$ and $b$ with
$$  
u = (a,-b),
\qquad
u' = w-u = (-a,2+b).
$$

Since $2x$, $2u'$ and $x+u'$ differ from $w$,
Lemma~\ref{multfface}~(i) and~(v) show that 
$u'$ and $x$ generate the lattice.
This contradicts the explicit presentation of
$u'$. Thus, step~6 is done.

\medskip

So, we proved that the $\mathfrak{F}$-bunches
corresponding to the smooth intrinsic quadrics
$X$ with $\mathcal{O}(X) = \KK$
arise from one of the figures listed in the 
assertion.
Quasiprojectivity and the statements on 
being Fano are due to our general results.
Moreover, Proposition~\ref{quasiproj2proj} 
gives completeness of the varieties in question.

In order to see that different figures belong 
to non-isomorphic varieties, note first that
the left and the right hand side have nothing to do
with each other, because the one are Fano and 
the other not.
Moreover, different choices of the weights $\mu$ 
on the left hand side give rise to varieties 
of different dimension.

The right hand side is a little more subtle.
First note that the occuring weight vectors form
a Hilbert basis for the cone in $K = \Cl(X)$ 
generated by the degrees that admit nontrivial 
homogeneous elements in $R = \mathcal{R}(X)$.
This Hilbert basis as well as the dimensions
of the corresponding spaces of homogeneous 
sections are invariants of $X$.
But the latter determine  
$\mu_{1}$ and $\mu_{2}$.
\endproof

\section{Surfaces}
In this section we focus on the case of surfaces. 
We first present some general results,
and then turn to explicit examples.
The first observation is that, in our setup,
complete surfaces are 
determined by their total coordinate ring.

\begin{proposition}\label{isoring2isosurf} 
Let $X_1$ and $X_2$ be normal $A_2$-maximal surfaces 
with finitely generated free divisor class groups 
and finitely generated total coordinate rings 
$R_1$ and $R_2$ respectively, and suppose that
$X_2$ is complete.
Then we have $X_1 \cong X_2$ if and only if 
$R_1 \cong R_2$ as graded rings.
\end{proposition}

\proof
Only for the ``if'' part there is something to show.
By Theorem~\ref{firstresult}, we may assume that 
$X_2$ arises from a bunched ring 
$(R_2,\mathfrak{F}_2,\Phi_2)$.
Let $\alpha \colon R_2 \to R_1$ be a graded isomorphism,
and set $\mathfrak{F}_1 := \alpha(\mathfrak{F}_2)$.
Then, by Assertion~\ref{firstresult}~a),
we may assume $X_1  = X(R_1, \mathfrak{F}_1, \Phi_1)$
with an $\mathfrak{F}_1$-bunch $\Phi_1$.

Now, look at $\b{X}_i := \Spec(R_i)$, and the isomorphism
$\b{\varphi} \colon \b{X}_1 \to \b{X}_2$ defined by 
$\alpha \colon R_2 \to R_1$.
Moreover, consider the open subsets 
$\rq{X}_i \subset \b{X}_i$,
the quotient presentations $p_i \colon \rq{X}_i \to X_i$ 
of Proposition~\ref{xasgoodquot}, 
and the open sets $W_i \subset X_i$ as studied in 
Lemma~\ref{subsetWX}.

Then $\b{\varphi}$ maps $p_1^{-1}(W_{1})$ onto 
$p_2^{-1}(W_{2})$, and thus induces an isomorphism
$W_1 \to W_2$.
Since the surface $X_2$ is complete,
and the complements $X_i \setminus W_i$ 
are of codimension at least two and 
hence at most finite, 
the map $W_1 \to W_2$ extends to an 
open embedding $X_1 \to X_2$. 
Since $X_1$ is $A_2$-maximal,
$X_1 \to X_2$ is onto.
\endproof

As an application of the independence of a complete 
surface $X(R,\mathfrak{F},\Phi)$ from the particular 
choice of $\mathfrak{F}$ and $\Phi$, we obtain a general 
projectivity criterion:

\begin{proposition}
\label{compl2proj}
Let $X$ be a normal $A_{2}$-maximal surface with 
finitely generated free divisor class group and
finitely generated total coordinate ring. 
If $\mathcal{O}(X) = \KK$ holds, 
then $X$ is projective.
\end{proposition}

\proof 
By Theorem~\ref{firstresult},
we may assume that $X = X(R,\mathfrak{F},\Phi)$
holds with a bunched ring $(R,\mathfrak{F},\Phi)$.
As in Lemma~\ref{Phi2Theta}, we extend $\Phi$ 
to a bunch $\Theta$ in the projected cone 
$(E \topto{Q} K,\gamma)$ 
associated to $\mathfrak{F}$.
Consider the dual projected cone
$(F \topto{P} N,\delta)$,
the maximal projectable fan $\rq{\Delta}$ 
in $F$ correponding to $\Theta$, 
and its quotient fan
$\Delta$ in $N$, see~\cite[Sec.~4]{BeHa2}.

According to Proposition~\ref{globfunconst}, 
the rays of $\Delta$ generate the vector space 
$N_\QQ$ as a cone.
Consequently, Lemma~\ref{einsgeruest} provides 
a polytopal fan $\Delta'$ in $N$ having 
the same rays as $\Delta$.
The fan $\rq{\Delta}'$ in $F$ generated
by those faces of $\delta$ that map onto a cone 
of $\Delta'$ is maximal projectable, 
and hence $\rq{\Delta}'$ 
corresponds to a bunch $\Theta'$ in
$(E \topto{Q} K,\gamma)$, 
see again~\cite[Sec.~4]{BeHa2}.

Picking out the $\mathfrak{F}$-faces from $\rel(\Theta')$,
and taking the minimal cones among their images under $Q$, 
we arrive at an $\mathfrak{F}$-bunch $\Phi'$, and hence
a bunched ring $(R,\mathfrak{F}, \Phi')$.
The associated variety $X'$ is by construction 
embedded into $Z_{\Theta'}$, and hence is projective.
Proposition~\ref{isoring2isosurf} gives $X \cong X'$.
\endproof

We give some applications of this proposition. 
Goodman showed in~\cite{Go}
that a normal complete surface is projective,
if and only if the set of its non-factorial
singularities admits a common affine neighbourhood.
For surfaces with a finitely generated total coordinate
ring this can be sharpened:

\begin{corollary}
\label{surfproj}
A normal complete surface with 
finitely generated free divisor class group
and finitely generated
total coordinate ring is projective if and only if
any two of its non-factorial singularities admit a
common affine neighbourhood.
\end{corollary}

\proof
The ``if'' part is the interesting one. We have to show
that the surface, call it $X$, then is $A_2$. Let 
$x,x' \in X$. If both points are non-factorial, then
they have a common affine neighbourhood by assumption.
So, we may assume that $x$ is factorial. 

Let $X' \to X$ resolve all the non-factorial 
singularities different from $x'$. 
Then $X'$ is projective, 
hence $x$ and $x'$ have a common affine 
neighbourhood $U' \subset X'$. After suitably shrinking
$U'$, we achieve that it does not meet exceptional curves,
and hence maps onto an affine 
neighbourhood $U$ of $x$ and $x'$ in $X$.
\endproof

\begin{corollary}
If $X$ is a full intrinsic quadric of dimension two,
then $X$ is projective and $\Cl(X)$ is of rank at most $5$.
\end{corollary}

\proof 
The estimate for the rank of 
$\Cl(X)$ is a direct consequence 
of Proposition~\ref{rankest}, 
and projectivity of $X$ is due to
Proposition~\ref{compl2proj}.
\endproof

We now turn to explicit examples.
The first one is a quotient of an open
subset of the Grassmannian $G(2,4)$ 
by the action of a twodimensional
torus. Consider 
the ring 
$$
R 
\; := \; 
\KK[T_1, \ldots, T_6] / 
\bangle{ T_{1}T_{6} - T_{2}T_{5} + T_{3}T_{4} },
$$
and define a $\ZZ^3$-grading on $R$ by 
setting $\deg(T_i) := w_i$, where the 
vectors $w_i \in \ZZ^3$ are the following
$$
\begin{array}{lll} 
w_{1} := (1,0,1), 
& 
w_{2} := (1,1,1), 
& 
w_{3} := (0,1,1),
\\
w_{4} := (0,-1,1), 
&
w_{5} := (-1,-1,1),
& 
w_{6} := (-1,0,1).
\end{array}
$$
Note that in this setup,
the condition~\ref{quadricsetup}~(ii) is fulfilled.
By Proposition~\ref{isoring2isosurf},
these data already uniquely determine a surface,
and the results of the paper give insight
to the geometry:

\begin{theorem}
Let $X$ be an $A_2$-maximal normal variety having 
the $\ZZ^3$-graded ring $R$ as its total coordinate 
ring.
\begin{enumerate}
\item The variety $X$ is a $\QQ$-factorial, projective surface. 
It has precisely five singular points.
\item The Picard group $\Pic(X)$ is a sublattice of index 
72 in the divisor class group $\Cl(X) \cong \ZZ^{3}$.
\item The semiample cone $\SAmple(X)$
equals the moving cone $\Mov(X)$; it
is of full dimension in  $\Cl_{\QQ}(X)$, 
and has six extremal rays.
\item The canonical divisor of $X$ is of the form $4D$ with 
a Weil divisor $D$. The lowest locally principal multiple
is $12D$, and this is an ample divisor.   
\end{enumerate}
\end{theorem}

\proof 
We have to represent $X$ as the variety arising from 
a bunched ring $(R,\mathfrak{F},\Phi)$. 
As usual, let $\mathfrak{F}$ consist of the 
classes of the $T_i$ in $R$. 
Let $(E \topto{Q} K, \gamma)$ be the associated 
projected cone, i.e., we have $K = \ZZ^3$ and 
$E = \ZZ^6$, the map $Q$ sends the $i$-th canonical 
base vector $e_{i}$ to $w_{i}$,
and  $\gamma = \QQ^{6}_{\ge 0}$ is the positive 
orthant.

According to Assertion~\ref{firstresult}~a)
and 
Proposition~\ref{isoring2isosurf},
we may take any $\mathfrak{F}$-bunch~$\Phi$.
Write $\tau(i_{1}, \dots, i_{r})$ 
for the cone in $K_{\QQ}$
generated by the vectors
$w_{i_{1}}, \dots, w_{i_{r}}$, and
consider the cones
$$ 
\tau(1,3,5), 
\quad 
\tau(2,4,6),
\quad
\tau(1,6,2,5), 
\quad 
\tau(1,6,3,4), 
\quad 
\tau(2,5,3,4).
$$
As in Proposition~\ref{quadricprop}~(ii) we see that 
these cones form an  $\mathfrak{F}$-bunch~$\Phi$.
Indicating the two simplicial cones by dashed lines,
we obtain the following picture:

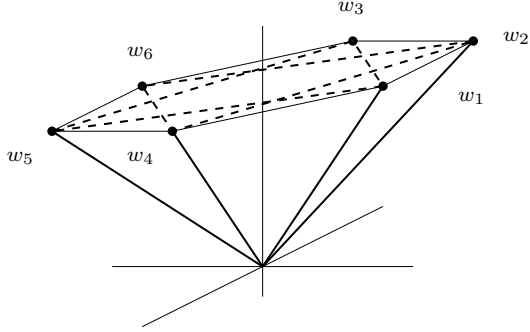
\begin{center}
\begin{picture}(0,0)%
\includegraphics{grassquot.pstex}%
\end{picture}%
\setlength{\unitlength}{1657sp}%
\begingroup\makeatletter\ifx\SetFigFont\undefined%
\gdef\SetFigFont#1#2#3#4#5{%
  \reset@font\fontsize{#1}{#2pt}%
  \fontfamily{#3}\fontseries{#4}\fontshape{#5}%
  \selectfont}%
\fi\endgroup%
\begin{picture}(7425,4929)(226,-4573)
\put(7651,-286){\makebox(0,0)[lb]{\smash{\SetFigFont{8}{9.6}{\familydefault}{\mddefault}{\updefault}{\color[rgb]{0,0,0}$w_2$}%
}}}
\put(2026,-511){\makebox(0,0)[lb]{\smash{\SetFigFont{8}{9.6}{\familydefault}{\mddefault}{\updefault}{\color[rgb]{0,0,0}$w_6$}%
}}}
\put(5176,164){\makebox(0,0)[lb]{\smash{\SetFigFont{8}{9.6}{\familydefault}{\mddefault}{\updefault}{\color[rgb]{0,0,0}$w_3$}%
}}}
\put(226,-2086){\makebox(0,0)[lb]{\smash{\SetFigFont{8}{9.6}{\familydefault}{\mddefault}{\updefault}{\color[rgb]{0,0,0}$w_5$}%
}}}
\put(2026,-2086){\makebox(0,0)[lb]{\smash{\SetFigFont{8}{9.6}{\familydefault}{\mddefault}{\updefault}{\color[rgb]{0,0,0}$w_4$}%
}}}
\put(6976,-1186){\makebox(0,0)[lb]{\smash{\SetFigFont{8}{9.6}{\familydefault}{\mddefault}{\updefault}{\color[rgb]{0,0,0}$w_1$}%
}}}
\end{picture}

\end{center}

Similarly as for the projected faces, we denote the 
faces of $\gamma$ by $\gamma(i_{1}, \dots, i_{r})$. 
Then the collection $\rel(\Phi)$ of relevant faces
consists of $\gamma$ itself, all the facets of $\gamma$, 
and the faces
\begin{equation}\label{eq:singfac}
\gamma(1,3,5),
\quad
\gamma(2,4,6),
\quad
\gamma(1,6,2,5),
\quad
\gamma(1,6,3,4),
\quad
\gamma(2,5,3,4).
\end{equation}

The images of the relevant cones are of
full dimension, 
hence $X$ is $\QQ$-factorial. 
Moreover, all strata associated to the faces of dimension 
five or more consist of smooth points, 
whereas the strata associated to the faces listed 
in (\ref{eq:singfac}) are (isolated) singularities.
This proves~(i).

In order to prove the second assertion,
recall from Corollary~\ref{picarddescr} that the 
Picard group of $X$ is given by
$$ 
\Pic(X) 
\; = \; 
\bigcap_{\gamma_{0} \in \cov(\Phi)} 
        Q(\lin(\gamma_{0}) \cap E).
$$ 
Since the images of the faces of dimension five or six 
generate $K$, 
we may restrict the intersection to the faces 
listed  in (\ref{eq:singfac}).
Thus, we obtain that $\Pic(X)$ is generated as a 
sublattice of $K$ by the vectors
$$
(2,4,0), 
\qquad
(0,6,0),
\qquad
(0,0,6).
$$ 

We turn to (iii).
Theorem~\ref{amplecone} gives us the semiample cone of $X$.
Namely, we have to intersect the members of $\Phi$.
Explicitly, this means that $\SAmple(X)$
has the primitive generators 
$$
(1,-1,3),
\quad
(-1,-2,3),
\quad 
(-1,1,3),
\quad 
(-2,-1,3),
\quad 
(2,1,3),
\quad
(1,2,3).
$$ 

According to Proposition~\ref{candivdescr},
the canonical divisor class of $X$ is represented
by the vector $(0,0,-4)$. 
In particular, $X$ is a $\QQ$-Fano variety. 
More precisely, from our description of the Picard group
of $X$, we see that $(0,0,-12)$ is the least multiple,
which is Cartier.
\endproof

Finally, we consider a surface studied by 
Hassett and Tschinkel in~\cite{Tsch},
namely the minimal resolution $X$ of 
the $E_{6}$ cubic surface 
$$
X' 
\; := \; 
\{[z_{0}, \ldots, z_{3}] \in \PP^{3}; \; z_{1}z_{2}^{2}+
z_{2}z_{0}^{2} + z_{3}^{3} = 0\}
\; \subset \;
\PP^{3}. 
$$

As shown in~\cite{Tsch}, 
the surface $X$ has a free divisor
class group $\Cl(X)$ of rank $7$ with a 
(quite) canonical basis 
$$F_{1}, F_{2}, F_{3}, F_{4}, F_{5}, F_{6}, \ell$$
where the $F_{i}$ are the exceptional curves of the 
resolution $X \to X'$ and $\ell$ is the strict transform
of a certain distinguished line on $X'$, 
see~\cite[Prop.~4.4]{Tsch}.

Moreover, Hassett and Tschinkel determined the total 
coordinate ring of $X$. As a ring, it is given by
$$
\mathcal{R}(X)
\; = \; 
\CC[T_{1}, \ldots, T_{10}] 
/
\bangle{T_{10}T_{7}^{3}T_{4}^{2}T_{5}+T_{9}^{2}T_{2}+T_{8}^{3}T_{1}^{2}T_{3}}.
$$

The $\Cl(X)$-grading of $\mathcal{R}(X)$ is defined 
by sending
$T_{i}$ to the vector in $\Cl(X)$ having as its 
coordinates with respect to the above canonical 
basis the $i$-th column of the following matrix:
\begin{equation}\label{eq:matrixdef}
\left[
\begin{array}{cccccccccc}
1 & 0 & 0 & 0 & 0 & 0 & 0 & 0 & 1 & 2 \\
0 & 1 & 0 & 0 & 0 & 0 & 0 & 1 & 1 & 3 \\
0 & 0 & 1 & 0 & 0 & 0 & 0 & 1 & 2 & 4 \\
0 & 0 & 0 & 1 & 0 & 0 & 0 & 2 & 3 & 4 \\
0 & 0 & 0 & 0 & 1 & 0 & 0 & 2 & 3 & 5 \\
0 & 0 & 0 & 0 & 0 & 1 & 0 & 2 & 3 & 6 \\
0 & 0 & 0 & 0 & 0 & 0 & 1 & 2 & 3 & 3 \\
\end{array}
\right]
\end{equation}
Using Remark \ref{primetest}, it is easy to see 
that the classes of the $T_{i}$ in $\mathcal{R}(X)$
are indeed pairwise nonassociated primes.

In order to apply our results, 
we have to describe $X$ 
as the variety arising from a bunched ring 
$(R,\mathfrak{F},\Phi)$, where 
$R := \mathcal{R}(X)$ with the grading 
lattice $K = \ZZ^{7}$.
Let $\mathfrak{F}$ consist of the
classes of the variables $T_{i}$ in $R$.

The associated projected cone $(E \topto{Q} K, \gamma)$
is given by $E := \ZZ^{10}$, the map $Q \colon E \to K$
determined by the matrix~\ref{eq:matrixdef}, and
$K = \ZZ^{7}$. 
Similar as above, we write 
$\tau(i_{1}, \dots, i_{s})$ for the cone in $K_{\QQ}$ 
generated by $Q(e_{i_{1}}), \ldots, Q(e_{i_{s}})$. 
Consider
$$
\begin{array}{ll}
\tau(1,2,3,4,7,8,9,10),
&
\tau(1,2,3,5,6,8,9,10),
\\
\tau(1,2,3,6,7,8,9,10),
&
\tau(1,2,4,5,7,8,9,10),
\\
\tau(2,4,5,6,7,8,9,10),
&
\tau(1,3,4,5,7,8,9,10).
\end{array}
$$ 
Using Proposition~\ref{hypersfrel}, one sees that these 
are the cones of an $\mathfrak{F}$-bunch $\Phi$.
Thus, Propositions~\ref{amplecone} and~\ref{candivdescr}
give us the semiample and the ample cone of $X$
and its canonical divisor class, 
compare~\cite[Prop.~4.5]{Tsch}:

\begin{proposition}
\label{hatschi}
The semiample cone $\SAmple(X)$ 
is the regular cone in $\ZZ^{7}$
having as its primitive generators the 
columns of the matrix:
$$ 
\left[
\begin{array}{ccccccc}
0 & 2 & 1 & 2 & 2 & 2 & 1 \\
1 & 3 & 1 & 3 & 3 & 3 & 2 \\ 
1 & 4 & 2 & 4 & 4 & 4 & 2 \\
2 & 4 & 3 & 4 & 5 & 6 & 4 \\
2 & 5 & 3 & 5 & 5 & 6 & 4 \\
2 & 6 & 3 & 6 & 6 & 6 & 4 \\
2 & 3 & 3 & 4 & 5 & 6 & 4 \\
\end{array}
\right]
$$
Moreover, $\SAmple(X)$ equals the moving cone, 
the ample cone is the relative interior of  
$\SAmple(X)$, and the canonical divisor class 
of $X$ is given by the vector
$-(2,3,4,4,5,6,3)$.
\endproof
\end{proposition}

Finally, we provide the fan describing 
the minimal toric ambient variety 
$Z := Z(R,\mathfrak{F},\Phi)$ 
of $X$ as discussed in Proposition~\ref{minimalambient}. 
Note that in the present example, 
$Z$ is not complete, moreover it is not 
even $A_2$-maximal.

\begin{proposition}
\label{hatschi2}
The fan corresponding to the minimal toric 
ambient variety $Z = Z(R,\mathfrak{F},\Phi)$ 
of $X$ has ten maximal cones;
one of them is the negative orthant\/
$\cone((-1,0,0), (0,-1,0), (0,0,-1))$, 
and the remaining nine are of dimension 
two:
$$ 
\begin{array}{lll}
\cone((2,3,3), (0,0,-1)), &
\cone((1,1,3), (0,-1,0)), & 
\cone((0,1,2),(-1,0,0)), 
\\
\cone((2,3,4),(2,3,3)), & 
\cone((2,3,5), (2,3,6)), & 
\cone((1,2,4), (2,3,6)),
\\ 
\cone((1,1,3),(2,3,6)), & 
\cone((2,3,4),(2,3,5)), & 
\cone((0,1,2),(1,2,4)).
\end{array}
$$
\end{proposition}

\proof
First we have to determine the covering collection of
the $\mathfrak{F}$-bunch $\Phi$. It consists of the following
ten faces of $\gamma = \cone(e_1, \ldots, e_{10})$:
$$
\begin{array}{lll}
\gamma(1,2,3,4,7,8,9,10), & \gamma(1,2,3,5,6,8,9,10), &
\gamma(1,2,3,6,7,8,9,10), 
\\
\gamma(1,2,4,5,7,8,9,10), & \gamma(2,4,5,6,7,8,9,10), &
\gamma(1,3,4,5,7,8,9,10), 
\\
\gamma(1,2,3,4,5,6,8,9), & \gamma(1,3,4,5,6,7,8,10), & 
\gamma(2,3,4,5,6,7,9,10), 
\\ 
\gamma(1,2,3,4,5,6,7), & & 
\end{array}
$$
where, as usual, we denote 
$\cone(e_{i_{1}}, \dots, e_{i_{r}})$ 
by $\gamma(i_{1}, \dots, i_{r})$.
Then, the maximal cones of the fan of $Z$ 
are the images $P(\gamma_0^*)$, 
where $\gamma_0 \in \cov(\Theta)$ and $P$ 
is the map dual to the inclusion of $M := \ker(Q)$ 
in $E$.
\endproof

\end{document}